\documentclass[runningheads]{llncs}
\usepackage[T1]{fontenc}
\usepackage{graphicx}
\usepackage{color}
\usepackage{hyperref} 

\usepackage{multirow}%
\usepackage{amsmath,amssymb}%
\usepackage{mathrsfs}%
\usepackage[title]{appendix}%
\usepackage{xcolor}%
\usepackage{textcomp}%
\usepackage{manyfoot}%
\usepackage{booktabs}%
\usepackage{algorithm}%
\usepackage{algorithmicx}%
\usepackage{algpseudocode}%
\usepackage{listings}%

\usepackage[all]{xy}

\usepackage{float}
\usepackage{tikz}
\usetikzlibrary{
  knots,
  hobby,
  decorations,
  decorations.pathreplacing,
  shapes.geometric,
  calc,
  intersections,
  graphs
}

\tikzset{
  every path/.style={
    ultra thick,
    black,
  },
  show curve controls/.style={
    postaction=decorate,
    decoration={show path construction,
      curveto code={
        \draw [blue, dashed]
        (\tikzinputsegmentfirst) -- (\tikzinputsegmentsupporta)
        node [at end, draw, solid, red, inner sep=2pt]{};
        \draw [blue, dashed]
        (\tikzinputsegmentsupportb) -- (\tikzinputsegmentlast)
        node [at start, draw, solid, red, inner sep=2pt]{}
        node [at end, fill, blue, ellipse, inner sep=2pt]{}
        ;
      }
    }
  },
  show curve endpoints/.style={
    postaction=decorate,
    decoration={show path construction,
      curveto code={
        \node [fill, blue, ellipse, inner sep=2pt] at (\tikzinputsegmentlast) {}
        ;
      }
    }
  }
}

\usepackage{caption} 
\usepackage{subcaption}
\usepackage{comment}
\usepackage{enumitem}

\begin{document}
\title{The untangling number of 3-periodic tangles}
%
%
\author{Toky Andriamanalina\inst{1} \and
Sonia Mahmoudi\inst{2,3}\orcidID{0000-0003-2095-4066} \and
Myfanwy E. Evans\inst{1,*}\orcidID{0000-0002-0161-6523}}
\authorrunning{T. Andriamanalina et al.}
%
\institute{
Institute for Mathematics, University of Potsdam, Karl-Liebknecht-Str. 24-25, Potsdam Golm 14476, Germany \and
Advanced Institute for Materials Research, Tohoku University, 2-1-1 Katahira, Aoba-ku, Sendai 980-8577, Japan \and
RIKEN iTHEMS, 2-1 Hirosawa, Wako, 351-0198, Japan \\
\inst{*}\email{evans@uni-potsdam.de}
}
\maketitle          
\begin{abstract}
The entanglement of curves within a 3-periodic box provides a model for complicated space-filling entangled structures occurring in biological, chemical and physical systems. Quantifying the complexity of the entanglement within these models enhances the characterisation of these structures. In this paper, we introduce a new measure of entanglement complexity through the \emph{untangling number}, reminiscent of the unknotting number in knot theory. The untangling number quantifies the minimum distance between a given 3-periodic structure and its least tangled version, called \emph{ground state}, through a sequence of operations in a diagrammatic representation of the structure. For entanglements that consist of only infinite open curves, we show that the generic ground states are crystallographic rod packings, well known in structural chemistry.


\end{abstract}
\section{Introduction}\label{sec:1}

Tangling is a fundamental feature of three-dimensional structures composed of multiple long, curved filaments in space. Quantifying tangling is key to describing a multitude of entangled biological, chemical and physical systems, such as the mesoscale structure of mammalian skin cells (corneocytes) \cite{evans2011,PhysRevLett.112.038102}, polymer melts \cite{MICHELETTI20111,tzoumanekas_theodorou}, liquid crystals \cite{Ertman2023PeriodicLC}, DNA origami crystals \cite{dna_tensegrity_triangle,topology_based_lattice_engineering}, and especially volumetric textiles and mechanical metamaterials \cite{volumetric_bravais_weaves,portela_volumetric_textile,PESCIALLI2025114974}.

\begin{figure}[hbtp]
    \centering
    \begin{subfigure}[b]{0.33\textwidth}
    \centering
        \includegraphics[width=\textwidth]{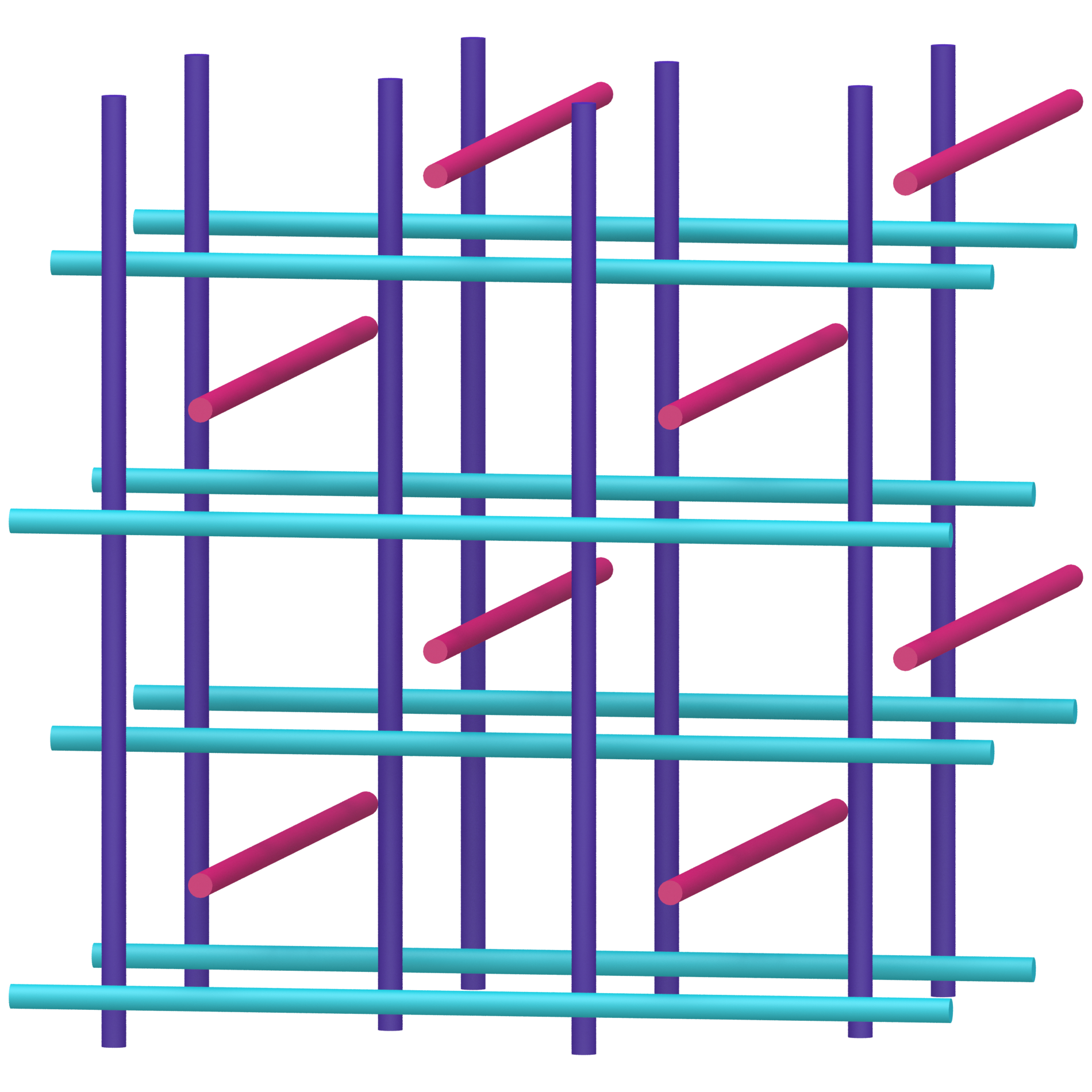}
        \caption{}
        \label{fig:pi_plus_extended}
    \end{subfigure}
    \begin{subfigure}[b]{0.33\textwidth}
    \centering
        \includegraphics[width=\textwidth]{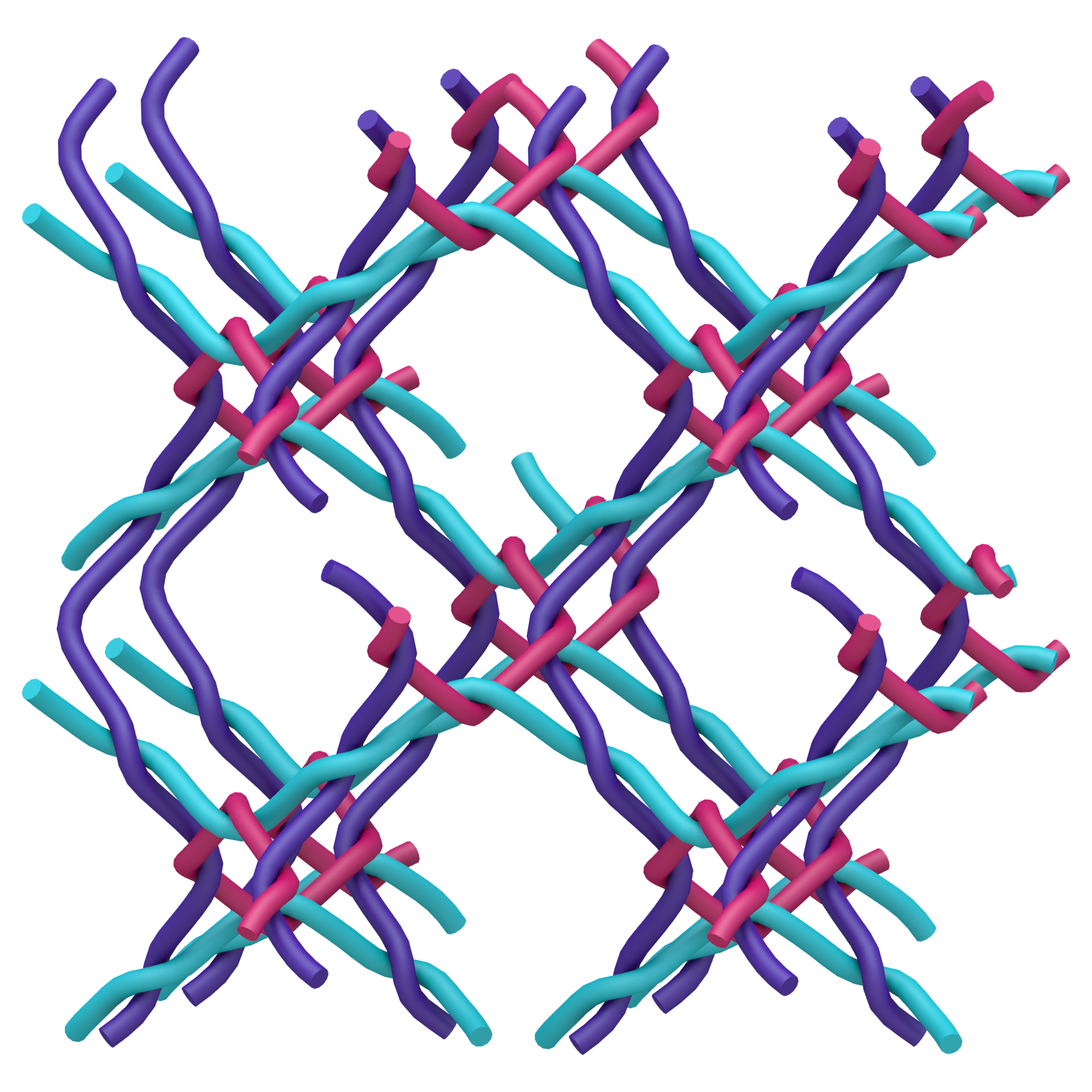}
        \caption{}
        \label{fig:1p6on2_to_the_6_extended}
    \end{subfigure}
    \caption{Two different 3-periodic arrangements of straight or curved filaments. (a) A 3-periodic arrangement of straight filaments, known in structural chemistry as the $\Pi^{+}$ rod packing \cite{OKeeffe:au0217}. (b) A 3-periodic arrangement of curved filaments that tangle around each other, where the invariant axes of the filaments is also the $\Pi^{+}$ shown to the left. This structure can be intuitively considered as a more tangled version of the $\Pi^{+}$ rod packing.}
    \label{fig:various_tangles}
\end{figure}

Modelling systems of entangled filaments using periodic boundary conditions (that is, using periodic unit cells) is a useful way to reduce the complexity of extended three-dimensional structures. This is a natural restriction for crystalline structures, but it remains useful for understanding less symmetric systems. Geometric and topological characterisations of 1-periodic, 2-periodic \cite{grishanov2007,grishanov2009_part_1,grishanov2009,diamantis2023equivalencedoublyperiodictangles,fukuda2023,fukuda2024constructionweavingpolycatenanemotifs,mahmoudi2024classificationperiodicweavesuniversal}, and to a limited extent, 3-periodic entanglements of filaments \cite{periodic_ent_II,purcell2024,purcell2025,PANAGIOTOU2015533,Barkataki_2024}, that we call \emph{3-periodic tangles}, have been investigated in different contexts. Key questions arise as to how tangled a structure is, and which configuration is considered to be the least entangled state. For example, Fig.~\ref{fig:various_tangles} shows two distinct 3-periodic arrangements of straight or curved filaments. Intuitively, the configuration in Fig.~\ref{fig:pi_plus_extended} is less tangled than the one in Fig.~\ref{fig:1p6on2_to_the_6_extended}, and it is this idea that we want to formulate in a rigorous setting.

The characterisation of entanglement in 3-periodic tangles is comparable to the mathematical descriptions of knots and links \cite{Adams.book}. A link is an embedding (a geometric realisation) of finitely many loops in space. A knot is a link with a single component. Knots and links are classified up to ambient isotopy, which means that two embeddings represent the same knot or link if one embedding can be transformed into the other by deforming it in space without cutting and gluing. A diagram of a knot is a projection of the knot onto a plane to which one adds the crossing information (see Fig.~\ref{fig:trefoil_a}). The \emph{crossing number} of a knot is the minimum number of crossings in any diagram of the knot. The simplest knot is a circle, which has crossing number 0 and is called the \emph{unknot}.
Any knot can be transformed into the unknot by applying a finite sequence of transformations called \emph{crossing changes}. These consist of changing the parity of a crossing, where a strand passing under another is changed to passing over, similar to allowing the knot to pass through itself (see Fig.~\ref{fig:trefoil_b}). The least number of crossing changes among all possible diagrams of the knot needed to obtain the unknot is called the \emph{unknotting number} \cite{Murasugi1996chap4}. The crossing number and the unknotting number are invariant under isotopy and are used to measure knot complexity and describe knotting in knotted systems. In particular, the conceptual machinery underlying the unknotting number has been central to understanding the action of type II topoisomerases on DNA molecules \cite{LIU1980697}.

\begin{figure}[hbtp]
    \centering
    \begin{subfigure}[b]{0.25\textwidth}
        \centering
        \resizebox{\linewidth}{!}{
            \begin{tikzpicture}
            \path[spath/save=trefoil]
            (0,2) .. controls +(2.2,0) and +(120:-2.2) ..
            (210:2) .. controls +(120:2.2) and +(60:2.2) ..
            (-30:2) .. controls +(60:-2.2) and +(-2.2,0) .. (0,2);
            \tikzset{
              every trefoil component/.style={draw},
              spath/knot={trefoil}{15pt}{1,3,5},
            }
            \end{tikzpicture}
        }
        \caption{Trefoil knot.}
        \label{fig:trefoil_a}
    \end{subfigure}
    \hspace{0.5cm}
    \begin{subfigure}[b]{0.5\textwidth}
        \centering
        \begin{minipage}{0.5\linewidth}
            \centering
            \resizebox{\linewidth}{!}{
                \begin{tikzpicture}
            \path[spath/save=trefoil]
            (0,2) .. controls +(2.2,0) and +(120:-2.2) ..
            (210:2) .. controls +(120:2.2) and +(60:2.2) ..
            (-30:2) .. controls +(60:-2.2) and +(-2.2,0) .. (0,2);
            \tikzset{
              every trefoil component/.style={draw},
              spath/knot={trefoil}{15pt}{3,4,5}, 
            }
            \end{tikzpicture}
            }
        \end{minipage}%
        \begin{minipage}{0.5\linewidth}
            \centering
            \begin{tikzpicture}
                \draw[line width=0.775pt] (0,0) circle (0.85); 
            \end{tikzpicture}
            \vspace{0.775cm} 
        \end{minipage}
        \caption{Unknotting the trefoil knot.}
        \label{fig:trefoil_b}
    \end{subfigure}
    \caption{(a) A crossing diagram of the trefoil knot. (b) A crossing change, which consists of changing the parity of a crossing. It transforms the trefoil knot into the unknot.}
    \label{fig:trefoil}
\end{figure}

\begin{figure}[hbtp]
    \centering
    \begin{subfigure}[b]{0.95\textwidth}
        \centering
        \includegraphics[width=0.6625\textwidth]{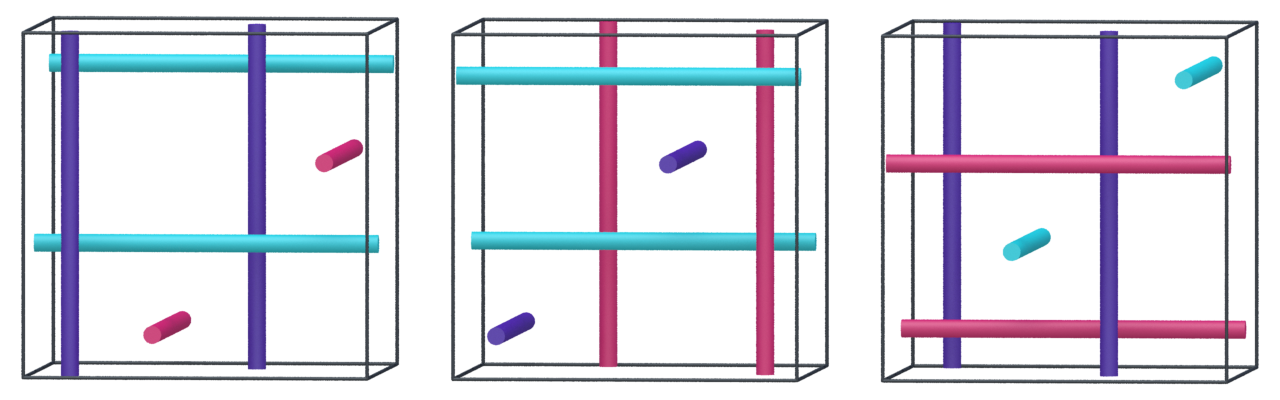}
        \caption{}
        \label{fig:pi_plus_uc_xyz}
    \end{subfigure}
    \vskip\baselineskip
    \begin{subfigure}[b]{0.95\textwidth}
        \centering
        \includegraphics[width=0.2\textwidth]{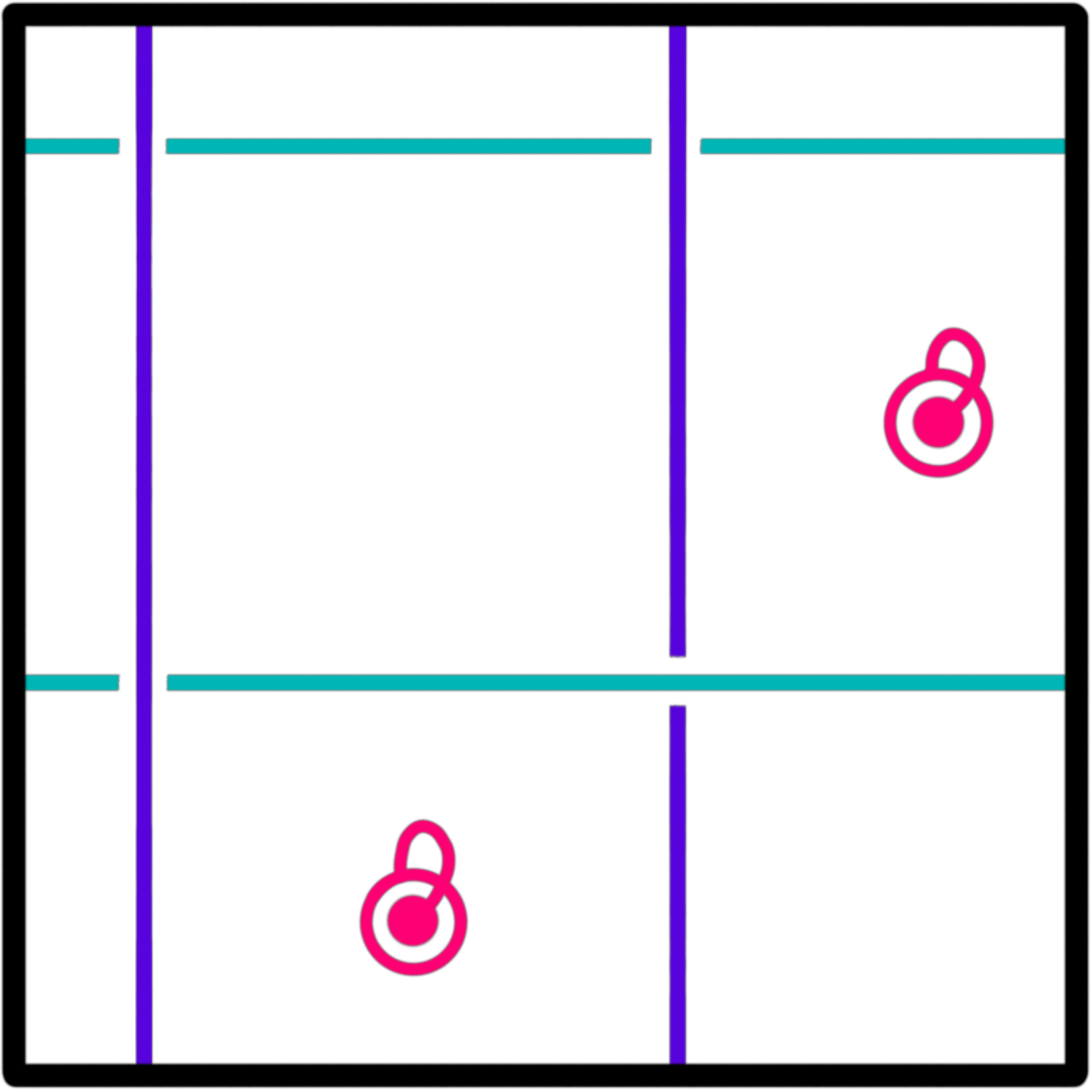}
        \hspace{0.0825cm}
        \includegraphics[width=0.2\textwidth]{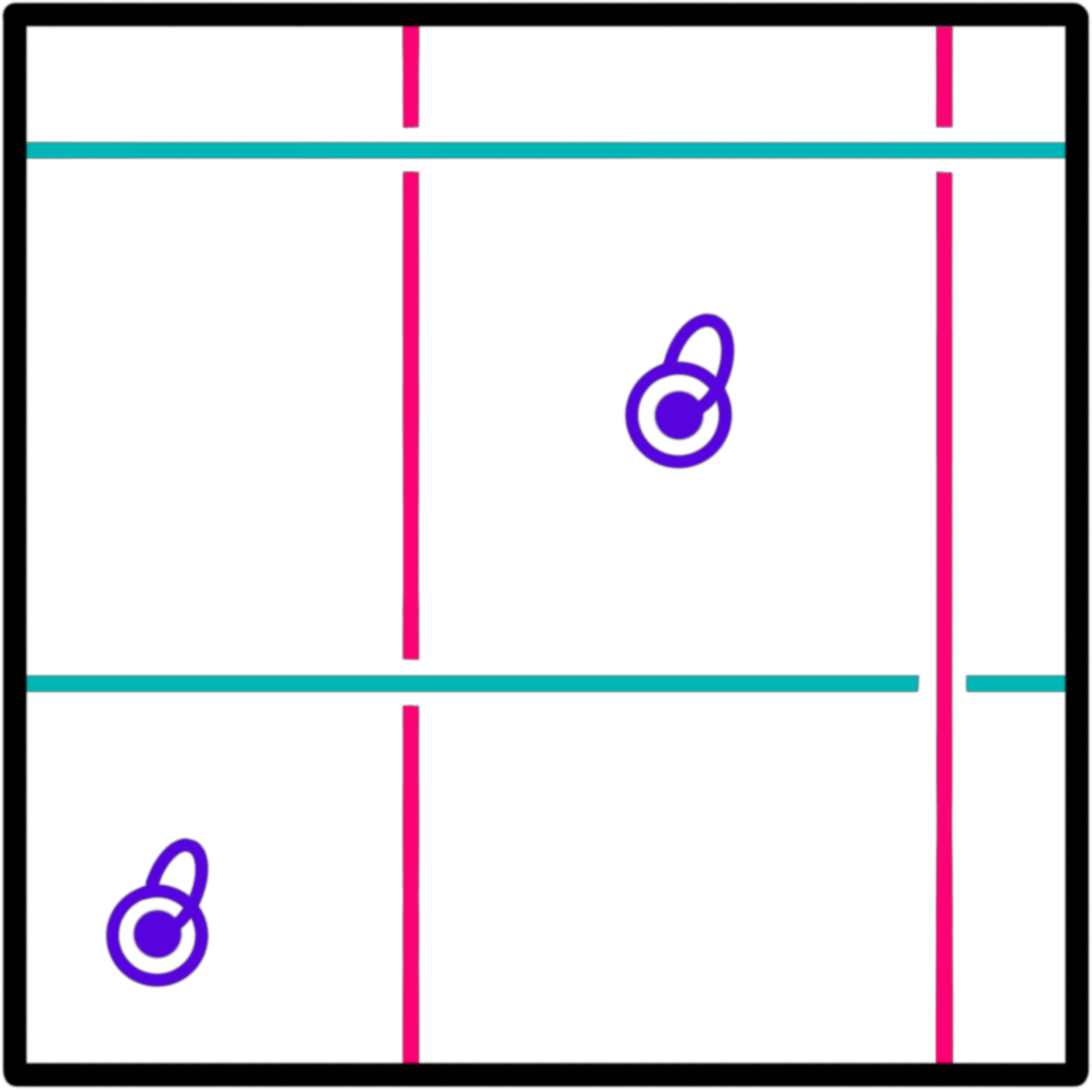}
        \hspace{0.0825cm}
        \includegraphics[width=0.2\textwidth]{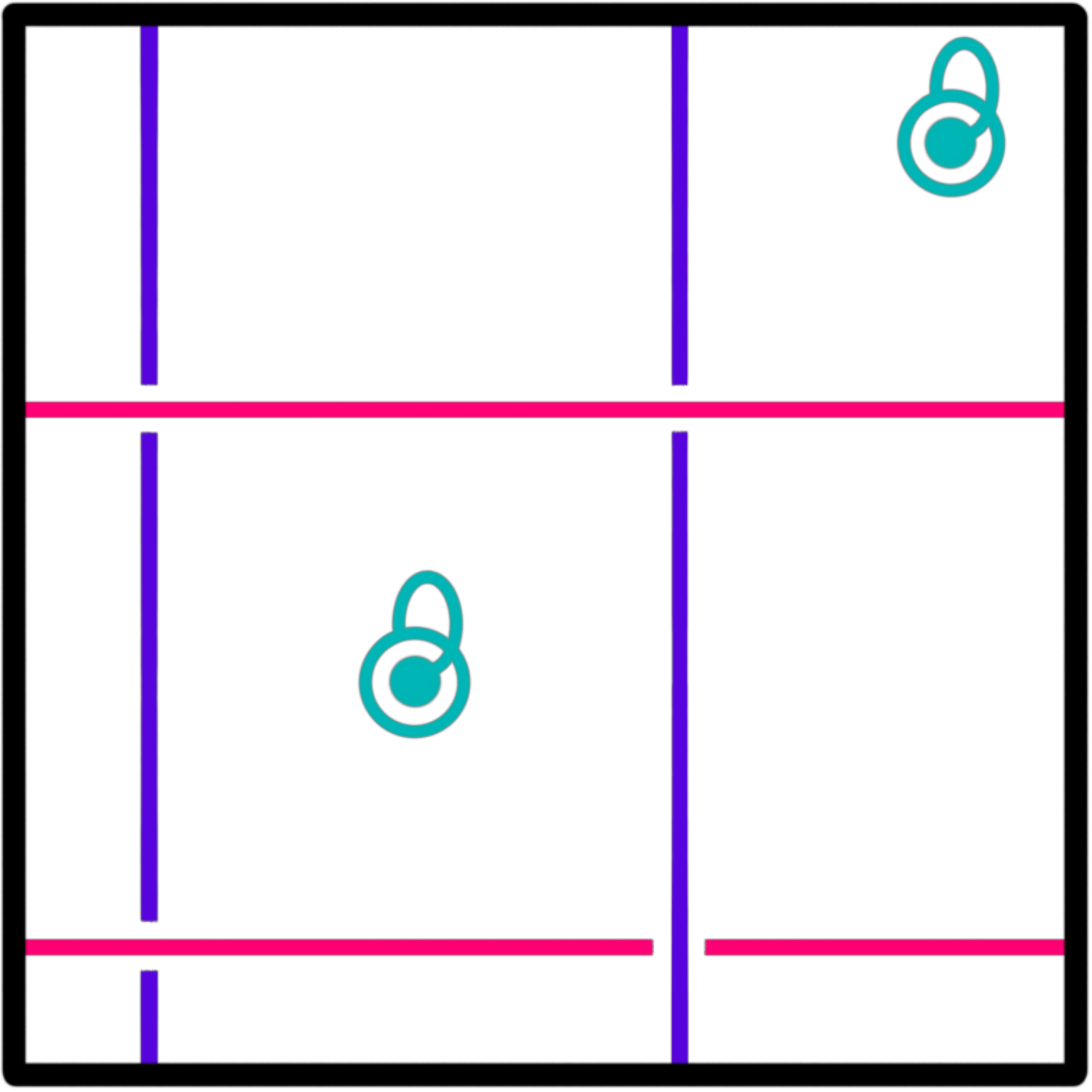}
        \caption{}
        \label{fig:pi_plus_tridia}
    \end{subfigure}
    \caption{Transformation of a unit cell of a 3-periodic tangle into a tridiagram. (a) One unit cell of the $\Pi^{+}$ rod packing, rotated in space so as to be viewed from the front, the top and the right sides of the structure. (b) A tridiagram obtained from (a), where the three different orientations of the unit cell are each projected to the plane.}
    \label{fig:pi_plus_uc_xyz_and_tridia}
\end{figure}

Knot-theoretic diagrammatic representations of 3-periodic tangles have been recently defined (within a unit cell) and rigorously investigated from a topological perspective in \cite{ANDRIAMANALINA2025109346}. A unit cell is a domain delimited by a parallelepiped with identified opposite faces that generates the entire 3-periodic structure under translations. A diagram is obtained from a projection of a unit cell along one of the vectors delimiting the parallelepiped. Three projections along three non-coplanar vectors constitute a \emph{tridiagram}. A unit cell and a tridiagram of the $\Pi^{+}$ rod packing, previously shown in Fig.~\ref{fig:pi_plus_extended}, is displayed in Fig.~\ref{fig:pi_plus_uc_xyz_and_tridia}. It was proved in \cite{ANDRIAMANALINA2025109346} that tridiagrams fully encode ambient isotopy transformations of 3-periodic tangles, provided that the periodicity of the structures is preserved. These tridiagrams allow the definition of the crossing number for 3-periodic tangles, which is a basic complexity ordering of such structures, and thus, they open the door to further investigations of entanglement in 3-periodic systems, including the definition of new measures of entanglement complexity.

In this paper, we define and characterise the least tangled 3-periodic tangles, that we call \emph{ground states}, as well as the \emph{untangling number} of 3-periodic tangles, which quantifies the shortest path to transform a given structure into its least tangled version. The untangling number resembles the idea of the unknotting number of a knot, and can be used as a measure of complexity of a periodic entanglement of filaments. Note that the notions of least tangled embeddings and the untangling number of 3-periodic networks are defined in another paper \cite{andriamanalina2026measuringentanglementcomplexity3periodic}, in which the techniques employed share some similarities with those described here. This paper is organised as follows. In Section \ref{sec:reminder_diagrams}, we recall the definitions and results that provide the foundation for the tridiagrams of a 3-periodic tangle \cite{ANDRIAMANALINA2025109346}. In Section \ref{sec:ground_states}, we define and characterise the least tangled structures, which we call \emph{ground states}, by analogy with the unknot. In Section \ref{sec:untangling_number}, we define the \emph{untangling number} of 3-periodic tangles with respect to a unit cell, as well as the \emph{minimum untangling number}, which we illustrate with examples. In Section \ref{sec:computability}, we discuss the computability of the untangling number. Finally, Section \ref{sec:conclusion} is devoted to the conclusion of this work.

\section{Framework for 3-periodic tangles}\label{sec:reminder_diagrams}

In this section, we present the basic mathematical framework for describing 3-periodic tangles. A rigorous treatment is given in \cite{ANDRIAMANALINA2025109346}; here, we summarise the most important notions for defining the untangling number.

\subsection{3-periodic tangles and their diagrams}
\begin{figure}[hbtp]
    \centering
    \begin{subfigure}[b]{0.25\textwidth}
    \centering
        \includegraphics[width=\textwidth]{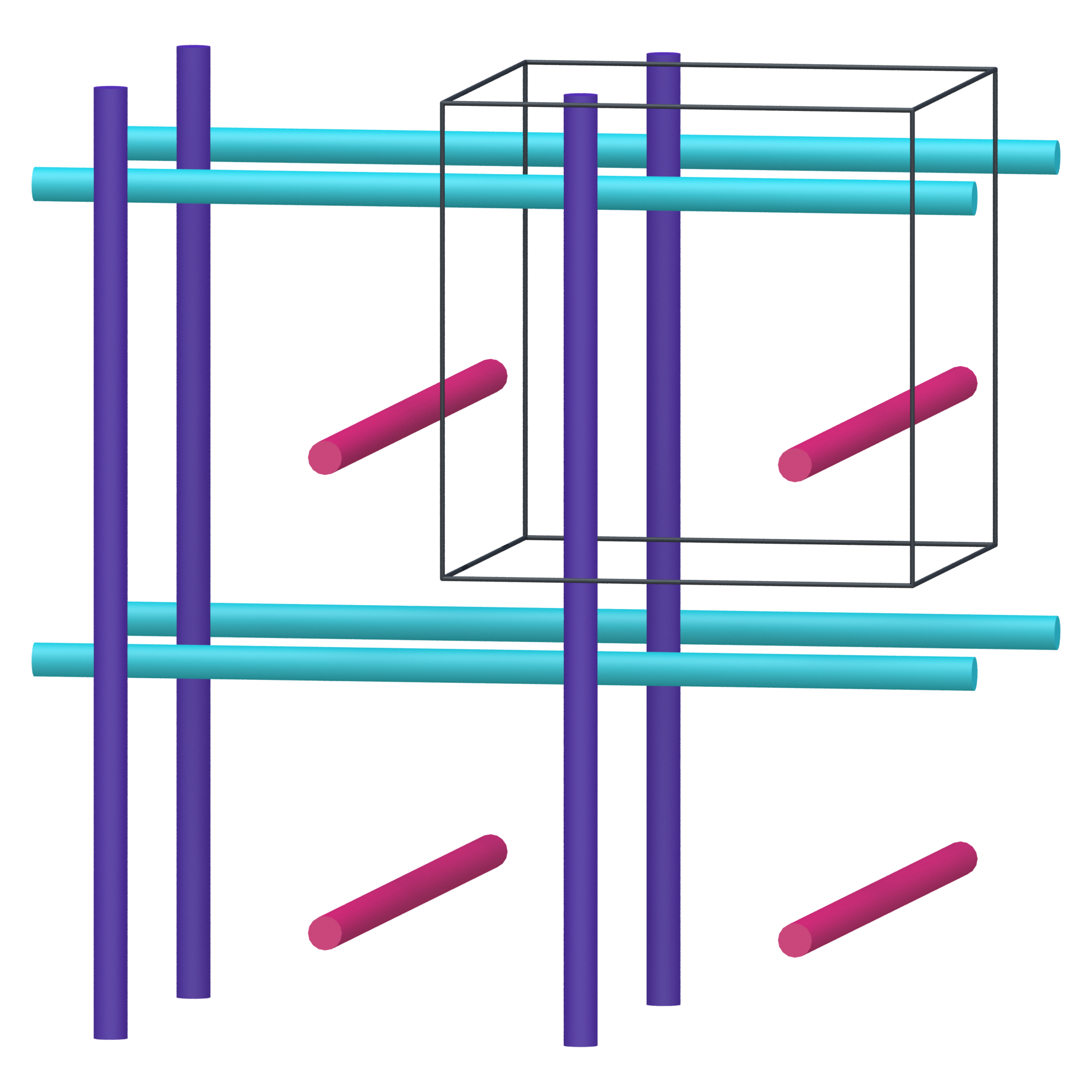}
        \caption{$\Pi^{\star}$}
        \label{fig:pi_star_extended_axis_100}
    \end{subfigure}
    \begin{subfigure}[b]{0.25\textwidth}
    \centering
        \includegraphics[width=\textwidth]{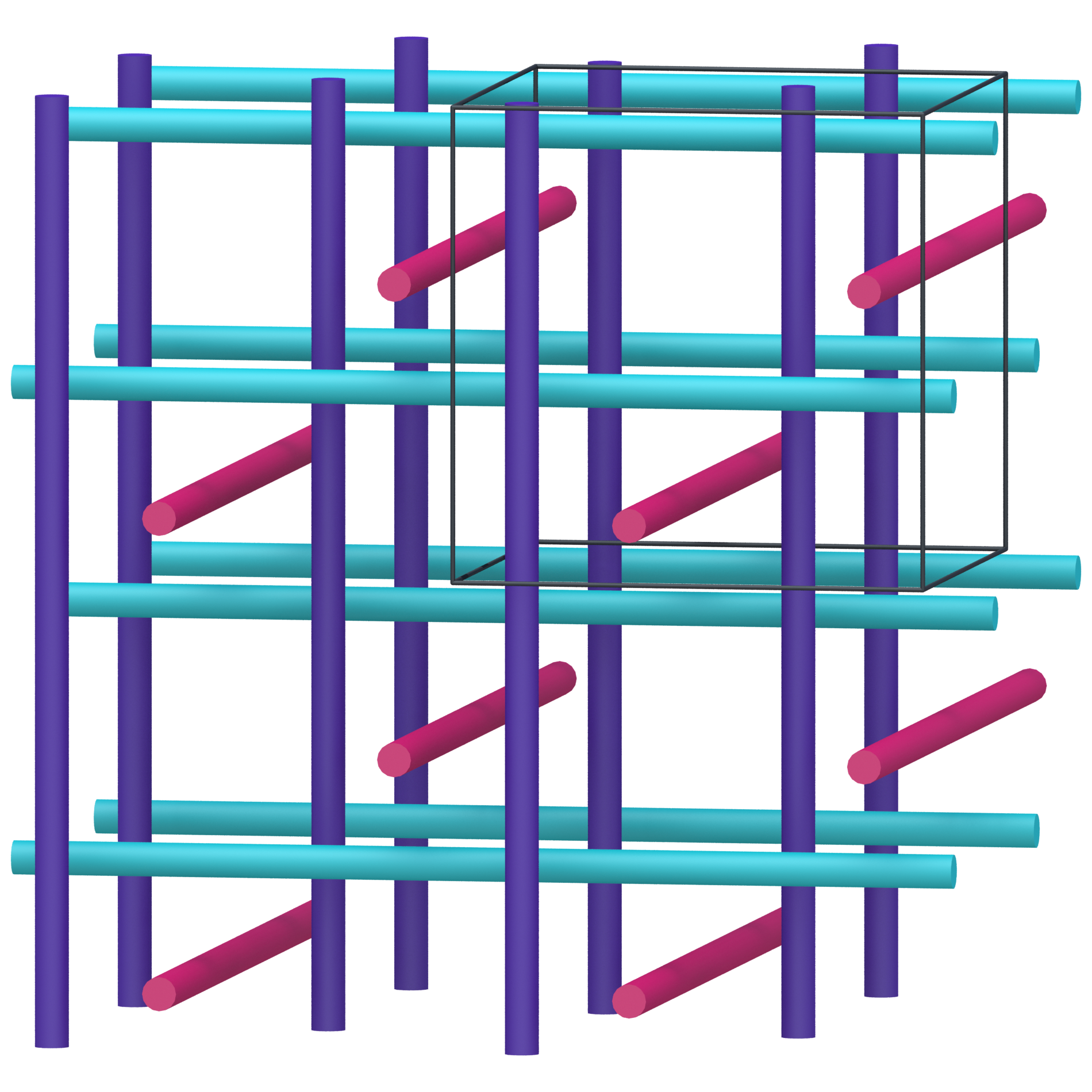}
        \caption{$\Pi^{+}$}
        \label{fig:pi_plus_extended_axis_100}
    \end{subfigure}
    \begin{subfigure}[b]{0.3\textwidth}
    \centering
        \includegraphics[width=\textwidth]{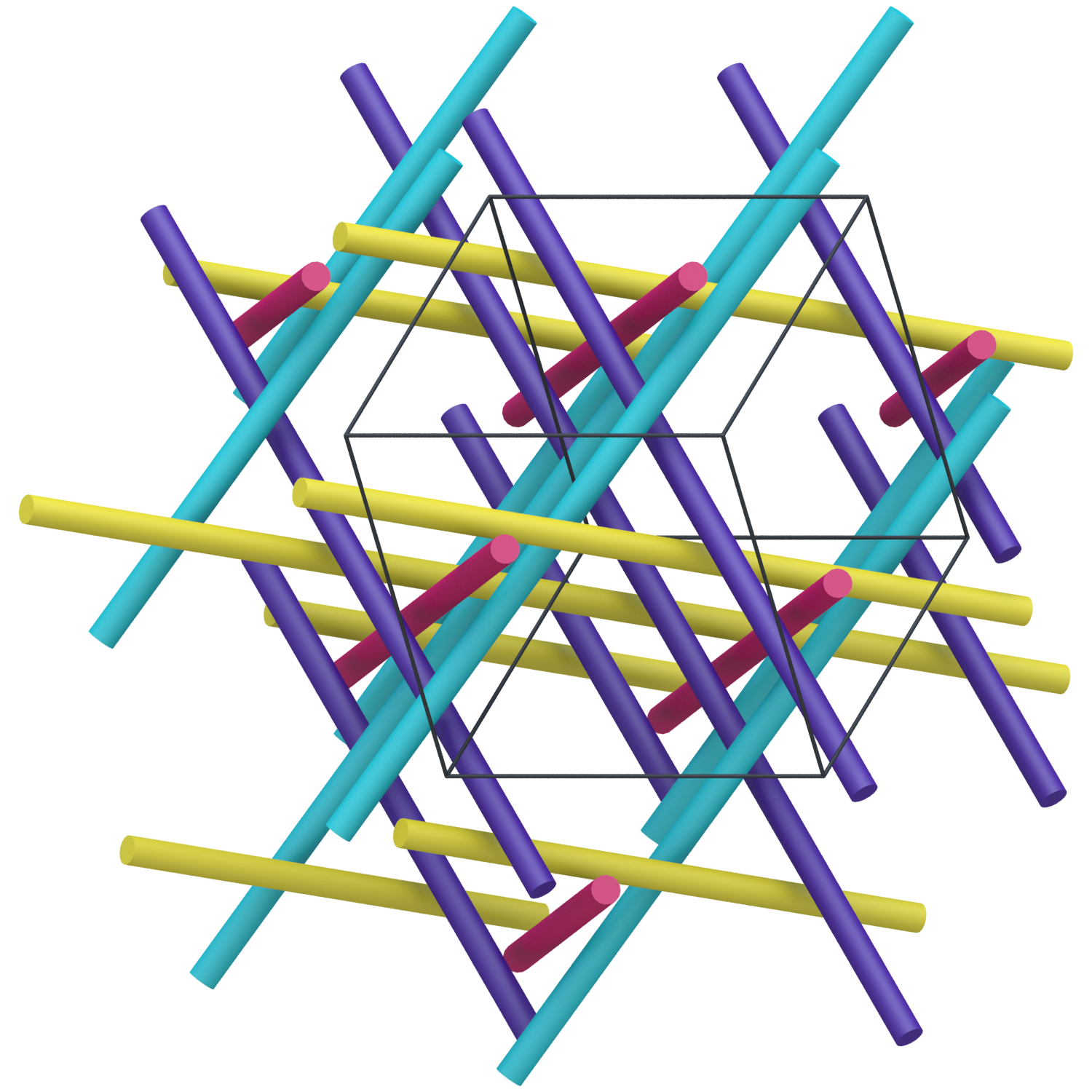}
        \caption{$\Sigma^{+}$}
        \label{fig:sigma_plus_extended_axis_111}
    \end{subfigure}
    \vskip\baselineskip
    \begin{subfigure}[b]{0.3\textwidth}
    \centering
        \includegraphics[width=\textwidth]{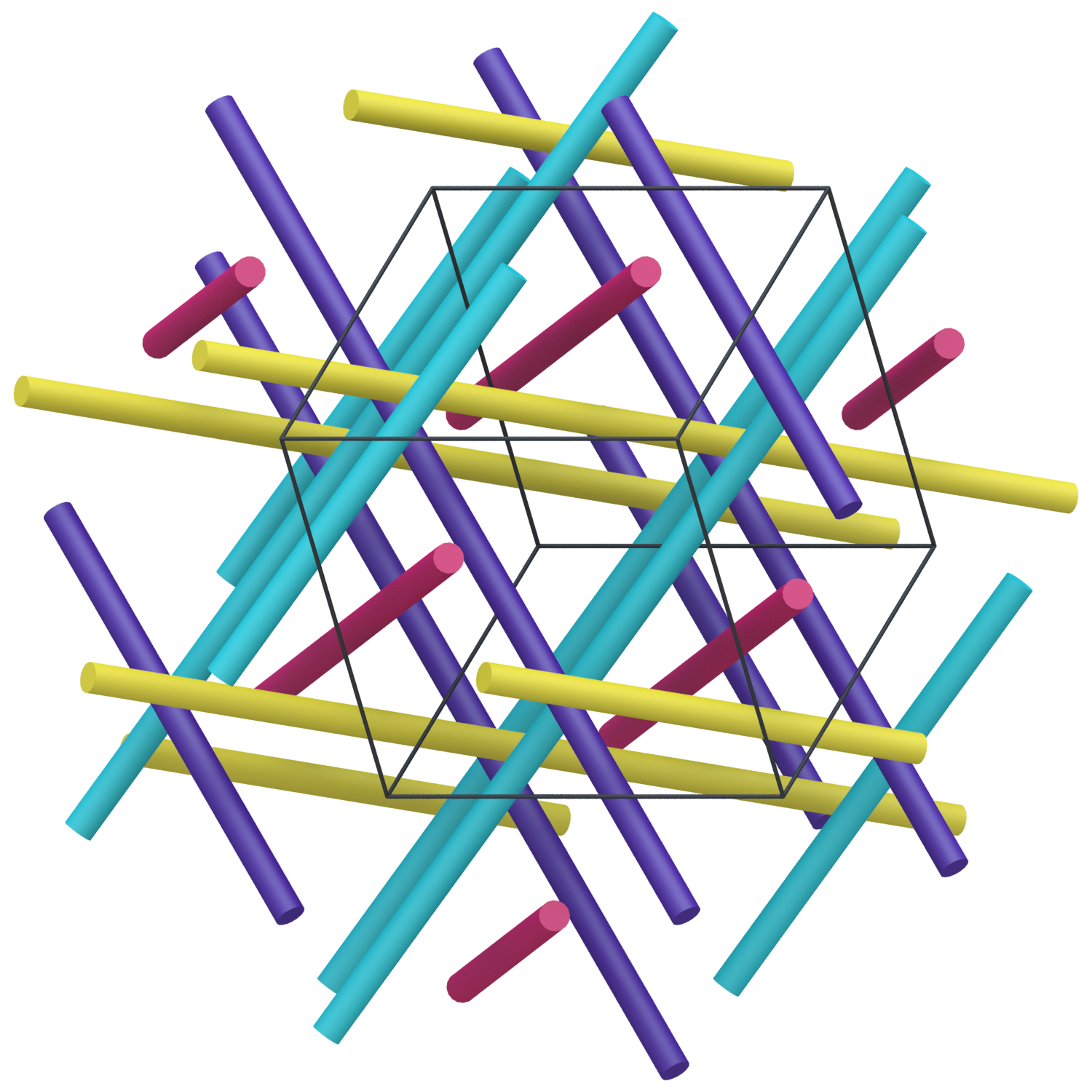}
        \caption{$\Omega^{+}$}
        \label{fig:omega_plus_extended_axis_111}
    \end{subfigure}
    \begin{subfigure}[b]{0.3\textwidth}
    \centering
        \includegraphics[width=\textwidth]{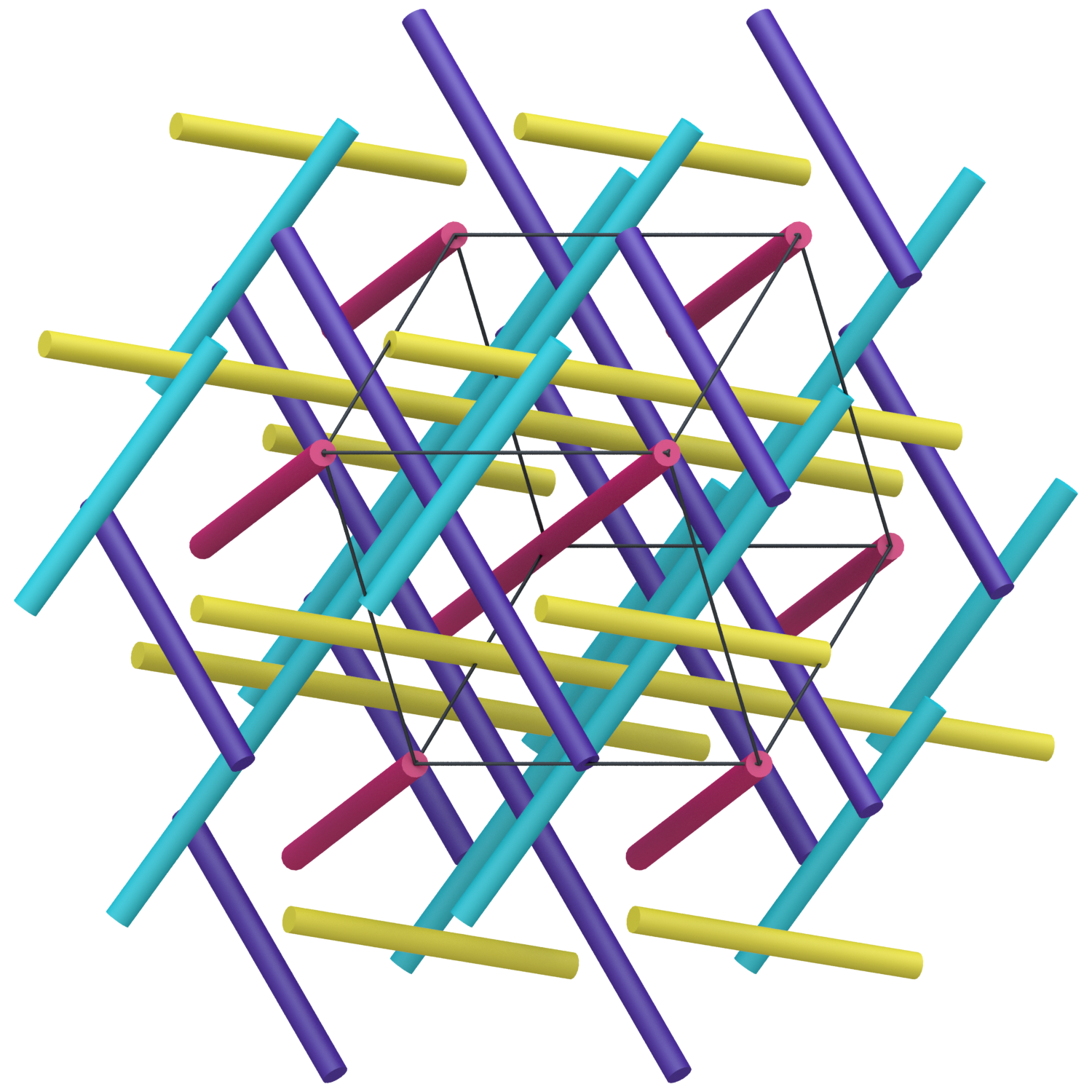}
        \caption{$\Gamma$}
        \label{fig:gamma_extended_axis_111}
    \end{subfigure}
    \begin{subfigure}[b]{0.3\textwidth}
    \centering
        \includegraphics[width=\textwidth]{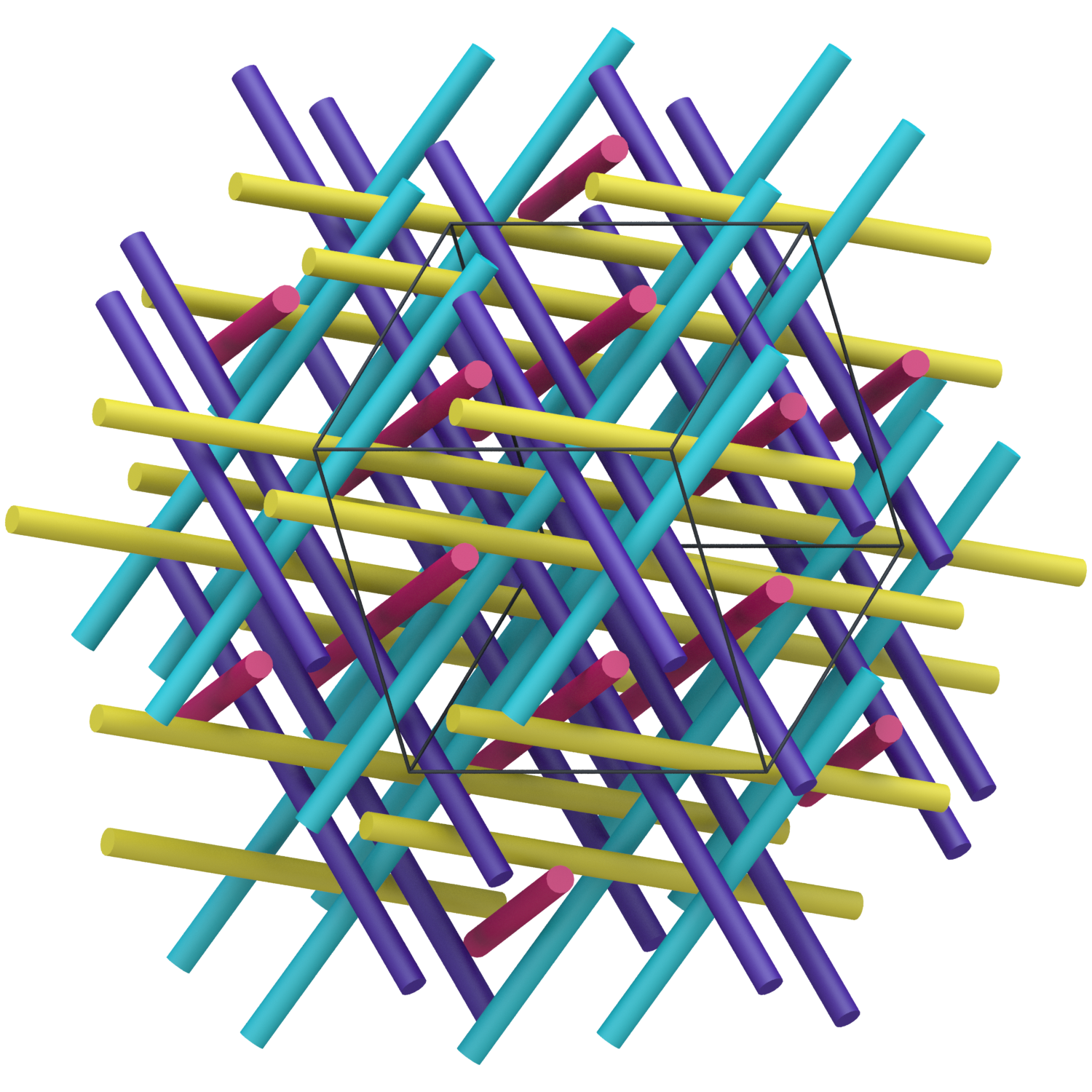}
        \caption{$\Sigma^{\star}$}
        \label{fig:sigma_star_extended_axis_111}
    \end{subfigure}
    \caption{The six invariant cubic rod packings, as described in \cite{OKeeffe:au0217}. They are the rod packings with cubic crystallographic symmetry, where rods are placed at the invariant axes of the space groups. They constitute simple examples of 3-periodic tangles, where the filaments are straight. The rods are coloured by those that are parallel; the colours have no meaning beyond visual clarification.}
    \label{fig:six_rod_packings}
\end{figure}

A \emph{3-periodic tangle} $K$ is a disjoint collection of simple curves that is invariant under translations along three vectors generating the space $\mathbb{R}^3$. Simple examples of 3-periodic tangles, where all of the filaments are straight, are the six invariant rod packings with cubic crystallographic symmetry $\Pi^{\star}$, $\Pi^{+}$, $\Sigma^{+}$, $\Omega^{+}$, $\Gamma$, and $\Sigma^{\star}$ studied in \cite{OKeeffe:au0217} and shown in Fig.~\ref{fig:six_rod_packings}.

\begin{figure}[hbtp]
    \centering
    \begin{subfigure}[b]{0.95\textwidth}
        \centering
        \includegraphics[width=0.6625\textwidth]{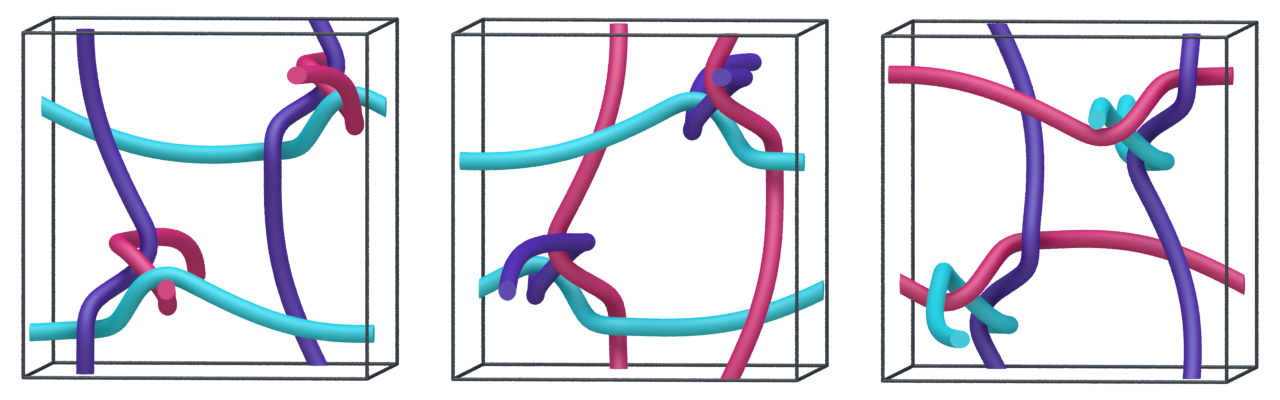}
        \caption{}
        \label{fig:pi_plus_u2_uc_xyz}
    \end{subfigure}
    \vskip\baselineskip
    \begin{subfigure}[b]{0.95\textwidth}
        \centering
        \includegraphics[width=0.2\textwidth]{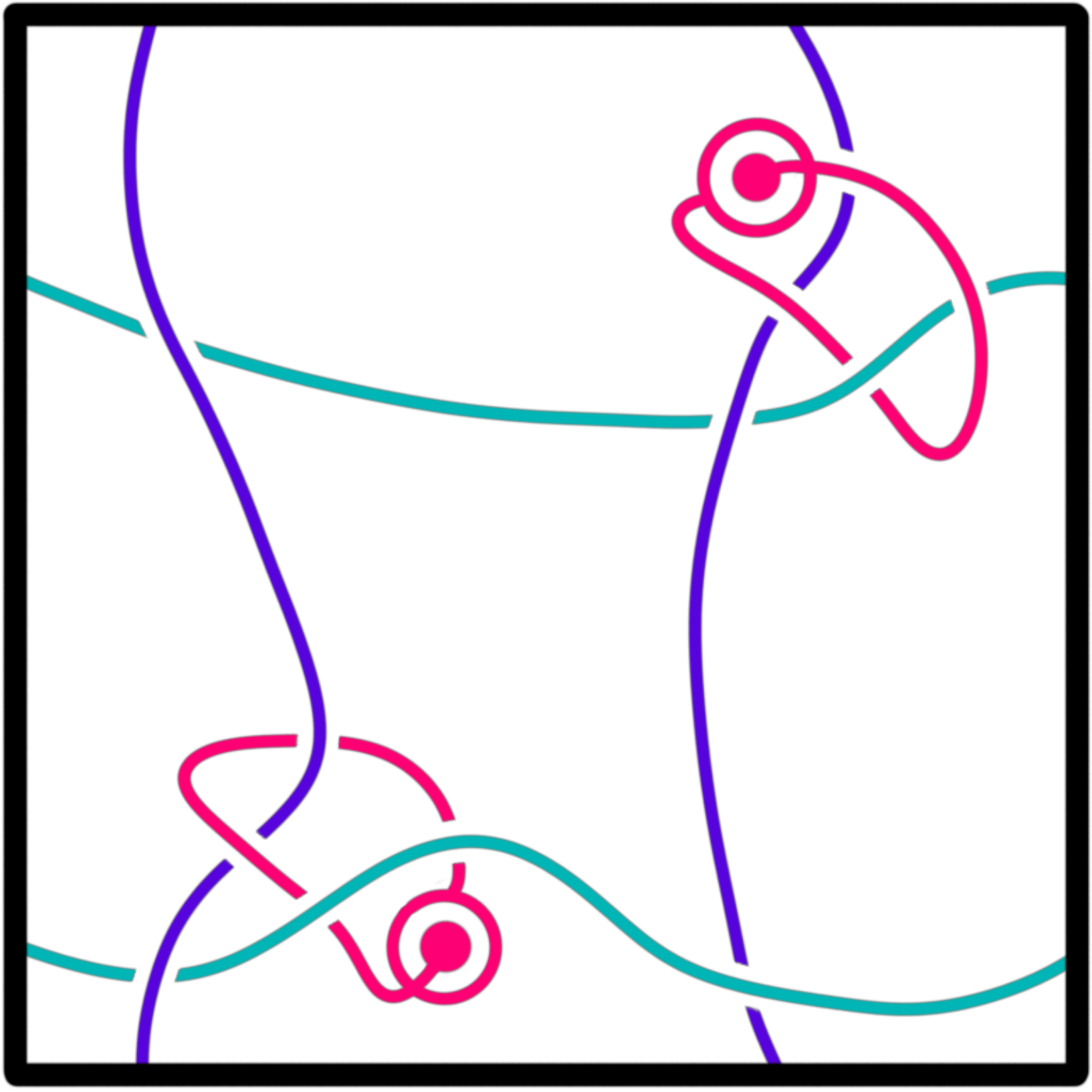}
        \hspace{0.0825cm}
        \includegraphics[width=0.2\textwidth]{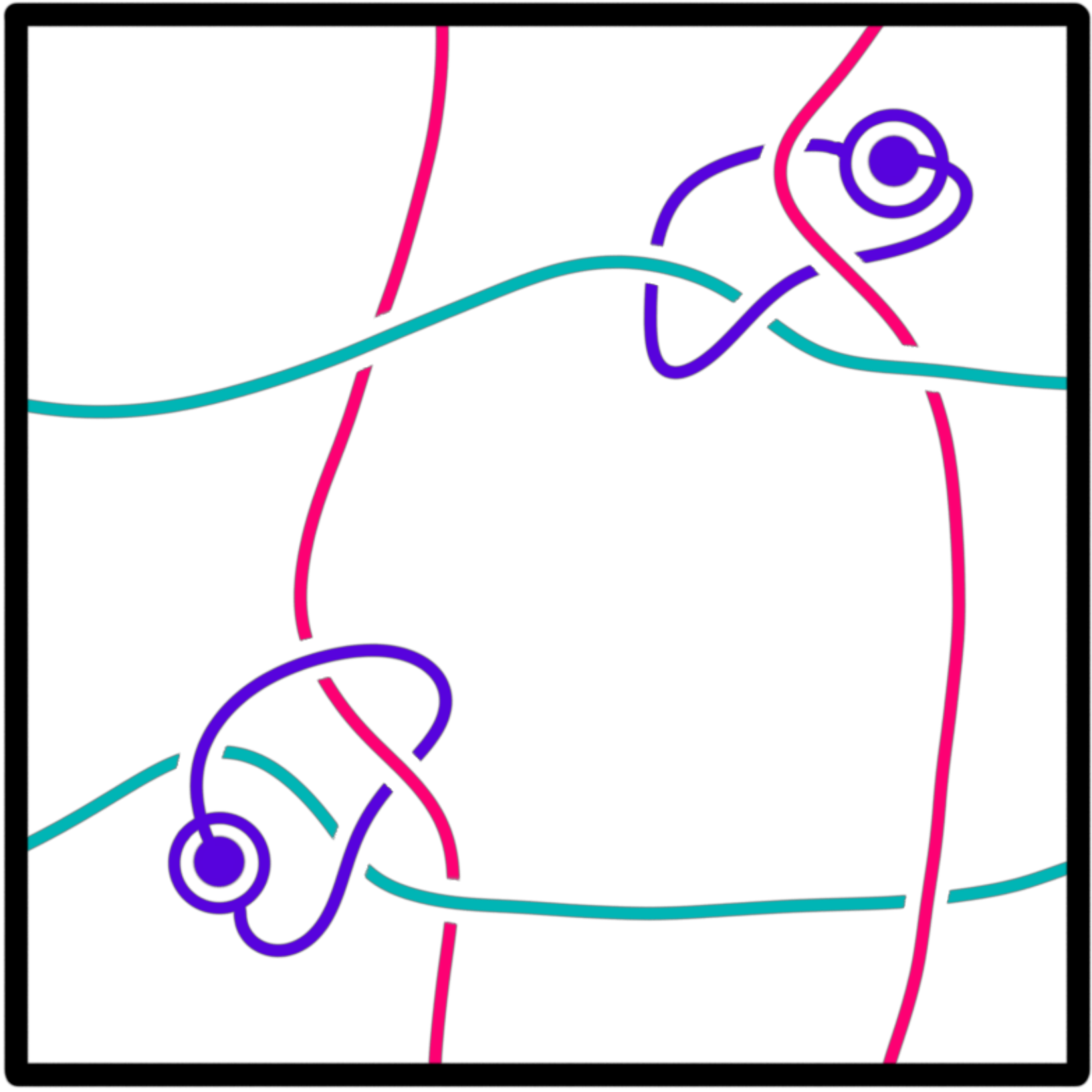}
        \hspace{0.0825cm}
        \includegraphics[width=0.2\textwidth]{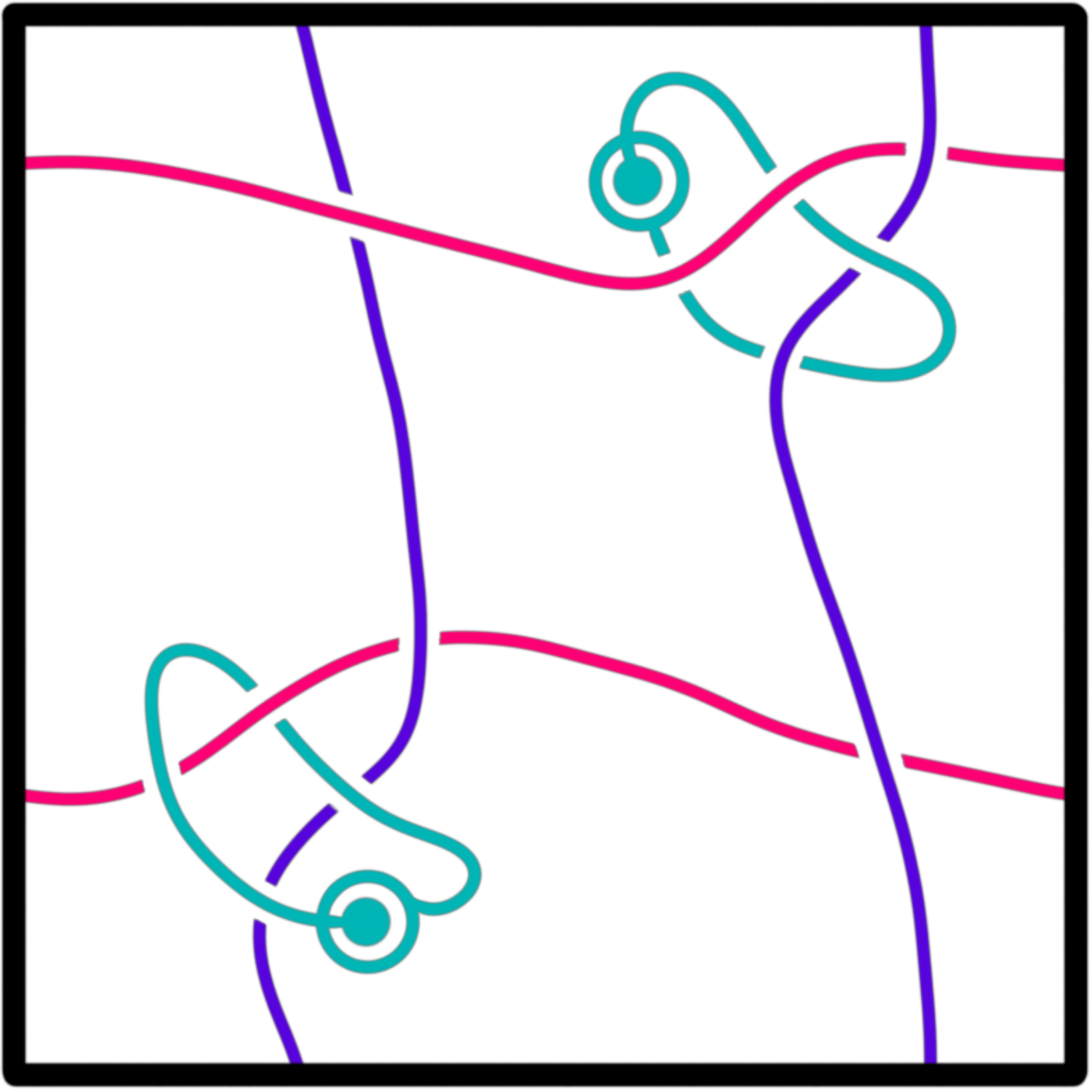}
        \caption{}
        \label{fig:pi_plus_u2_tridia}
    \end{subfigure}
    \caption{A unit cell and a tridiagram of a 3-periodic tangle. (a) One unit cell of a 3-periodic tangle viewed along three directions: the front, the top and the right sides of the cube. (b) A set of three diagrams that result from projections of the unit cell shown in (a). An intersection of a curve with the front face is encoded using a thick dot and an intersection with the back face is encoded with an open circle. Such a set of three diagrams is called a tridiagram.}
    \label{fig:pi_plus_u2_uc_xyz_and_tridia}
\end{figure}

Due to the periodicity of a given 3-periodic tangle, there exists a domain, delimited by a parallelepiped with identified opposite faces (a 3-torus), that we call a \emph{unit cell} and that generates the entire 3-periodic tangle. An example of a unit cell is shown in Fig.~\ref{fig:pi_plus_u2_uc_xyz}. A parallelepiped delimiting a unit cell can always be rectified to a cube, which we do for simplicity in this paper.

A \emph{unit-diagram} of a 3-periodic tangle, or simply \emph{diagram}, is obtained from the projection of a unit cell onto a square with identified opposite edges (a 2-torus). It is always assumed that one can deform the curves slightly so that each point of the projection has at most two preimages. 
To each double point (a point with two preimages), one assigns a \emph{crossing type} indicating which strand passes over the other. In addition to the crossings, other symbols encode the intersections of a curve with the identified front and back faces of the cube: a thick dot represents an intersection with the front face, and an open circle represents an intersection with the back face. Fig.~\ref{fig:pi_plus_u2_tridia} shows three diagrams corresponding to the unit cell shown in Fig.~\ref{fig:pi_plus_u2_uc_xyz} along three non-coplanar axes. Such a set of three diagrams constitutes a \emph{tridiagram}. All three diagrams are necessary for converting the full spectrum of transformations in $3$-space into transformations in 2-dimensional diagrams, where the crossing information of each diagram is different and equally important. Unless otherwise mentioned, every tridiagram that we consider in this paper will be obtained from the projections from the front, the top and the right faces of a unit cell.

\subsection{Equivalence of 3-periodic tangles}
Deforming a filamentous structure in space without cutting or reconnecting the filaments does not change how entangled the structure is. Mathematically, such deformations are described by the notion of ambient isotopy and include bending, stretching, or smoothly moving the filaments. Two filamentous structures that can be transformed into one another by an ambient isotopy are said to belong to the same \emph{isotopy class}. The concepts of least tangled embeddings and the untangling number, introduced in a later section, are invariant under ambient isotopies performed within a unit cell. We therefore briefly introduce here some notions concerning ambient isotopy of 3-periodic tangles that are used throughout the remainder of the paper. Further details are published in \cite{ANDRIAMANALINA2025109346}.

\bigbreak

\begin{figure}[hbtp]
    \centering
    \begin{subfigure}[b]{4.6cm}
        \includegraphics[width=\textwidth]{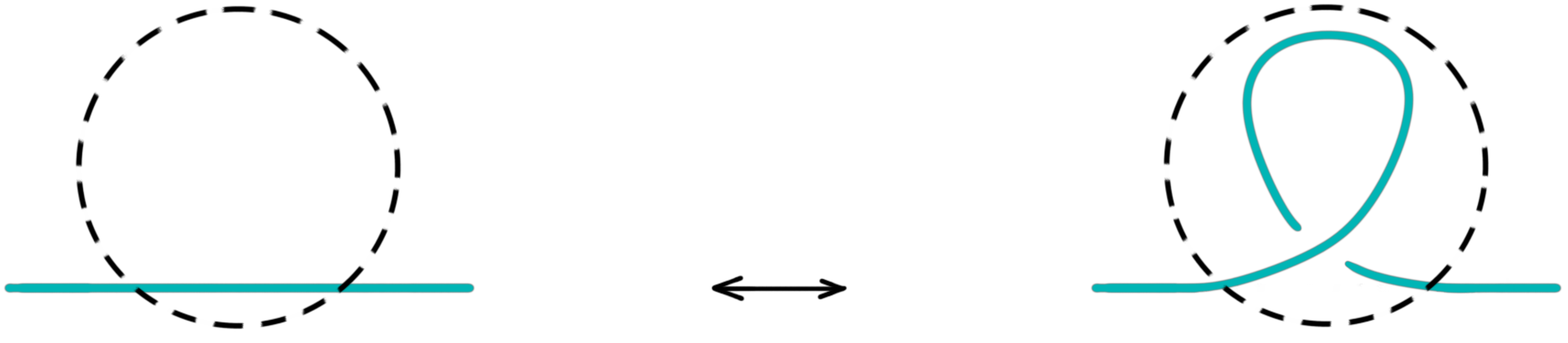}
        \caption{The $R_1$ move}
        \label{fig:Reid_move_I}
    \end{subfigure}
    \hspace{0.5cm}
    \begin{subfigure}[b]{5cm}
        \includegraphics[width=\textwidth]{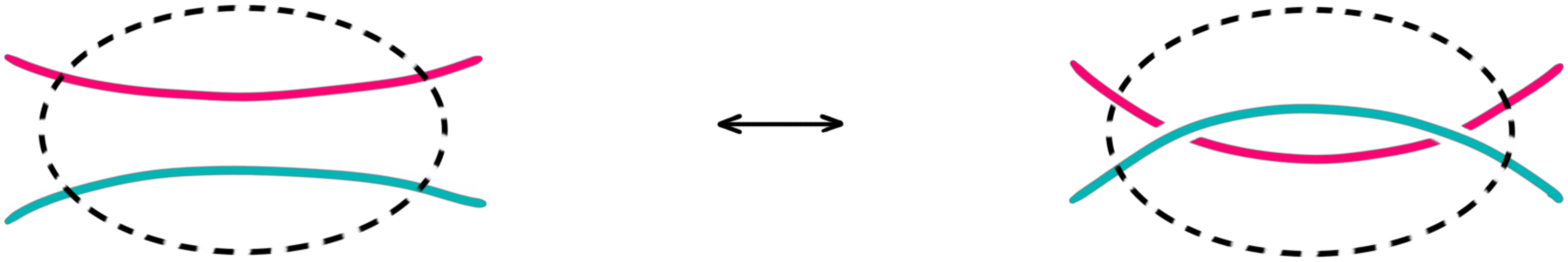}
        \caption{The $R_2$ move}
        \label{fig:Reid_move_II}
    \end{subfigure}
    \vskip\baselineskip
    \begin{subfigure}[b]{4.6cm}
        \includegraphics[width=\textwidth]{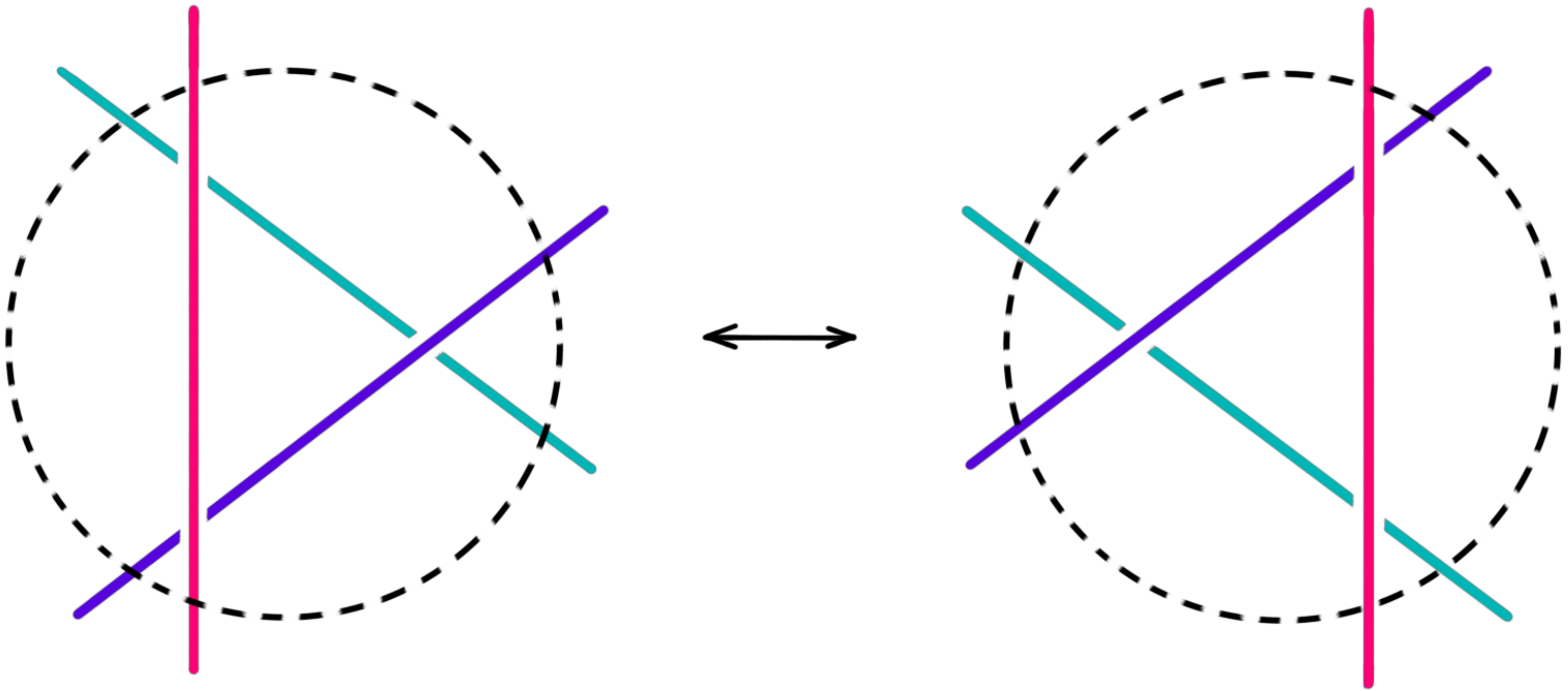}
        \caption{The $R_3$ move}
        \label{fig:Reid_move_III}
    \end{subfigure}
    \caption{The three Reidemeister moves, originally defined to capture isotopies of classical knots and links in $\mathbb{R}^3$.}
    \label{fig:usual_Reid_moves}
\end{figure}

\begin{figure}[hbtp]
    \centering
    \begin{subfigure}[b]{5cm}
    \centering
        \includegraphics[width=\textwidth]{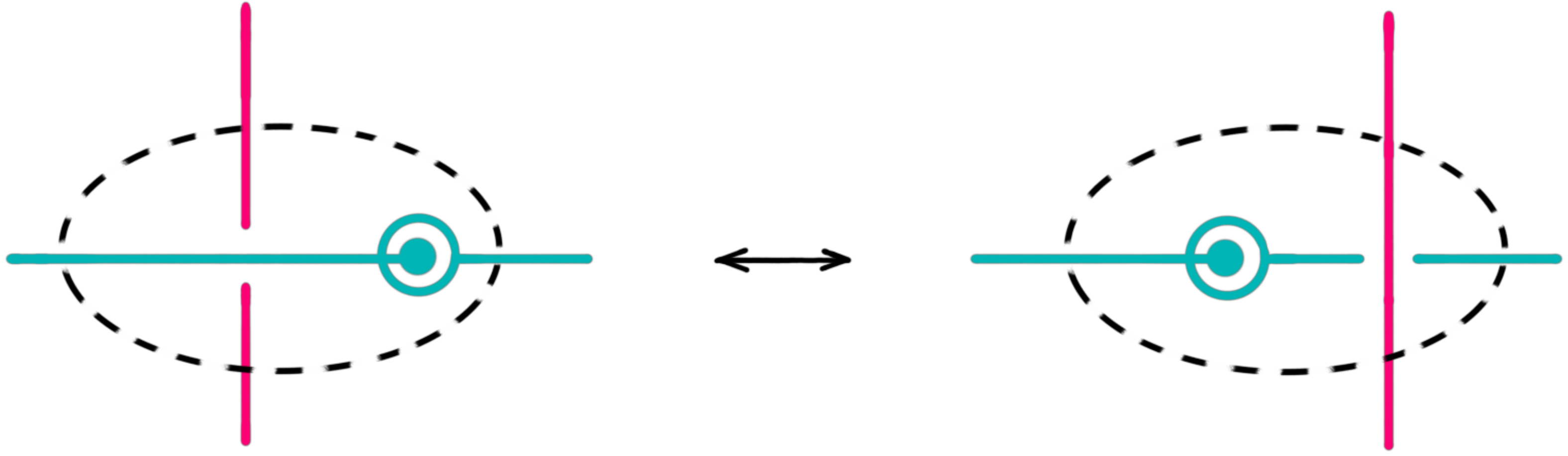}
        \caption{The $R_4$ move}
        \label{fig:Reid_move_IV}
        \vspace{0.25cm}
    \end{subfigure}
    \hspace{0.5cm}
    \begin{subfigure}[b]{5cm}
    \centering
        \includegraphics[width=\textwidth]{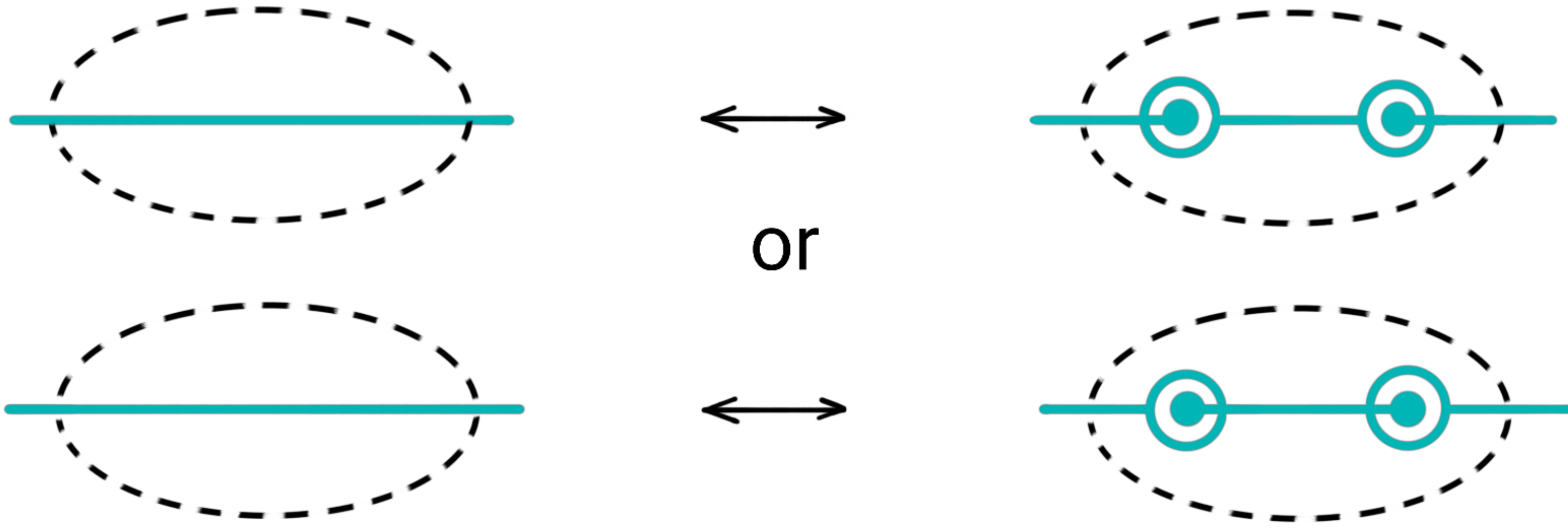}
        \caption{The $R_5$ move}
        \label{fig:Reid_move_V}
    \end{subfigure}
    
    \vskip\baselineskip
    
    \begin{subfigure}[b]{4.6cm}
    \centering
        \includegraphics[width=\textwidth]{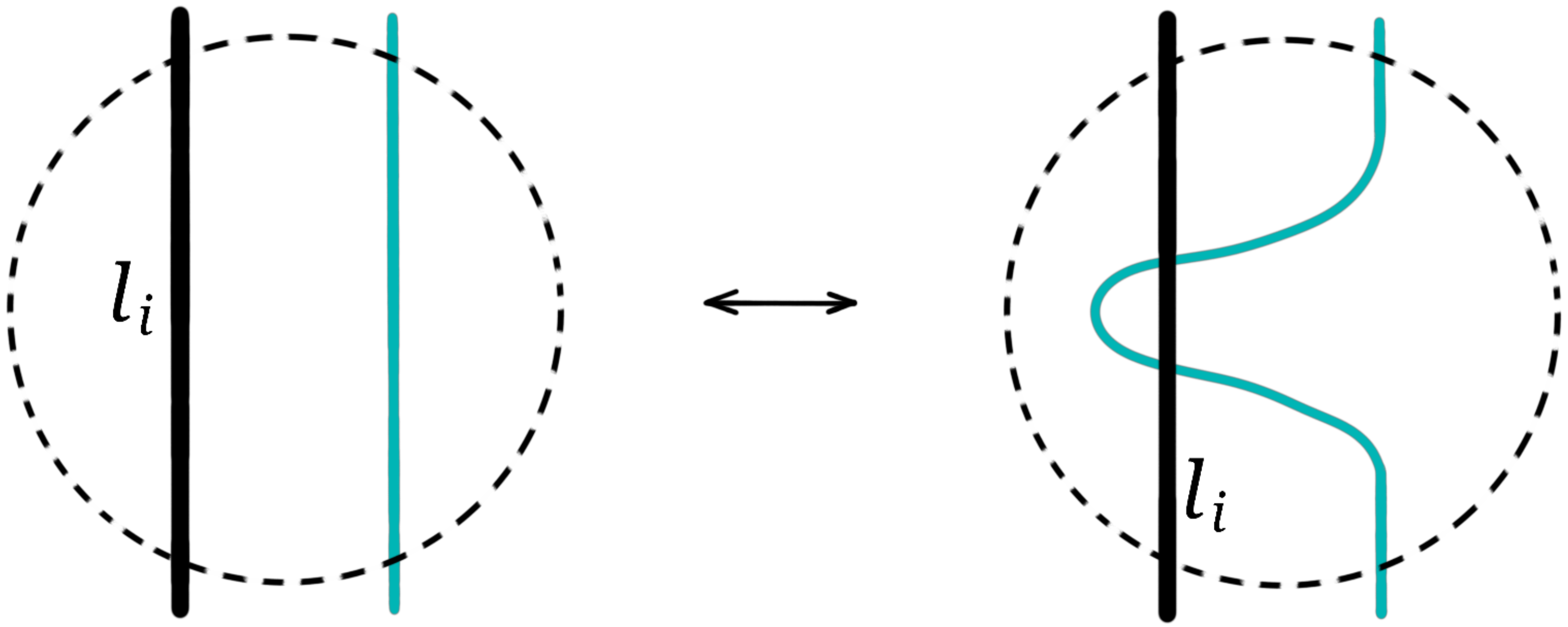}
        \caption{The $R_6$ move}
        \vspace{1cm}
        \label{fig:Reid_move_VI}
    \end{subfigure}
    \hspace{1cm}
    \begin{subfigure}[b]{4.6cm}
    \centering
        \includegraphics[width=\textwidth]{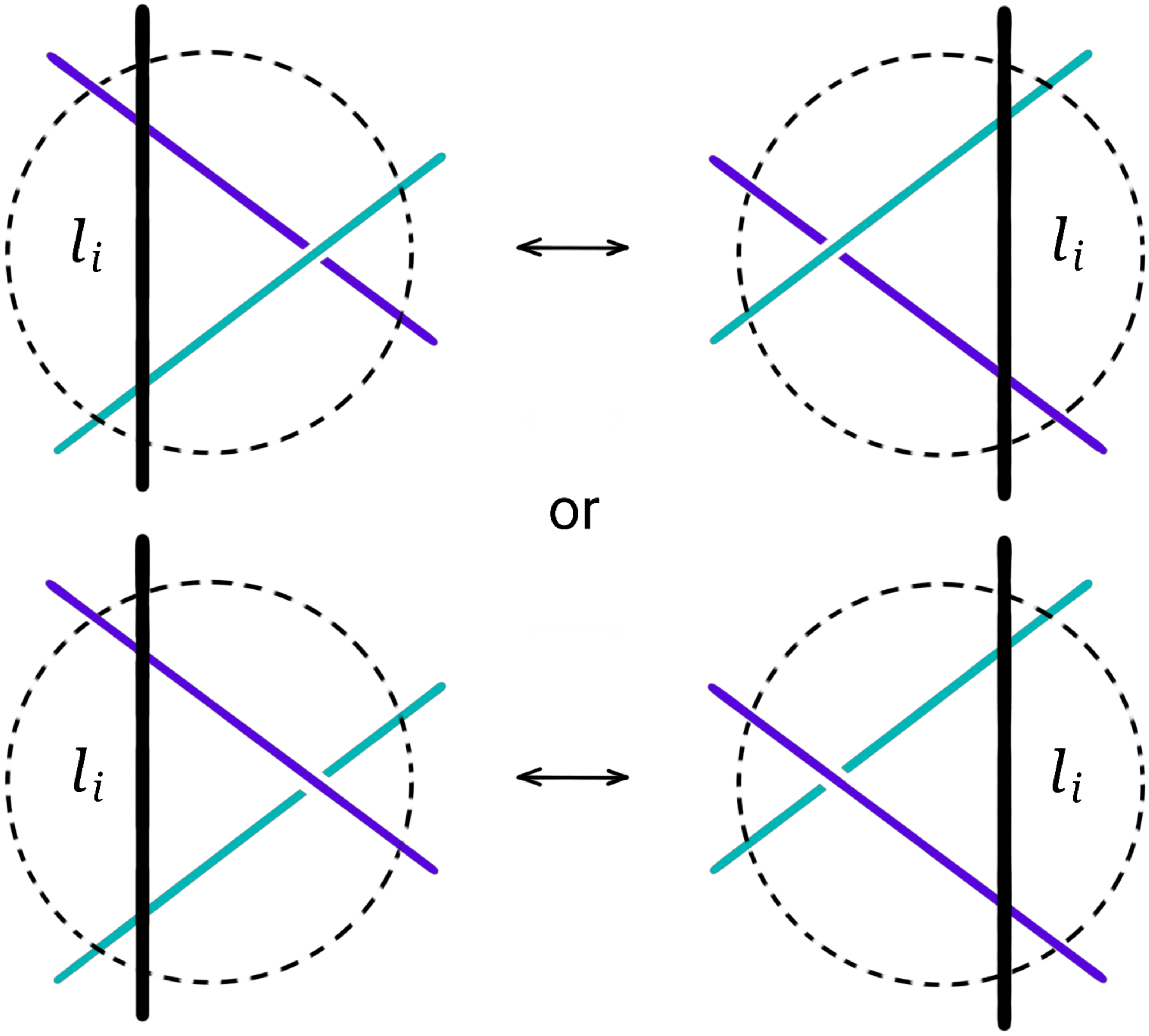}
        \caption{The $R_7$ move}
        \label{fig:Reid_move_VII}
    \end{subfigure}
    
    \vskip\baselineskip
    
    \begin{subfigure}[b]{4.6cm}
    \centering
        \includegraphics[width=\textwidth]{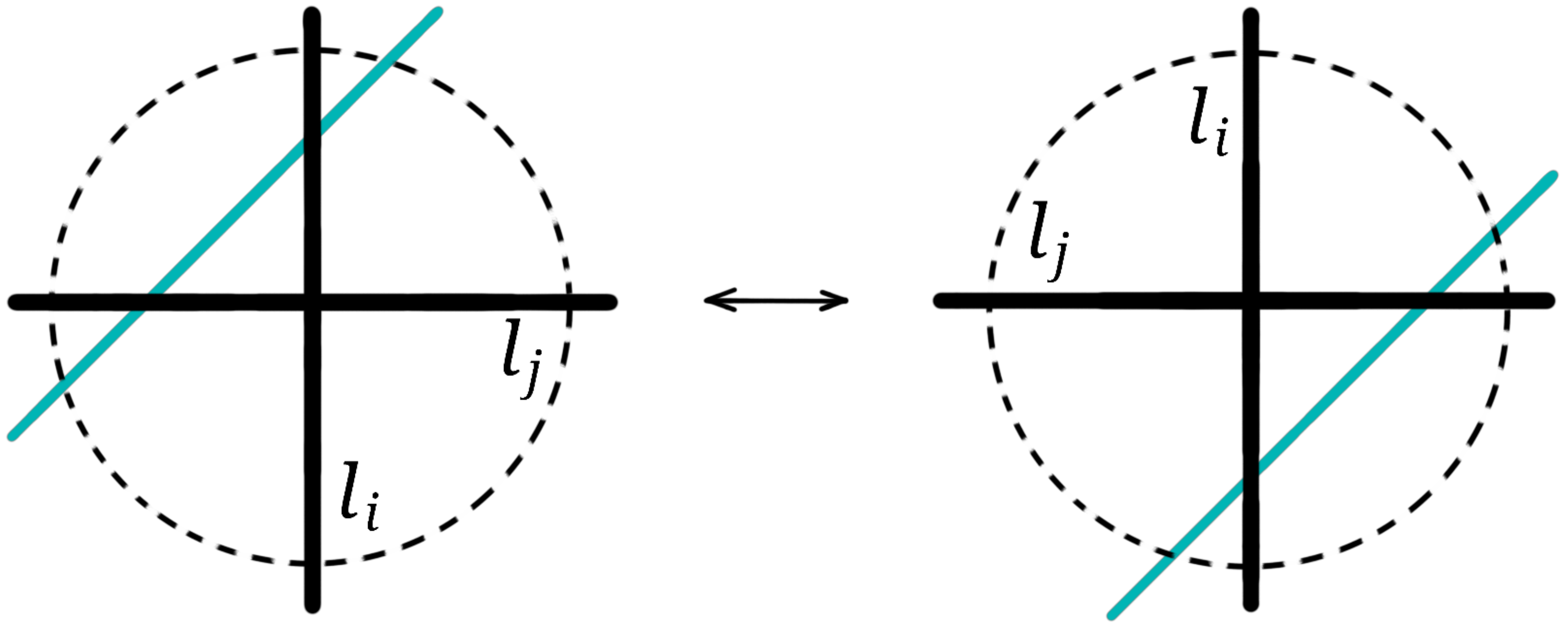}
        \caption{The $R_8$ move}
        \label{fig:Reid_move_VIII}
    \end{subfigure}
    \hspace{0.75cm}
    \begin{subfigure}[b]{5cm}
    \centering
        \includegraphics[width=\textwidth]{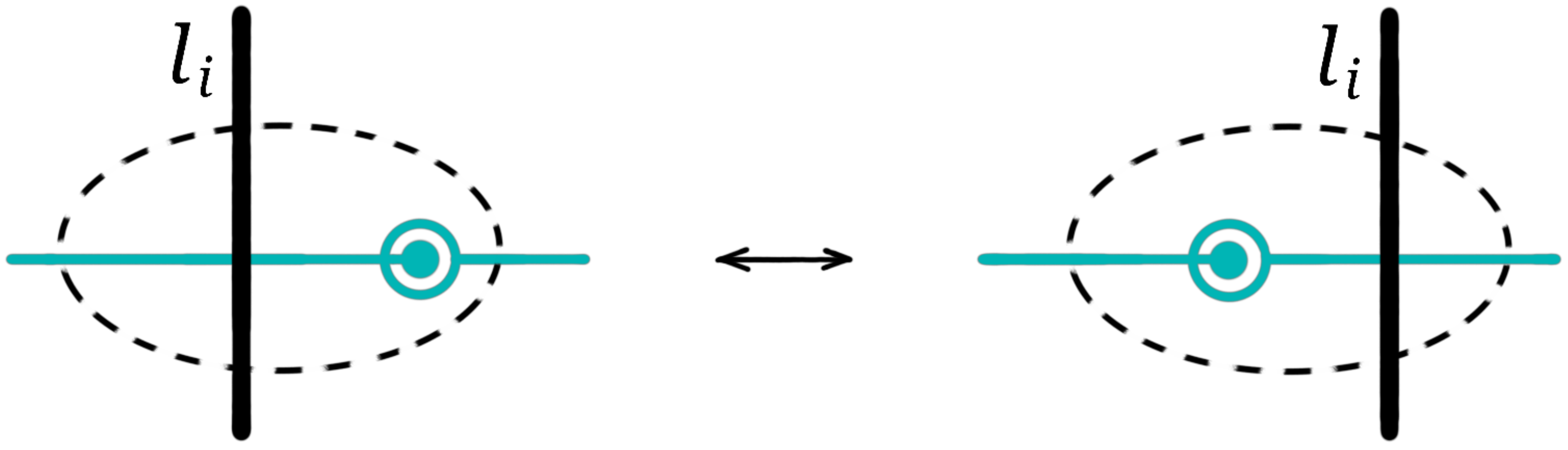}
        \caption{The $R_9$ move}
        \label{fig:Reid_move_IX}
        \vspace{0.35cm}
    \end{subfigure}
        \caption{The six new moves, necessary to capture ambient isotopies within unit cells of 3-periodic tangles. Here $l_i$ and $l_j$, where $i,j = 1,2$, represent the edges of the square delimiting a diagram.}
    \label{fig:new_Reid_moves}
\end{figure}

\begin{figure}[hbtp]
    \centering
        \includegraphics[width=0.18525\textwidth]{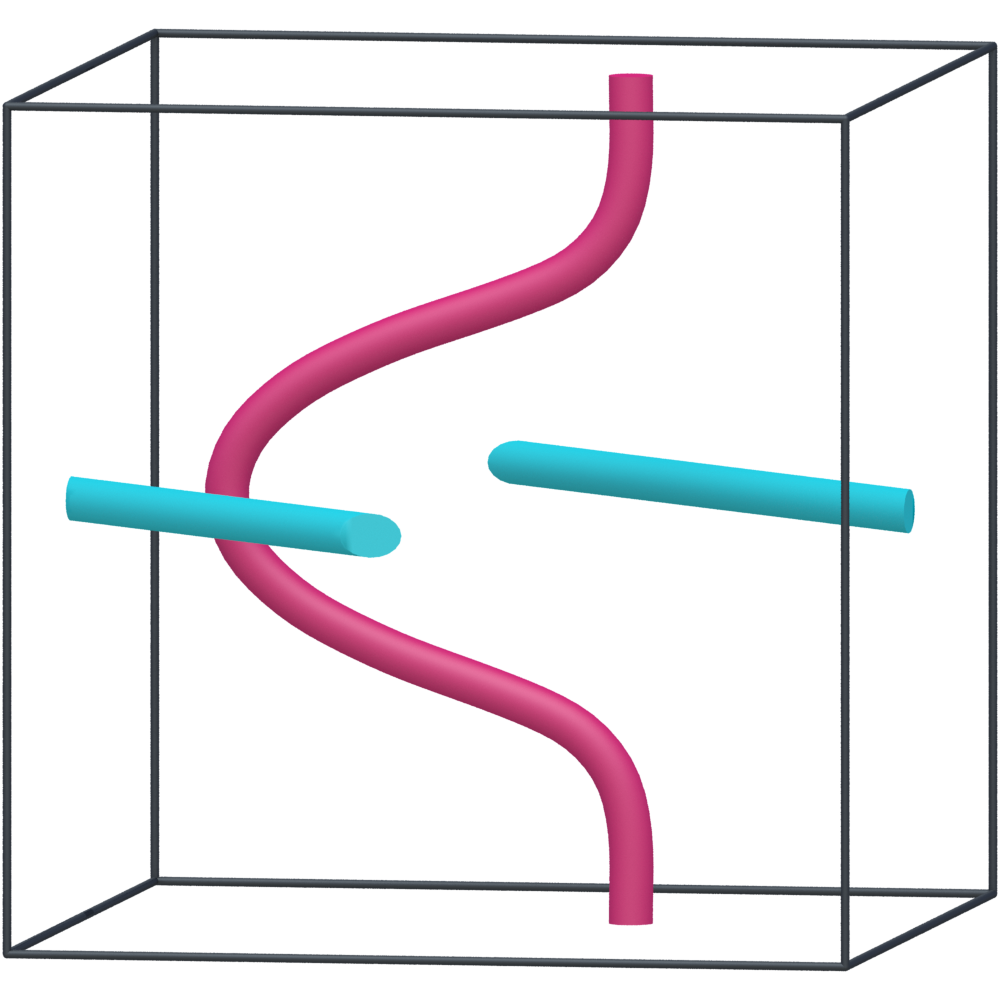}
        \hspace{0.05\textwidth}
        \includegraphics[width=0.18525\textwidth]{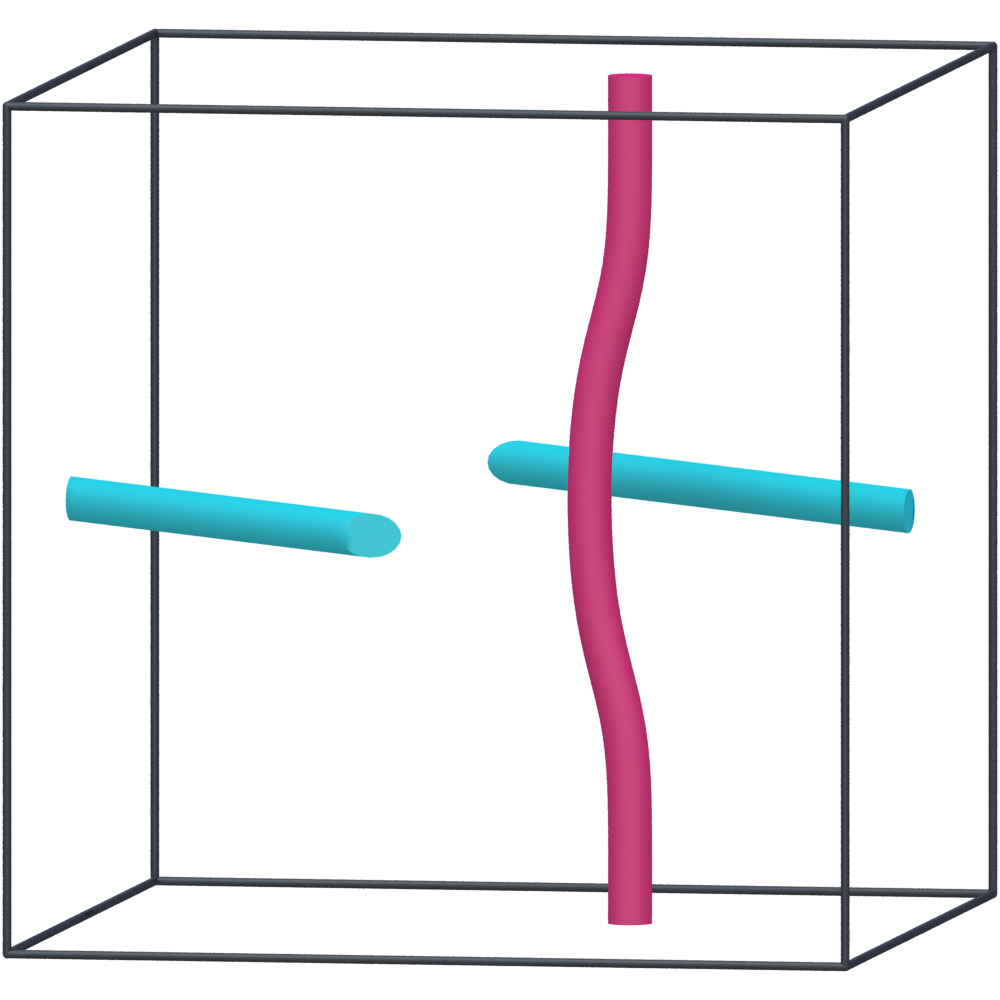}
    \caption{Two isotopic unit cells that differ by an $R_4$ move in diagrams obtained from projections from their front faces. The move reverses the type of a crossing while preserving the isotopy class of the 3-periodic tangle.}
    \label{fig:R4_3D}
\end{figure}

Isotopy transformations performed within a unit cell of a 3-periodic tangle are represented in a diagram by combinatorial moves called $R$-moves. The $R$-moves comprise the three Reidemeister moves $R_1$, $R_2$, and $R_3$, shown in Fig.~\ref{fig:usual_Reid_moves}, originally defined to characterise ambient isotopies of classical knots in space, together with six additional moves, $R_4, \dots, R_9$, shown in Fig.~\ref{fig:new_Reid_moves}, specific to 3-periodic tangles. In particular, the $R_4$ move, shown in Fig.~\ref{fig:Reid_move_IV}, reverses the type of a crossing without changing the isotopy class of a 3-periodic tangle. It corresponds to deformations such as the one shown in Fig.~\ref{fig:R4_3D}.

\bigbreak

From now on, unless otherwise mentioned, we regard each unit cell of a given 3-periodic tangle as a representative of its ambient isotopy class.

\bigbreak

\begin{figure}[hbtp]
    \centering    
    \begin{subfigure}[b]{\textwidth}
    \centering
        \includegraphics[width=0.37\textwidth]{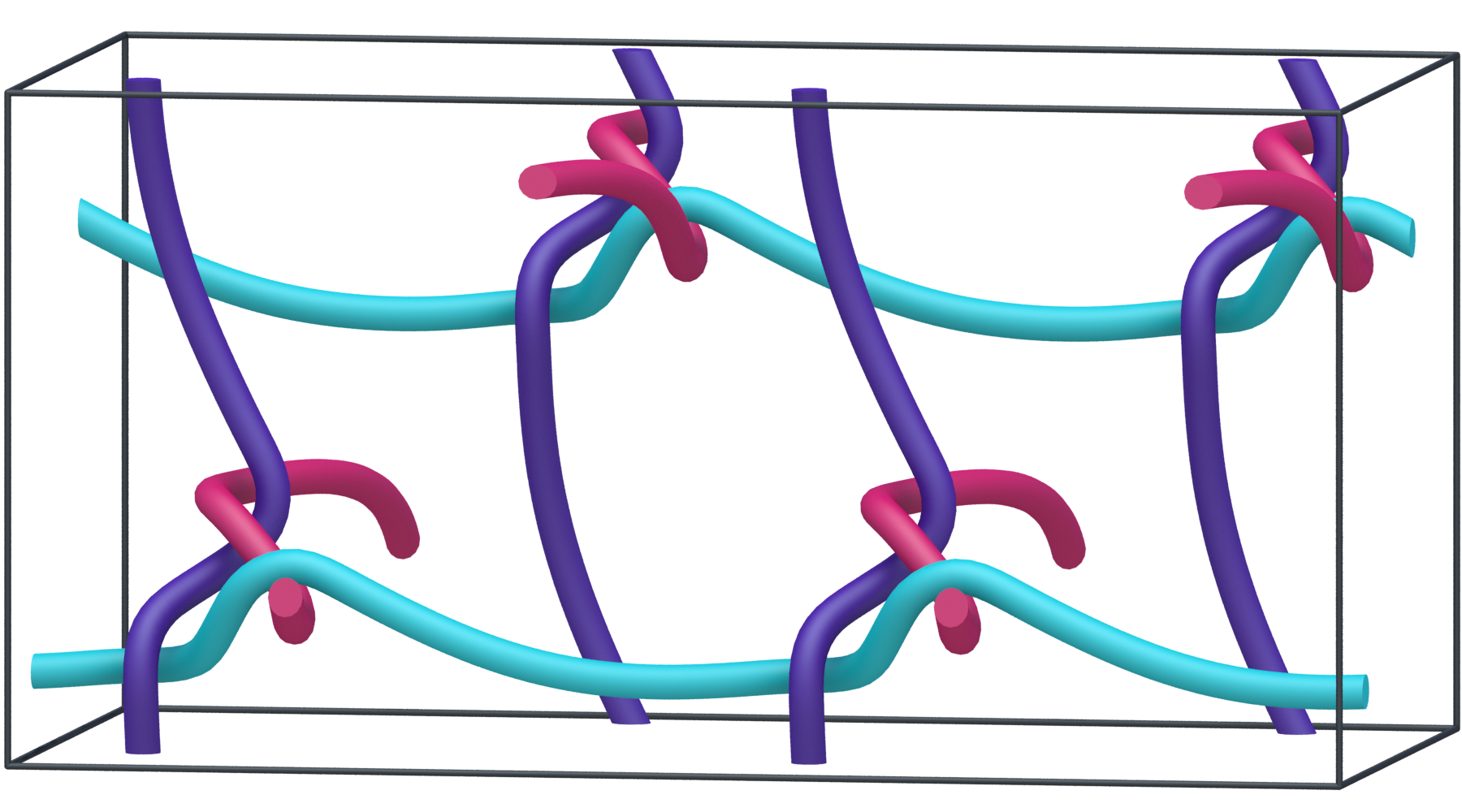}
        \caption{}
        \label{fig:pi_plus_u2_uc_vol_2}
    \end{subfigure}

    \vskip\baselineskip

    \begin{subfigure}[b]{\textwidth}
    \centering
        \includegraphics[width=0.37\textwidth]{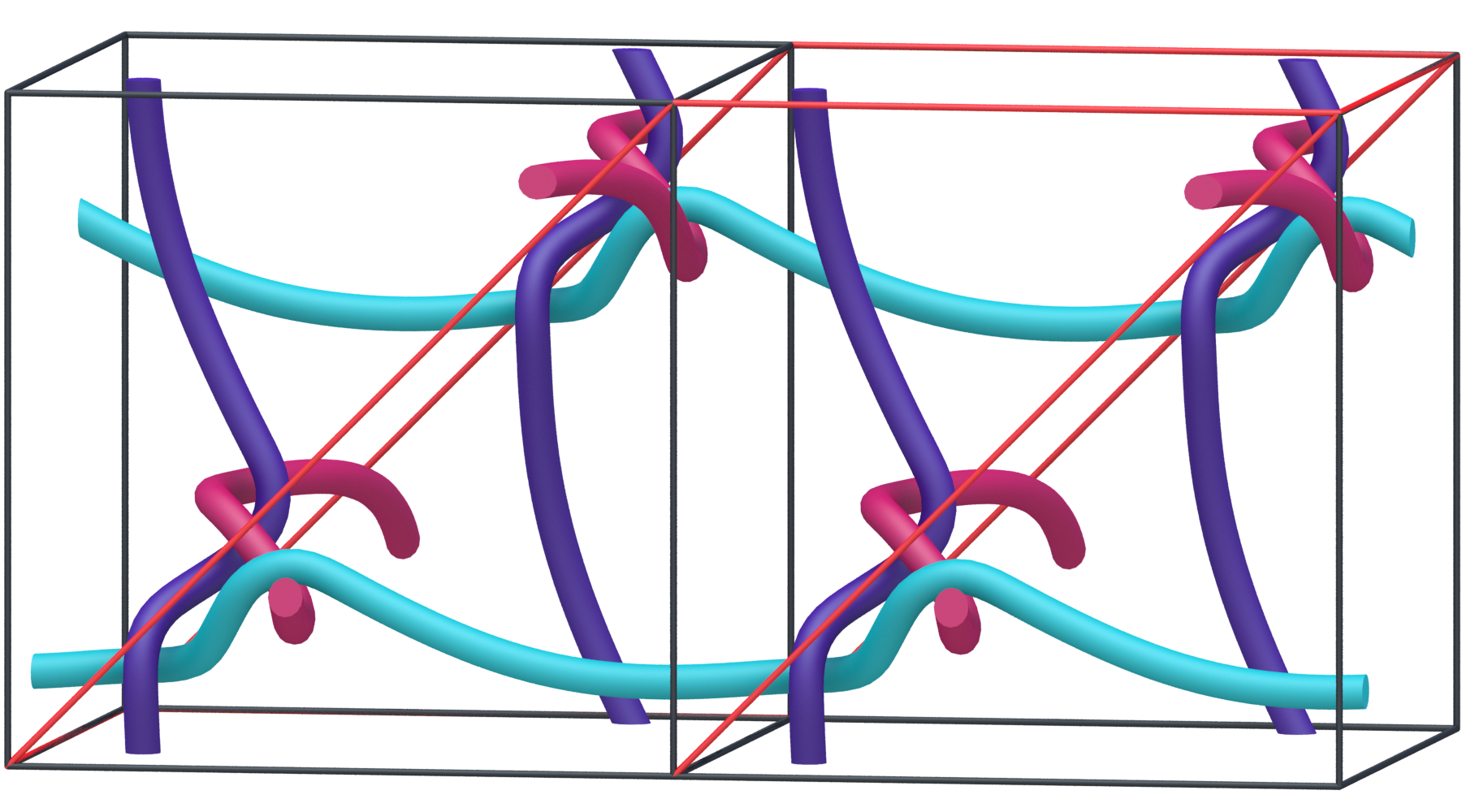}
        \caption{}
        \label{fig:pi_plus_u2_sheared_uc}
    \end{subfigure}
    \caption{Non-isotopic unit cells that represent the same 3-periodic tangle. (a)  A double of the unit cell of Fig.~\ref{fig:pi_plus_u2_uc_xyz}. (b) A red unit cell obtained from a different parallelepiped than that of the unit cell of Fig.~\ref{fig:pi_plus_u2_uc_xyz}.}
    \label{fig:pi_plus_u2_other_uc}
\end{figure}

We remark that there exist infinitely many non-isotopic unit cells that represent the same 3-periodic tangle, and hence naturally preserve its entanglement complexity. Examples are given in Fig.~\ref{fig:pi_plus_u2_other_uc}, where the unit cell shown in Fig.~\ref{fig:pi_plus_u2_uc_vol_2} is obtained by doubling the size of the unit cell in Fig.~\ref{fig:pi_plus_u2_uc_xyz}, and the red unit cell in Fig.~\ref{fig:pi_plus_u2_sheared_uc} has the same volume as that in Fig.~\ref{fig:pi_plus_u2_uc_xyz} but is obtained from a different parallelepiped.

\subsection{Crossing number}
Diagrams are valuable in that they allow the classification of three-dimensional tangled structures to be reduced to the classification of their two-dimensional projections. The latter is done by using \emph{invariants}, which are quantities that do not vary under changes performed on a diagram that preserve the isotopy class of the structure it represents. For 3-periodic tangles, an example of a diagrammatic invariant defined in \cite{ANDRIAMANALINA2025109346} is the \emph{crossing number}, which we recall here. Consider a 3-periodic tangle $K$ and a unit cell $U$ of $K$. To a tridiagram $T = \lbrace D_i \rbrace_{i=1,2,3}$ obtained from $U$, one can associate a triplet $(a,b,c)$, where $a$, $b$ and $c$ are the numbers of crossings of the diagrams $D_1$, $D_2$ and $D_3$, respectively. Such a triplet is called a \emph{triplet of crossings}. The \emph{crossing number of $K$ with respect to $U$}, denoted by $c(K,U)$, is the minimum of $c(T) = a^2 + b^2 +c^2$ among all tridiagrams $T$ obtained from all representatives of the isotopy class of the unit cell $U$. A tridiagram and its triplet of crossings realising $c(K,U)$ are respectively called a \emph{minimal tridiagram of $K$ with respect to $U$}, and a \emph{minimum crossing number triplet with respect to $U$}. For example, the tridiagram of the $\Pi^{+}$ rod packing shown in Fig.~\ref{fig:pi_plus_tridia} is a minimal tridiagram with respect to the unit cell displayed in Fig.~\ref{fig:pi_plus_uc_xyz}, and the minimum crossing number triplet is $(4,4,4)$. The \emph{crossing number} of $K$ is defined as the minimum of all $c(K,U)$ among all unit cells $U$ of $K$, and it is an invariant.

Note that, as explained in \cite{ANDRIAMANALINA2025109346}, we choose to sum the squares of the numbers of crossings rather than the numbers themselves in order to bias the unit cell towards a more symmetric configuration. In the 3-periodic setting, the crossing number is primarily relevant for ordering purposes; the minimum crossing number triplet, which indicates the numbers of crossings along the three directions, is the more significant quantity.

\section{Least tangled embeddings and crystallographic rod packings}\label{sec:ground_states}

Using the framework established in Section \ref{sec:reminder_diagrams}, namely crossing diagrams, we are now able to rigorously define and characterise the least tangled embeddings of 3-periodic tangles, which we refer to as ground states.

\bigbreak

\begin{figure}[hbtp]
    \centering
    \begin{subfigure}[b]{0.95\textwidth}
        \centering
        \includegraphics[width=0.65\textwidth]{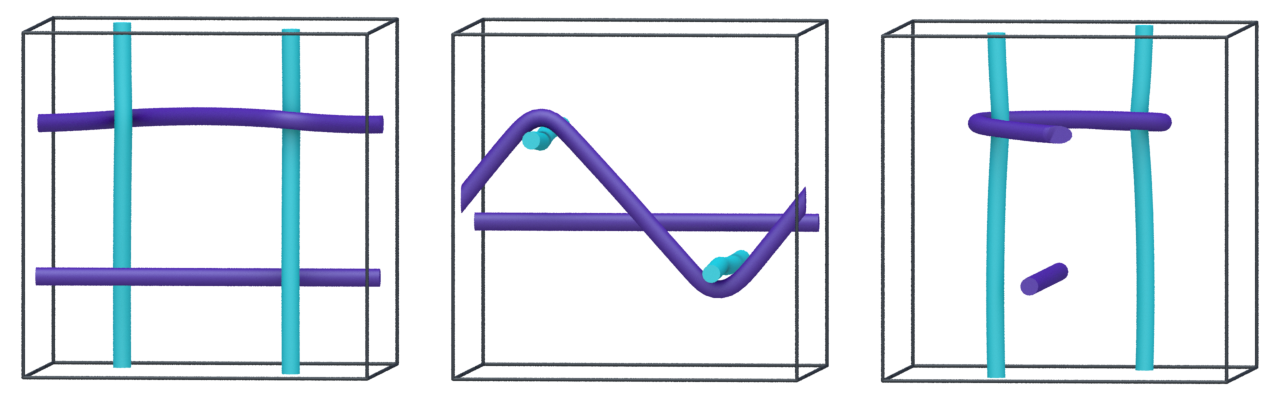}
        \caption{}
        \label{fig:plain_weave_xyz_uc}
    \end{subfigure}

    \vskip\baselineskip

    \begin{subfigure}[b]{0.2\textwidth}
        \centering
        \includegraphics[width=0.95\textwidth]{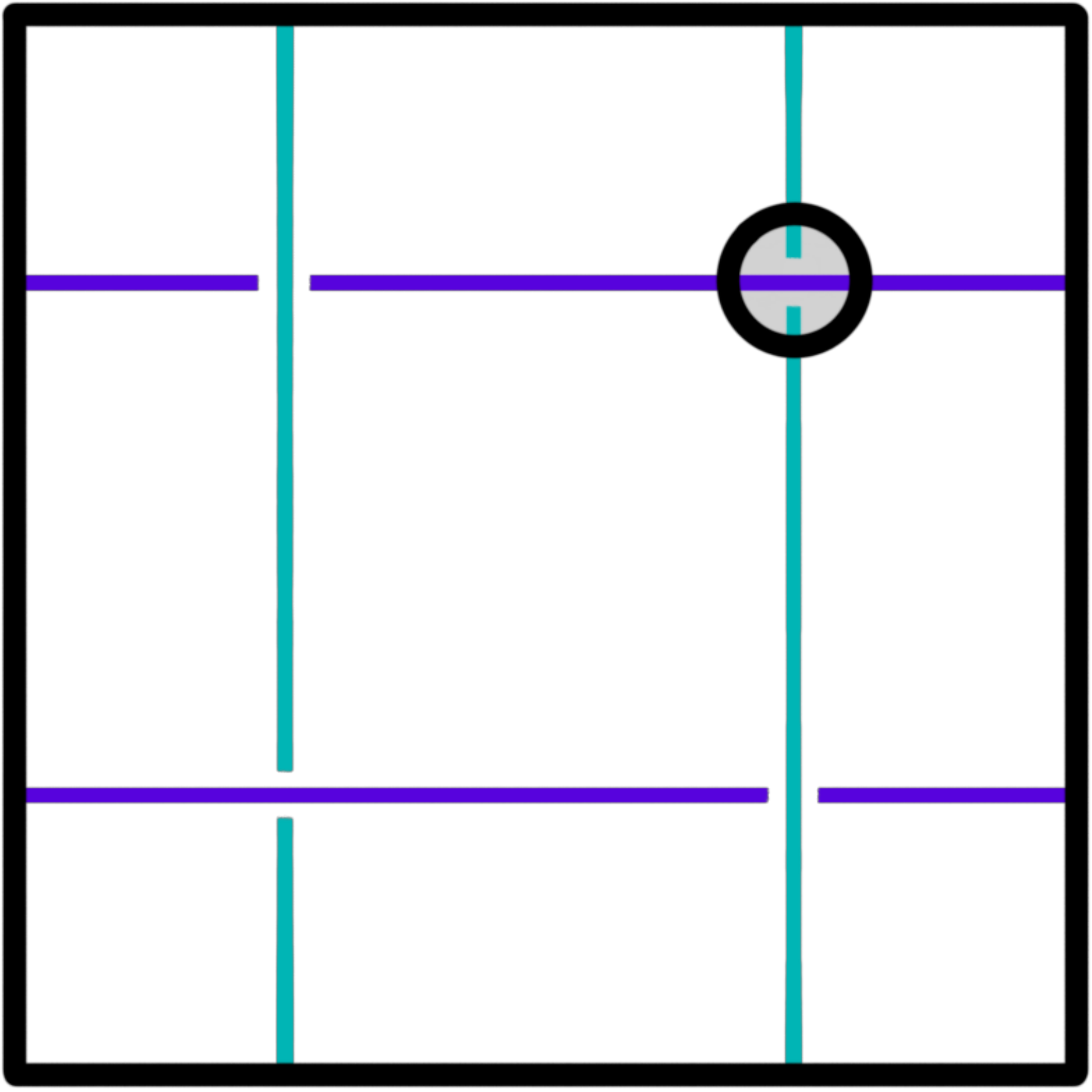}
        \caption{}
        \label{fig:plain_weave_dia_F}
    \end{subfigure}
    \hspace{0.025\textwidth}
    \begin{subfigure}[b]{0.2\textwidth}
        \centering
        \includegraphics[width=0.95\textwidth]{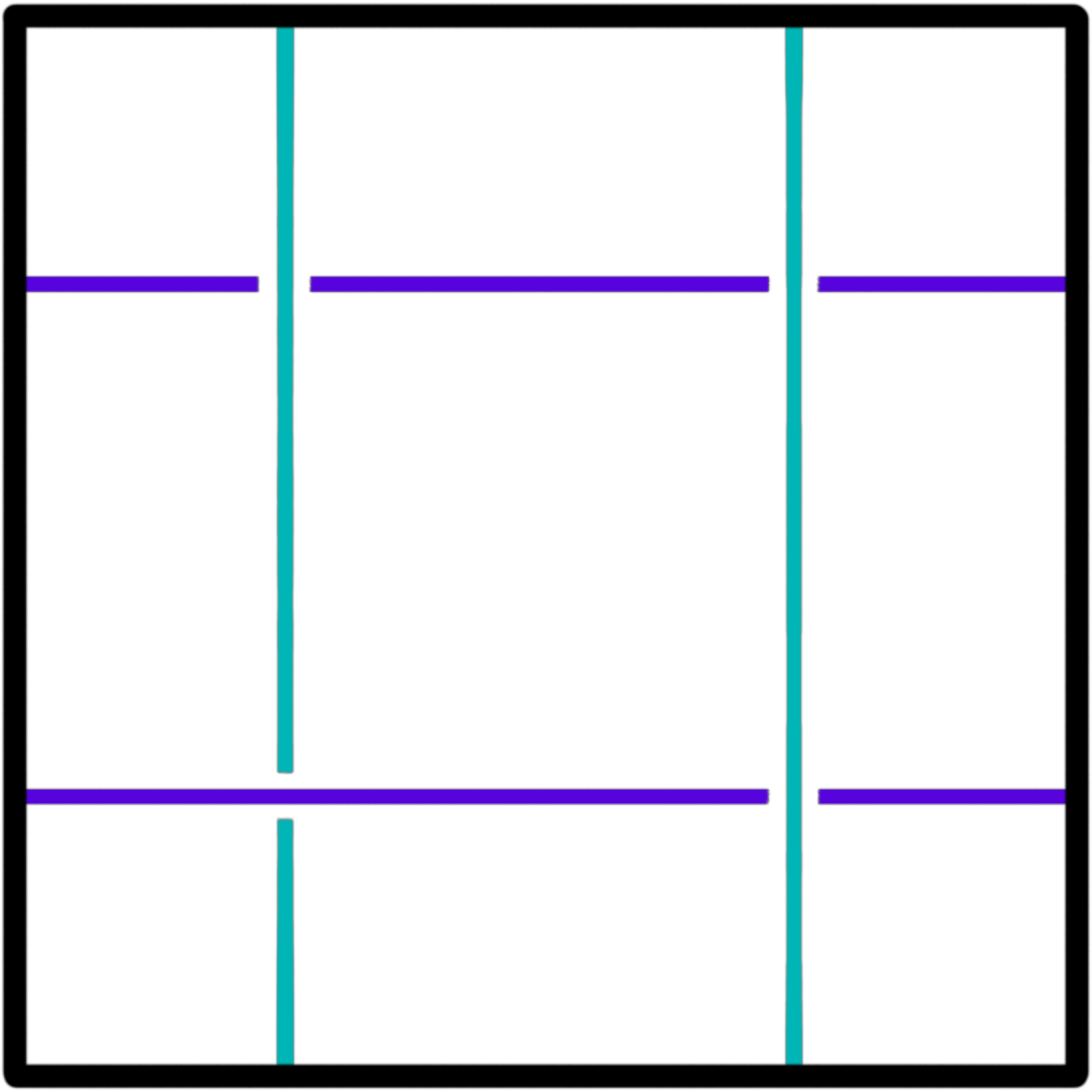}
        \caption{}
        \label{fig:plain_weave_ground_state_dia_F}
    \end{subfigure}

    \vskip\baselineskip

    \begin{subfigure}[b]{0.95\textwidth}
        \centering
        \includegraphics[width=0.65\textwidth]{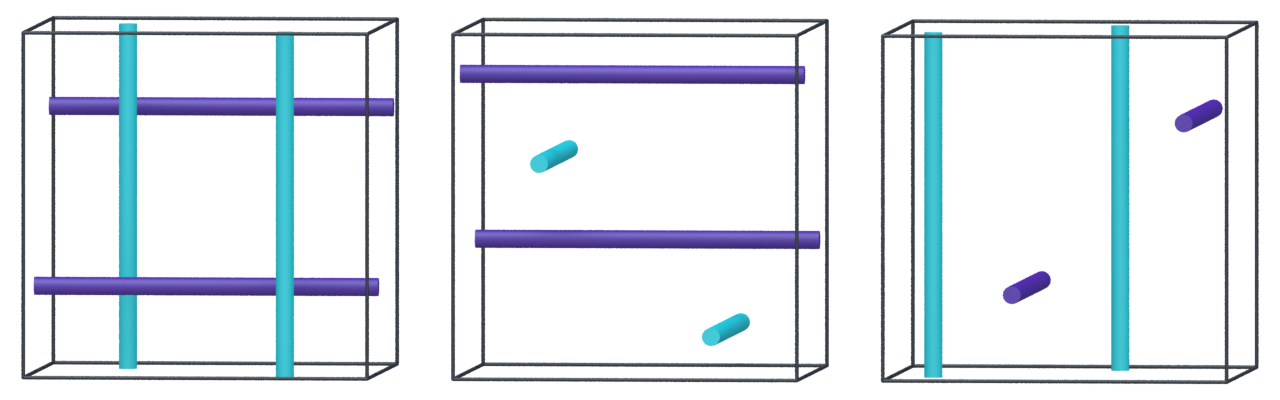}
        \caption{}
        \label{fig:plain_weave_ground_state_xyz_uc}
    \end{subfigure}
    \caption{Applying a crossing change on a diagram of a periodic structure. (a) One unit cell of the plain weave viewed along three directions: the front, the top and the right sides. (b) A diagram obtained from the unit cell shown in (a), where a crossing change is applied onto the highlighted crossing to obtain the diagram in (c). (d) A unit cell obtained from the diagram shown in (c), which presents fewer crossings (in a tridiagram) than the unit cell of the plain weave.}
    \label{fig:untangling_plain_weave}
\end{figure}

Classical knots and links are grouped by their number of components. For a given positive integer $n$, the \emph{unlink} is a disjoint collection of $n$ circles, and any link with $n$ components can be transformed into the unlink via crossing changes. See Fig.~\ref{fig:trefoil_b} for the case $n=1$, where the trefoil knot is unknotted. By analogy, given a 3-periodic tangle $K$ and a unit cell $U$, we define the \emph{$\mathcal{U}$-family of $K$ with respect to $U$} to be set of all 3-periodic tangles that possess a unit cell that is connected to $U$ by finitely many crossing changes. Within this family, we define the subfamily consisting of structures with the common crossing number that is the least among all elements of the family, as the \emph{$\mathcal{G}$-family of $K$ with respect to $U$}, and refer to its elements as the \emph{ground states}. For example, in Fig.~\ref{fig:untangling_plain_weave}, we apply a crossing change to a diagram of a given periodic structure, the plain weave, which subsequently untangles it into a structure with fewer crossings. We will later provide a sufficient condition that shows that this latter structure is indeed an element of the $\mathcal{G}$-family of the plain weave.

\begin{remark}
    The plain weave, whose unit cell is shown in Fig.~\ref{fig:plain_weave_xyz_uc}, is technically a 2-periodic tangle, but a 3-periodic tangle can be obtained from it by layering infinitely many copies along the third direction of space. Considering it as a 2-periodic or 3-periodic structure changes the $R$-moves that can be used in its diagrams.
\end{remark}

\bigbreak

In the following, we provide additional mathematical precision to the previous notions.

\begin{definition}
   Let $K$ be a 3-periodic tangle and $U$ a unit cell of $K$. The \emph{$\mathcal{U}$-family} of $K$ with respect to $U$, denoted by $\mathcal{U}(K,U)$, is the set containing $K$ and any 3-periodic tangle $K'$ that possesses a unit cell $U'$ connected to $U$ by finitely many crossing changes.
\end{definition}

Note that any 3-periodic tangle $K'$ belonging to $\mathcal{U}(K,U)$ is always understood to be paired with the unit cell $U'$ that is connected to $U$ by finitely many crossing changes, and not with any other unit cell representing $K'$. Within the $\mathcal{U}$-family $\mathcal{U}(K,U)$, there is a subfamily of 3-periodic tangles that share the same crossing number, which is the least for all 3-periodic tangles of $\mathcal{U}(K,U)$. We define this subfamily and its elements as follows.

\begin{definition}
    Let $K$ be a 3-periodic tangle and $U$ a unit cell of $K$. We define the \emph{$\mathcal{G}$-family} of $K$ with respect to $U$, denoted by $\mathcal{G}(K,U)$, as the subset of $\mathcal{U}(K,U)$ consisting of all 3-periodic tangles having the common crossing number that is the least among all elements of $\mathcal{U}(K,U)$. The elements of $\mathcal{G}(K,U)$ are called the \emph{ground states}.
\end{definition}

Notice that if $K'$ belongs to $\mathcal{U}(K,U)$ with respect to the unit cell $U'$, then $K$ belongs to $\mathcal{U}(K',U')$ with respect to $U$, and thus $\mathcal{U}(K,U) = \mathcal{U}(K',U')$. Notice also that, by the definition of $\mathcal{G}(K,U)$, if $K'$ belongs to $\mathcal{U}(K,U)$ with respect to the unit cell $U'$, then we have $\mathcal{G}(K',U') =\mathcal{G}(K,U)$.

\bigbreak

The definition of a ground state creates the need for a method to determine the least possible crossing number for a given $\mathcal{U}$-family. The mathematical notion of homotopy immediately provides such a method. Loosely speaking, two unit cells are freely homotopic or belong to the same free homotopy class, if one can be transformed into the other by deforming the curve components in space while allowing them to pass through themselves and each other. This leads to the following theorem on the minimum crossing number triplets of ground states.

\begin{theorem}\label{thm:homotopy_ground_states}
    Consider a 3-periodic tangle $K$ and a unit cell $U$ of $K$. Consider a tridiagram $\lbrace D_1,D_2,D_3 \rbrace$ obtained from $U$, and regard each diagram $D_i$ as a $2$-torus onto which curves lie. For each diagram, consider the minimum number of intersections among the free homotopy classes of the curves, including both intersections between distinct curves and self-intersections, and denote this number by $\mathrm{i}(D_i)$. If the minimum crossing number triplet of $K$ with respect to $U$ equals the triplet $(\mathrm{i}(D_1),\mathrm{i}(D_2),\mathrm{i}(D_3))$, then $K$ is a ground state with respect to $U$.
\end{theorem}

\begin{proof}
    Homotopy deformations transforming one unit cell into another correspond in a diagram to $R$-moves and crossing changes. By regarding the crossings in a diagram as intersections of curves sitting on the surface of the $2$-torus, determining the least crossing number of a given $\mathcal{U}$-family boils down to minimising the number of intersections in a tridiagram by applying homotopy deformations. Therefore, the ground states are precisely the structures that possess a minimal tridiagram whose crossings coincide with the intersections that cannot be eliminated under homotopy deformations, the characterisation of which is already well known for closed curves on the surface of the $2$-torus \cite{primer_mapping_cl_group_chap1,GU2022108205}.
\end{proof}

\begin{figure}[hbtp]
    \centering
    \begin{subfigure}[b]{0.19\textwidth}
        \centering
        \includegraphics[width=0.975\textwidth]{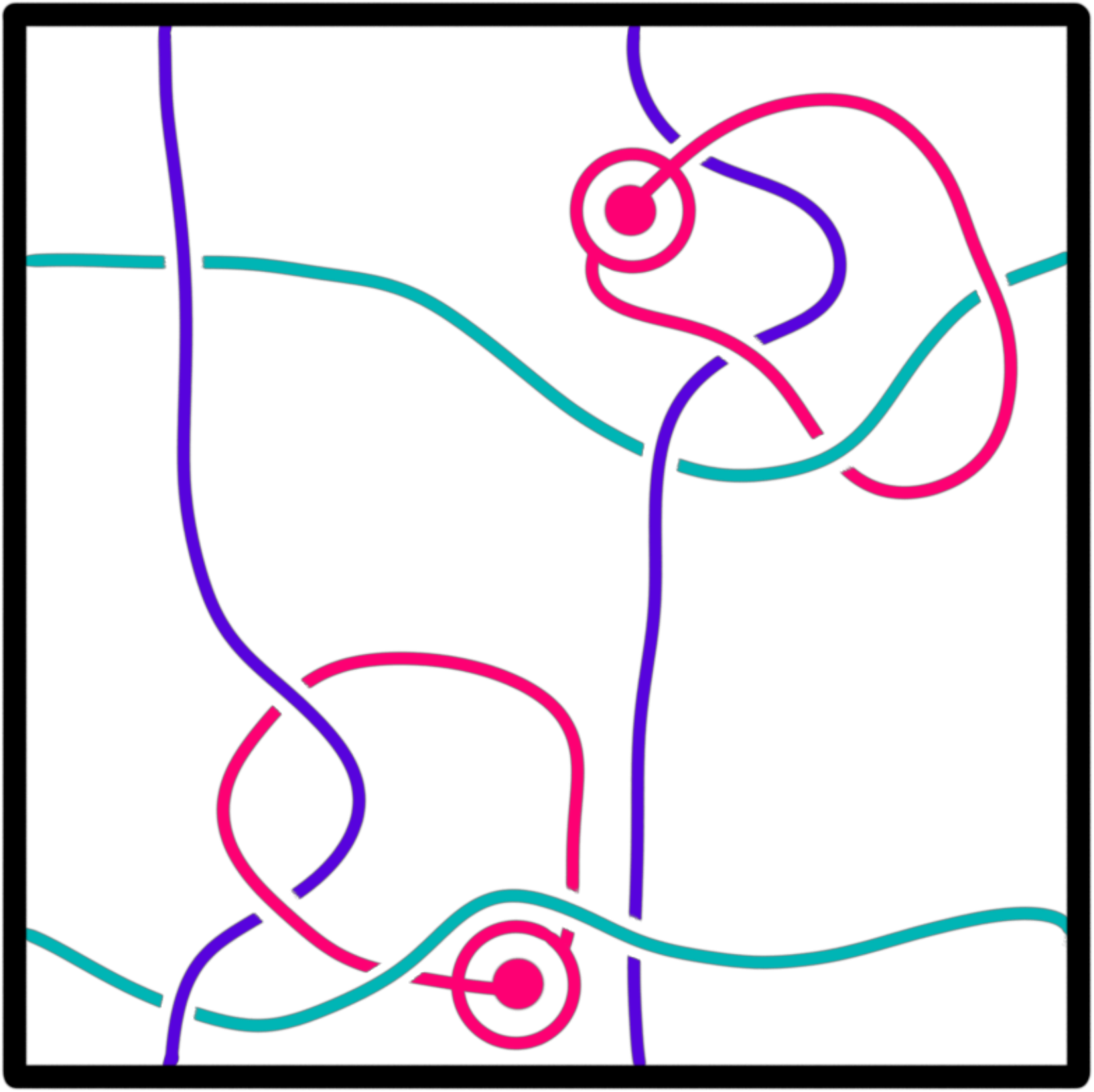}
        \caption{}
        \label{fig:homotopy_1}
    \end{subfigure}
    \begin{subfigure}[b]{0.19\textwidth}
        \centering
        \includegraphics[width=0.975\textwidth]{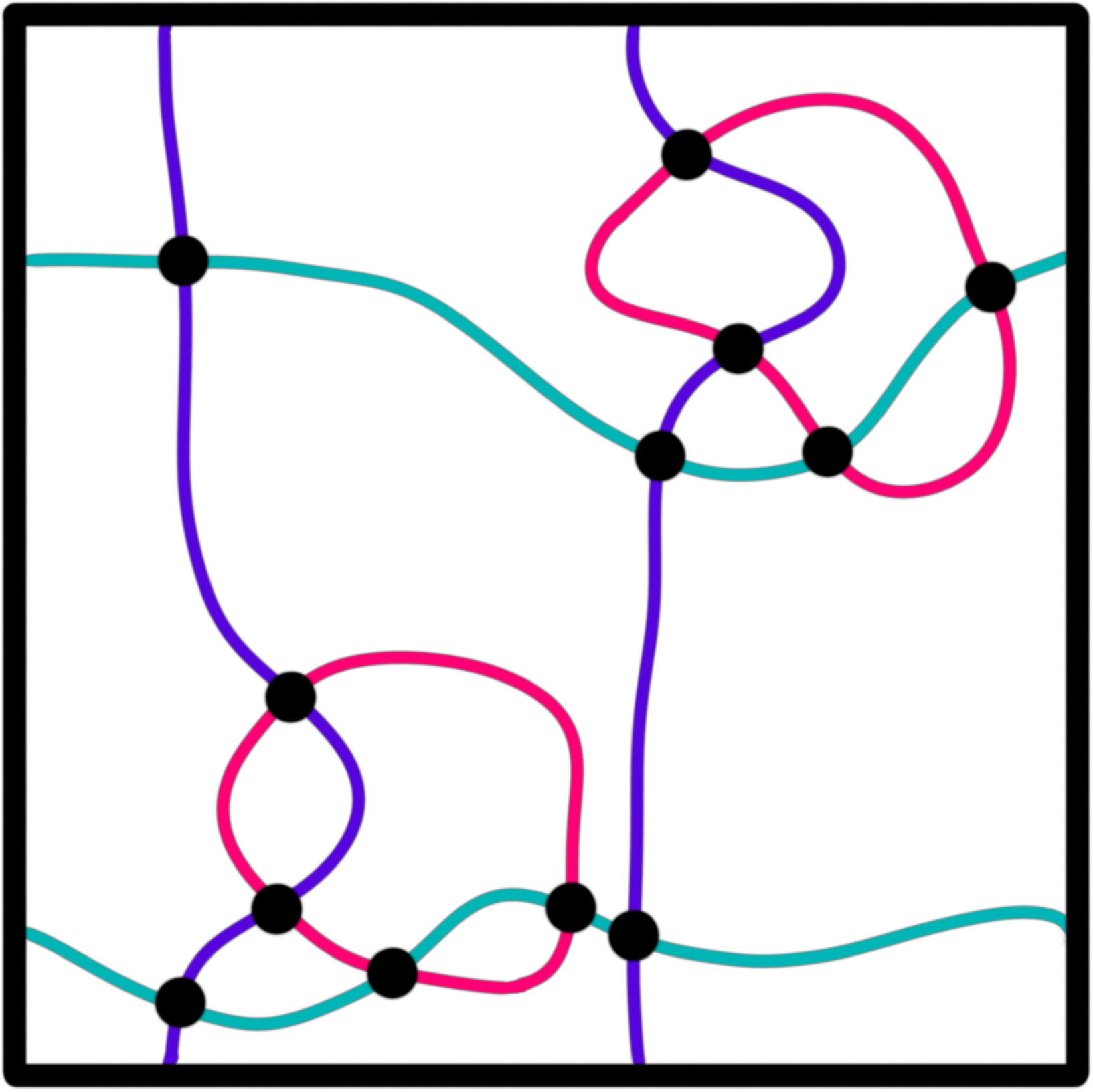}
        \caption{}
        \label{fig:homotopy_2}
    \end{subfigure}
    \begin{subfigure}[b]{0.19\textwidth}
        \centering
        \includegraphics[width=0.975\textwidth]{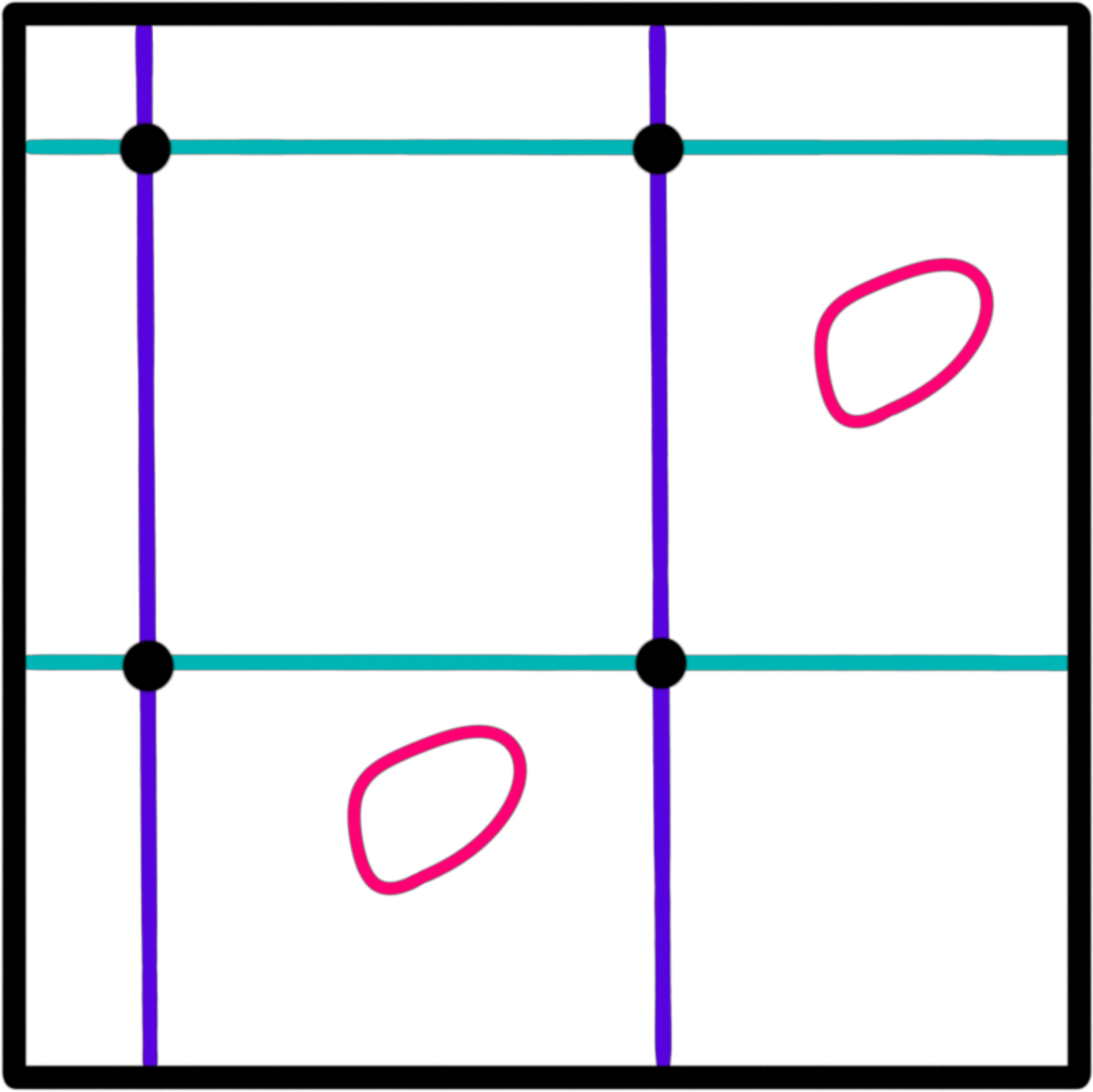}
        \caption{}
        \label{fig:homotopy_3}
    \end{subfigure}
    \begin{subfigure}[b]{0.19\textwidth}
        \centering
        \includegraphics[width=0.975\textwidth]{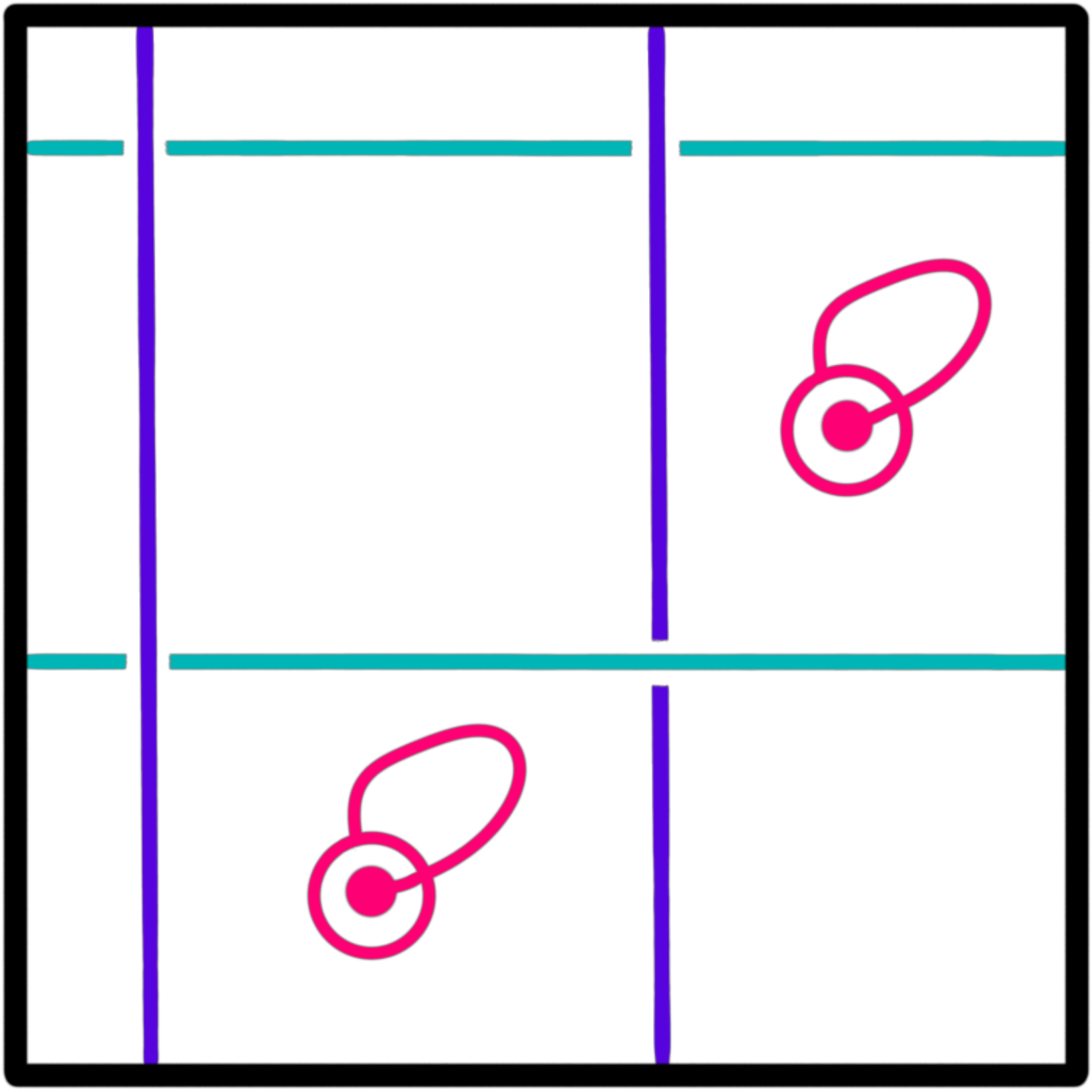}
        \caption{}
        \label{fig:homotopy_4}
    \end{subfigure}
    \caption{Relationship between crossings in diagrams of ground states and minimum number of intersections among free homotopy classes of curves. A diagram of the 3-periodic tangle of Fig.~\ref{fig:pi_plus_u2_uc_xyz_and_tridia} is displayed in (a). In (b) the over-under information is omitted as it is no longer of interest. By minimising the number of intersections of the curves of (b) under homotopy deformations, one may obtain the configuration of (c) among other possibilities. We display in (d) a diagram of $\Pi^{+}$ where its crossings coincide with the intersections of the curves of (c). Since all crossings in each diagram of the tridiagram of $\Pi^{+}$ correspond to intersections that cannot be eliminated under homotopy deformations, it is a ground state.}
    \label{fig:homotopy}
\end{figure}

Fig.~\ref{fig:homotopy} shows an example of these necessary crossings in a given tridiagram. In Fig.~\ref{fig:homotopy_1}, we consider the first diagram of the tridiagram displayed in Fig.~\ref{fig:pi_plus_u2_tridia}, from which the over-under information is omitted in Fig.~\ref{fig:homotopy_2} as it is no longer relevant in the context of homotopy. A possible configuration displaying the minimum number of intersections is given in Fig.~\ref{fig:homotopy_3}; loops can be placed in any region of the square while avoiding intersections, but perpendicular lines must intersect, giving a total of four intersections. The crossings in each diagram of $\Pi^{+}$'s tridiagram shown in Fig.~\ref{fig:pi_plus_tridia} coincide with these intersections, confirming that it is a ground state.

For a more precise computation of the minimum number of intersections between two simple curves, one need only consider the representatives $(k_1,k_2)$ and $(k_1',k_2')$ of the curves in the fundamental group $\pi_1(\mathbb{T}^2) \cong \mathbb{Z}^2$ of the 2-torus $\mathbb{T}^2$, and compute $|k_1k_2' - k_2k_1'|$ \cite{primer_mapping_cl_group_chap1}. The minimum number of intersections among all the curves is then obtained by summing the minimum number of intersections for each pair of curves. In our example in Fig.~\ref{fig:homotopy}, two curves are represented by $(1,0)$, two by $(0,1)$, and the remaining two by $(0,0)$. In total, they yield four intersections.

The computation of self-intersections of non-simple curves, as well as intersections between two non-simple curves, is also known \cite{GU2022108205}, although no simple formula is currently available.

It is important to note that, for a 3-periodic tangle to be a ground state, all the crossings in each diagram of its minimal tridiagram must coincide with the intersections that cannot be eliminated under homotopy deformations. Indeed, the diagram of the plain weave shown in Fig.~\ref{fig:plain_weave_dia_F} only has four crossings, which is the least possible for such an arrangement of curves (Compare with the diagram of \ref{fig:plain_weave_ground_state_dia_F}, which is that of a ground state). However, the plain weave is not a ground state as it cannot have the least number of crossings along three non-coplanar axes (Compare Fig.~\ref{fig:plain_weave_xyz_uc} and Fig.~\ref{fig:plain_weave_ground_state_xyz_uc}).

The example of the plain weave shows that, although the main focus of the paper is on 3-periodic tangles, the methods outlined here also apply to structures of lower periodicity. This is due to the use of additional information obtained by working with tridiagrams, in contrast to the traditional approach in knot-theoretic analysis, which typically relies on a single diagram, as in \cite{diamantis2023equivalencedoublyperiodictangles}, for example, for 2-periodic tangles.

Theorem~\ref{thm:homotopy_ground_states} provides a sufficient criterion for a 3-periodic tangle to be a ground state with respect to a given unit cell. Such a criterion is not yet available for 3-periodic networks \cite{andriamanalina2026measuringentanglementcomplexity3periodic}, which highlights the distinctive contribution of the present work.

\bigbreak

As suggested by the example of the $\Pi^{+}$ rod packing in Fig.~\ref{fig:homotopy}, the crystallographic rod packings are typical examples of structures that satisfy the homotopy condition presented in Theorem~\ref{thm:homotopy_ground_states}, and thus, are typical ground states. Some other possibilities are helices or loops.

\begin{figure}[hbtp]
    \centering
    \begin{subfigure}[b]{0.19\textwidth}
        \centering
        \includegraphics[width=0.975\textwidth]{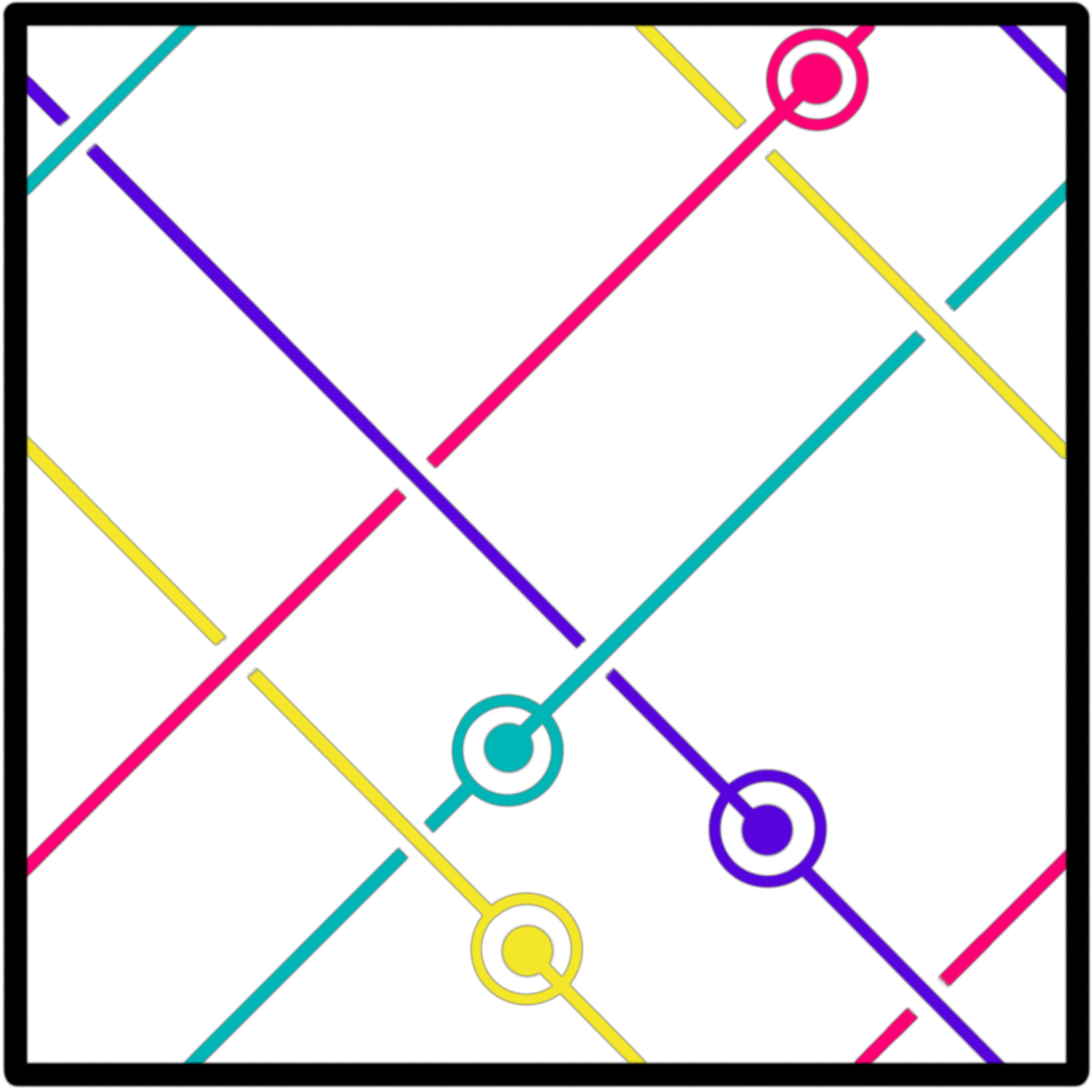}
        \caption{}
        \label{fig:sigma_plus_dia}
    \end{subfigure}
    \begin{subfigure}[b]{0.19\textwidth}
        \centering
        \includegraphics[width=0.975\textwidth]{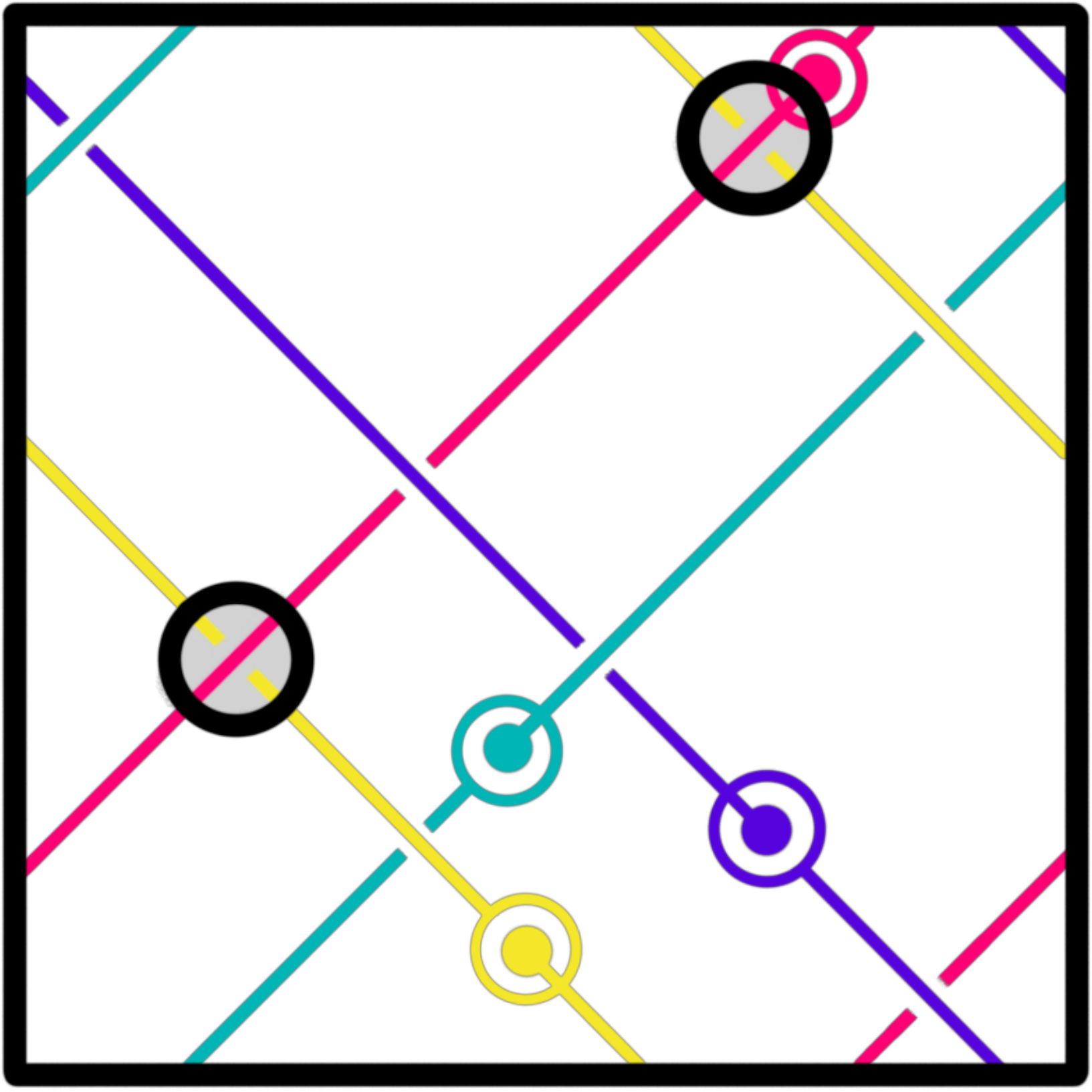}
        \caption{}
        \label{fig:sigma_plus_dia_crossing_changes}
    \end{subfigure}
    \begin{subfigure}[b]{0.19\textwidth}
        \centering
        \includegraphics[width=0.975\textwidth]{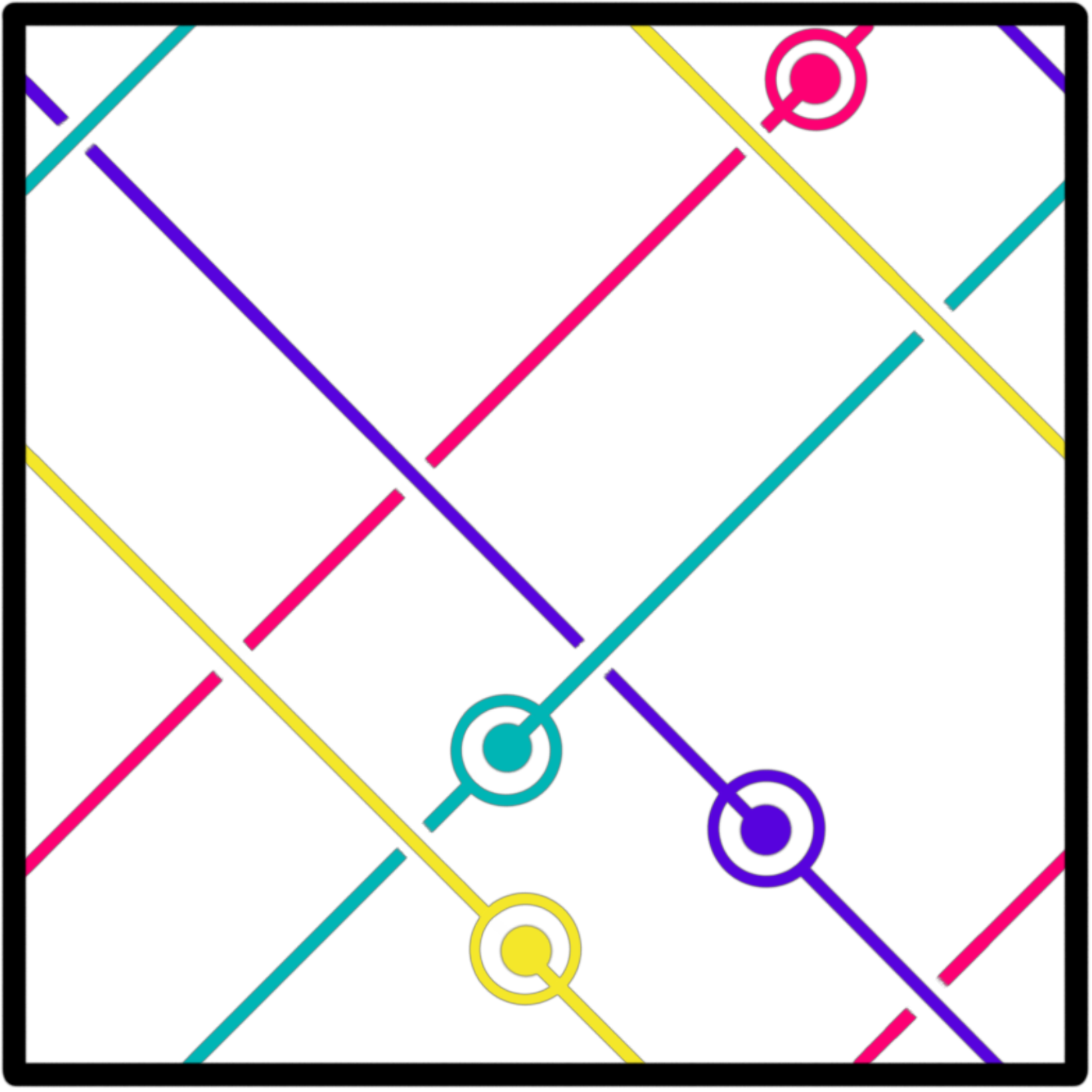}
        \caption{}
        \label{fig:omega_minus_dia}
    \end{subfigure}
    \begin{subfigure}[b]{0.19\textwidth}
        \centering
        \includegraphics[width=0.975\textwidth]{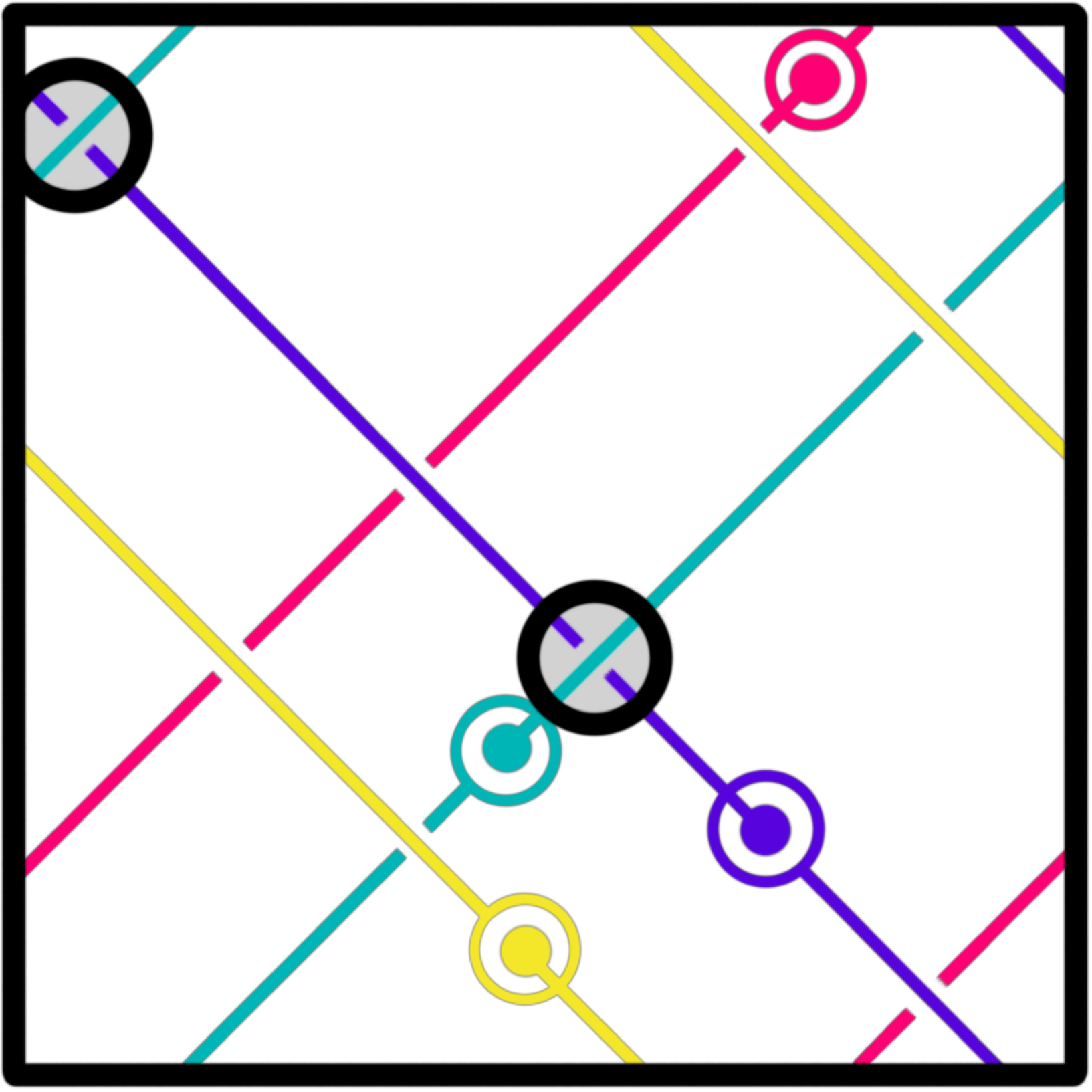}
        \caption{}
        \label{fig:omega_minus_dia_crossing_changes}
    \end{subfigure}
    \begin{subfigure}[b]{0.19\textwidth}
        \centering
        \includegraphics[width=0.975\textwidth]{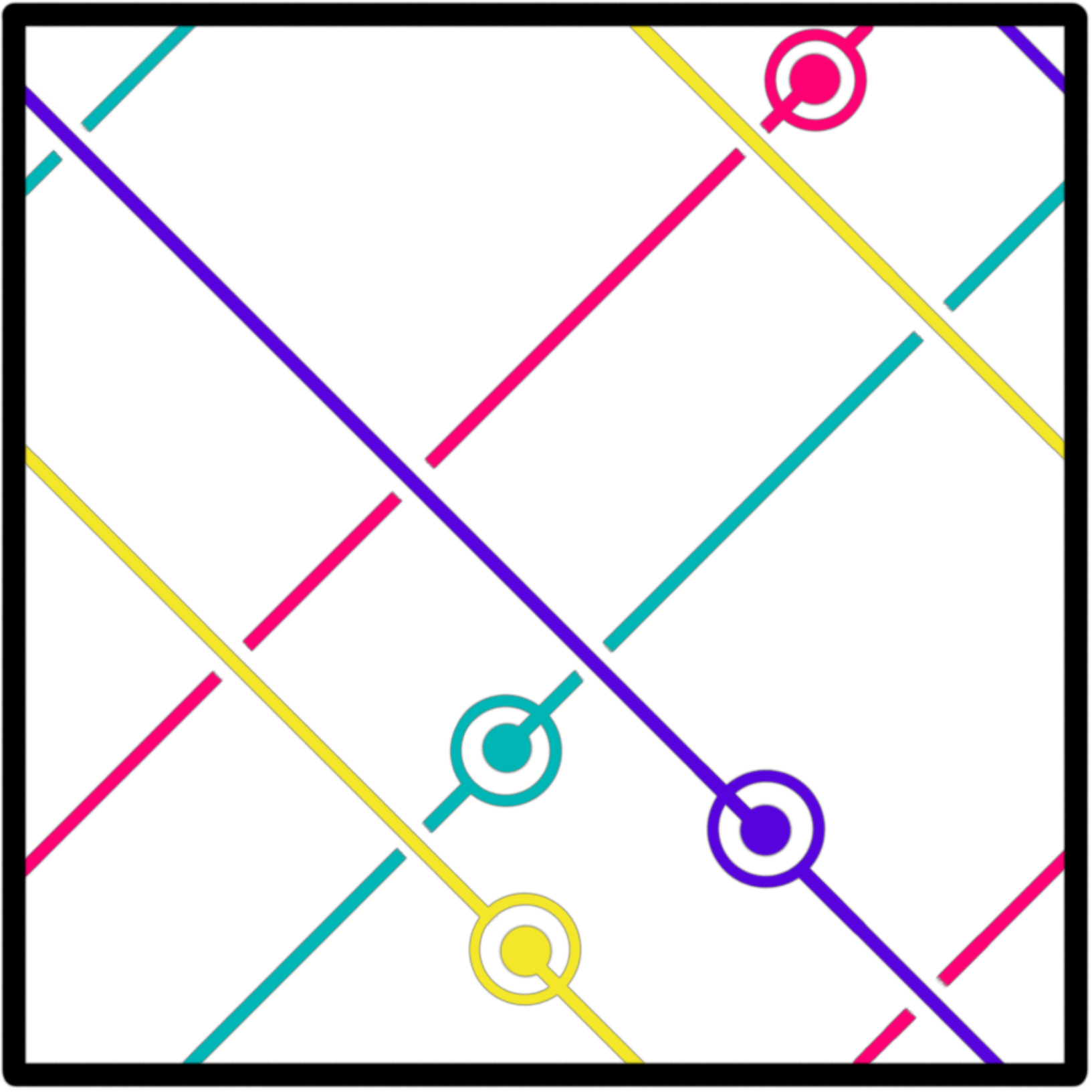}
        \caption{}
        \label{fig:gamma_dia}
    \end{subfigure}
    \caption{The three cubic rod packings with four axes connected by crossing changes. Diagrams of $\Sigma^{+}$, $\Omega^{-}$, and $\Gamma$ are respectively displayed in (a), (c), and (e). In (b) and (d), the black circles highlight the crossings onto which the crossing changes are applied, establishing the connection between the rod packings. The diagrams displayed here are not obtained from the usual embeddings of the rod packings but from isotopic ones.}
    \label{fig:sigma_omega_gamma}
\end{figure}

Of the six invariant cubic rod packings shown in Fig.~\ref{fig:six_rod_packings}, the $\Sigma^{+}$, $\Omega^{+}$ and $\Gamma$ rod packings as well as their mirror images are all connected by crossing changes with respect to their usual unit cells. In other words, they belong to the same $\mathcal{G}$-family. In Fig.~\ref{fig:sigma_plus_dia}, Fig.~\ref{fig:omega_minus_dia}, and Fig.~\ref{fig:gamma_dia}, we present diagrams representing $\Sigma^{+}$, $\Omega^{-}$ the mirror image of $\Omega^{+}$, and $\Gamma$, respectively. In Fig.~\ref{fig:sigma_plus_dia_crossing_changes} and Fig.~\ref{fig:omega_minus_dia_crossing_changes}, we highlight with black circles the crossings onto which the crossing changes are applied to transform one structure into another. Since $\Gamma$ is achiral, that is, its mirror image is itself, the connections to $\Sigma^{-}$ and $\Omega^{+}$ follow immediately. Note that the diagrams shown Fig.~\ref{fig:sigma_omega_gamma} are not obtained from the usual embeddings of the rod packings (displayed in Fig.~\ref{fig:six_rod_packings}), but rather from isotopic embeddings. However, it is with these diagrams that the connection between the three rod packings is most easily established.

In a similar fashion, $\Pi^{+}$ and its mirror image $\Pi^{-}$ are connected by crossing changes. However, $\Pi^{\star}$ and $\Sigma^{\star}$ are not connected to any of the other rod packings.

\section{Untangling number of 3-periodic tangles}\label{sec:untangling_number}
The definition of ground states given in Section \ref{sec:ground_states} allows us to define the untangling number of 3-periodic tangles, which can be used as a measure of entanglement complexity.

\bigbreak

The \emph{unknotting number} of a classical knot is the least number of crossing changes needed to transform the knot into the unknot \cite{Adams.book,Murasugi1996chap4}. By analogy, given a 3-periodic tangle $K$ and a unit cell $U$ of $K$, we define the \emph{untangling number of $K$ with respect to $U$}, denoted by $u(K,U)$, to be the least number of crossing changes needed to transform $U$ into a ground state belonging to the $\mathcal{G}$-family of $K$ with respect to $U$. Any ground state that realises the untangling number is called a \emph{nearest ground state of $K$ with respect to $U$}.

\bigbreak

In the following, we provide additional mathematical precision to the previous notions.

\begin{definition}\label{def:untangling_numb_U}
Let $K$ be a 3-periodic tangle and let $U$ be a unit cell of $K$. Consider all sequences $(U_{k})_{k=0,\dots,n}$  of unit cells of 3-periodic tangles belonging to $\mathcal{U}(K,U)$ satisfying:
\begin{itemize}
    \item $U_{0} = U$,
    \item $\forall k \in \lbrace 1, \dots,n\rbrace$,
    $U_{k}$ and $U_{k-1}$ are connected by one and only one crossing change,
    \item the 3-periodic tangle associated to $U_{n}$ belongs to $\mathcal{G}(K,U)$.
\end{itemize}
The \emph{untangling number of $K$ with respect to $U$}, denoted by $u(K,U)$, is the least integer $n$ over all such sequences. Any element of $\mathcal{G}(K,U)$ realising  $u(K,U)$ is called a \emph{nearest ground state of $K$ with respect to $U$}.
\end{definition}

\begin{remark}
    In Definition~\ref{def:untangling_numb_U}, we choose to enforce the condition that every pair of consecutive unit cells in a given sequence $(U_{k})_{k=0,\dots,n}$ is connected by one and only one crossing change, to ensure that the number $n$ indexing the sequence of unit cells corresponds to the total number of crossing changes.
\end{remark}

An immediate result obtained from Definition~\ref{def:untangling_numb_U} is that, the ground states are 3-periodic tangles with untangling number 0, and conversely, every 3-periodic tangle that has untangling number 0 with respect to a chosen unit cell is a ground state.

\bigbreak

We note that, although well defined, the untangling number is in practice difficult to compute, and only an upper bound can be effectively determined. This is similar to the limitation of the unknotting number, for which additional criteria are required to establish a lower bound. Fig.~\ref{fig:untangling_plain_weave} provides an example of the computation of such an upper bound for the plain weave. In particular, it shows that the untangling number, as described here, is well defined for structures with lower periodicity. In fact, for loops that do not intersect the faces of the unit cell, the value of the untangling number coincides with that of the classical unknotting number.

\bigbreak

As a second example, we compute an upper bound for the untangling number of the 3-periodic tangle whose unit cell is shown in Fig.~\ref{fig:pi_plus_u2_uc_xyz_and_tridia}. We project the unit cell to obtain a diagram, shown in Fig.~\ref{fig:untangling_pi_plus_4,4,5_1_usual_1}. Onto the latter, we apply a crossing change on the crossing highlighted by a black circle, which untangles some curve components as shown in Fig.~\ref{fig:untangling_pi_plus_4,4,5_1_usual_2} to Fig.~\ref{fig:untangling_pi_plus_4,4,5_1_usual_4_1st}. Onto the diagram shown in Fig.~\ref{fig:untangling_pi_plus_4,4,5_1_usual_4_1st}, a second crossing change is applied, which transforms the initial diagram into the diagram of Fig.~\ref{fig:untangling_pi_plus_4,4,5_1_usual_9_1st}. This final diagram represents a unit cell of the $\Pi^{+}$ rod packing, which we know to be a ground state. Two crossing changes were needed to reach a ground state, which means that the untangling number is at most 2.

\begin{figure}[hbtp]
    \centering
    \begin{subfigure}[b]{0.19\textwidth}
        \centering
        \includegraphics[width=0.975\textwidth]{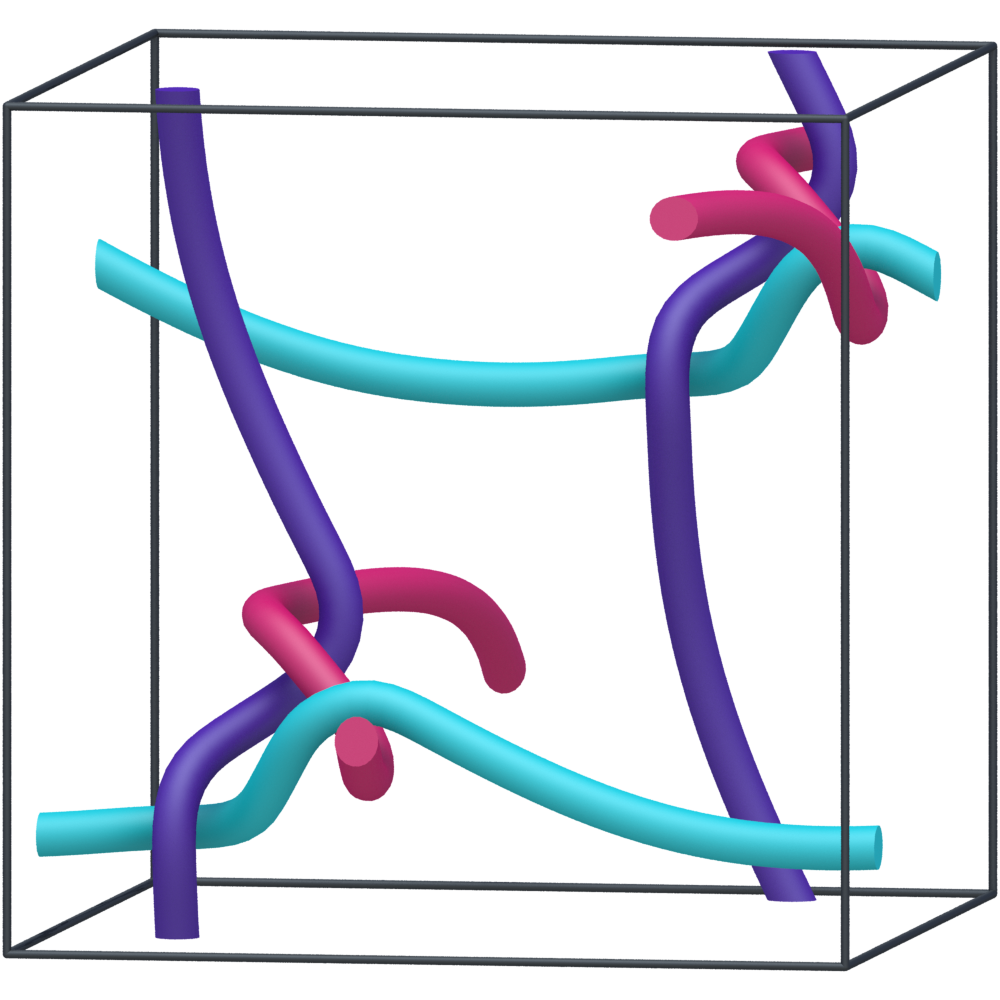}
        \caption{}
        \label{fig:pi_plus_u2_uc}
    \end{subfigure}
    \begin{subfigure}[b]{0.19\textwidth}
        \centering
        \includegraphics[width=0.975\textwidth]{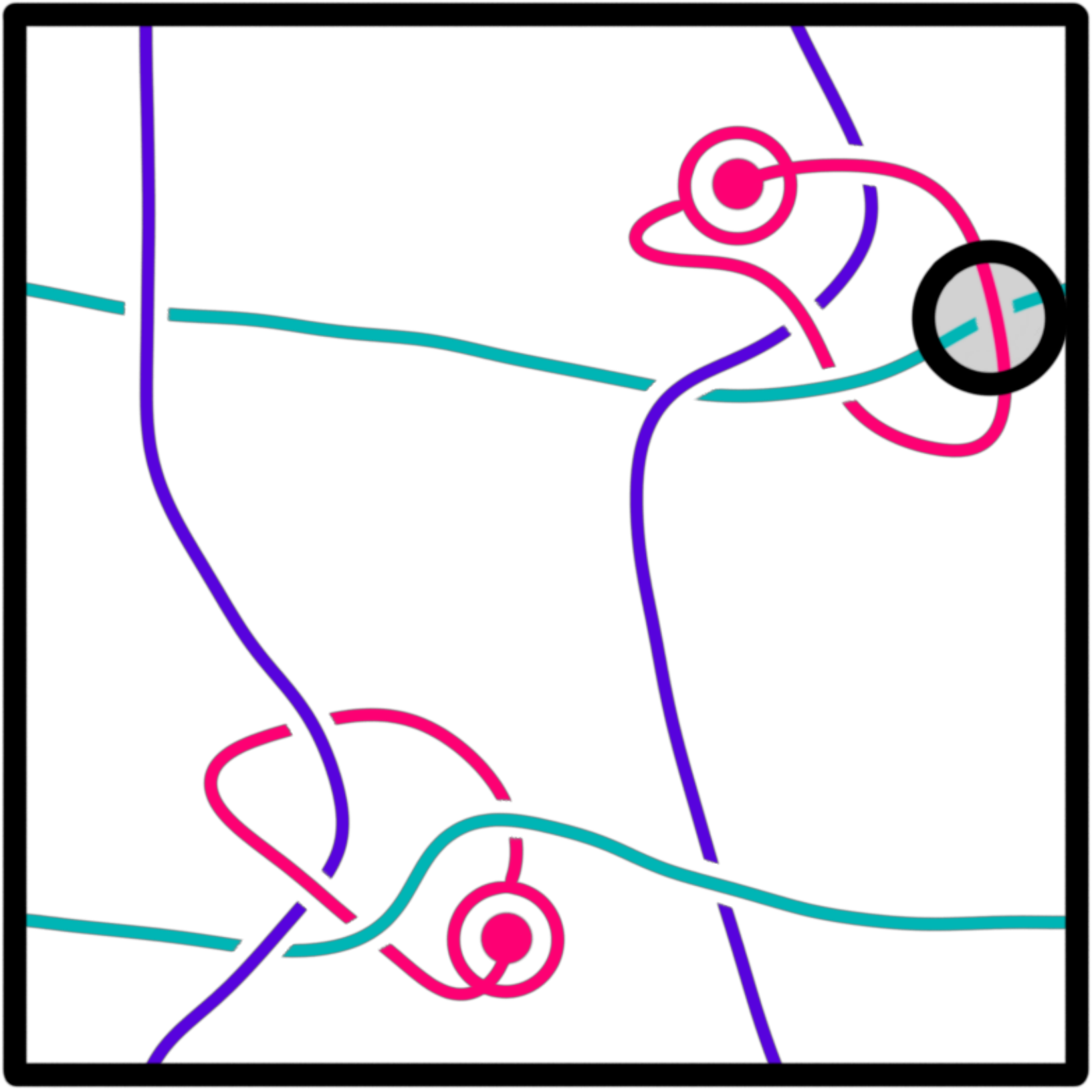}
        \caption{}
        \label{fig:untangling_pi_plus_4,4,5_1_usual_1}
    \end{subfigure}
    \begin{subfigure}[b]{0.19\textwidth}
        \centering
        \includegraphics[width=0.975\textwidth]{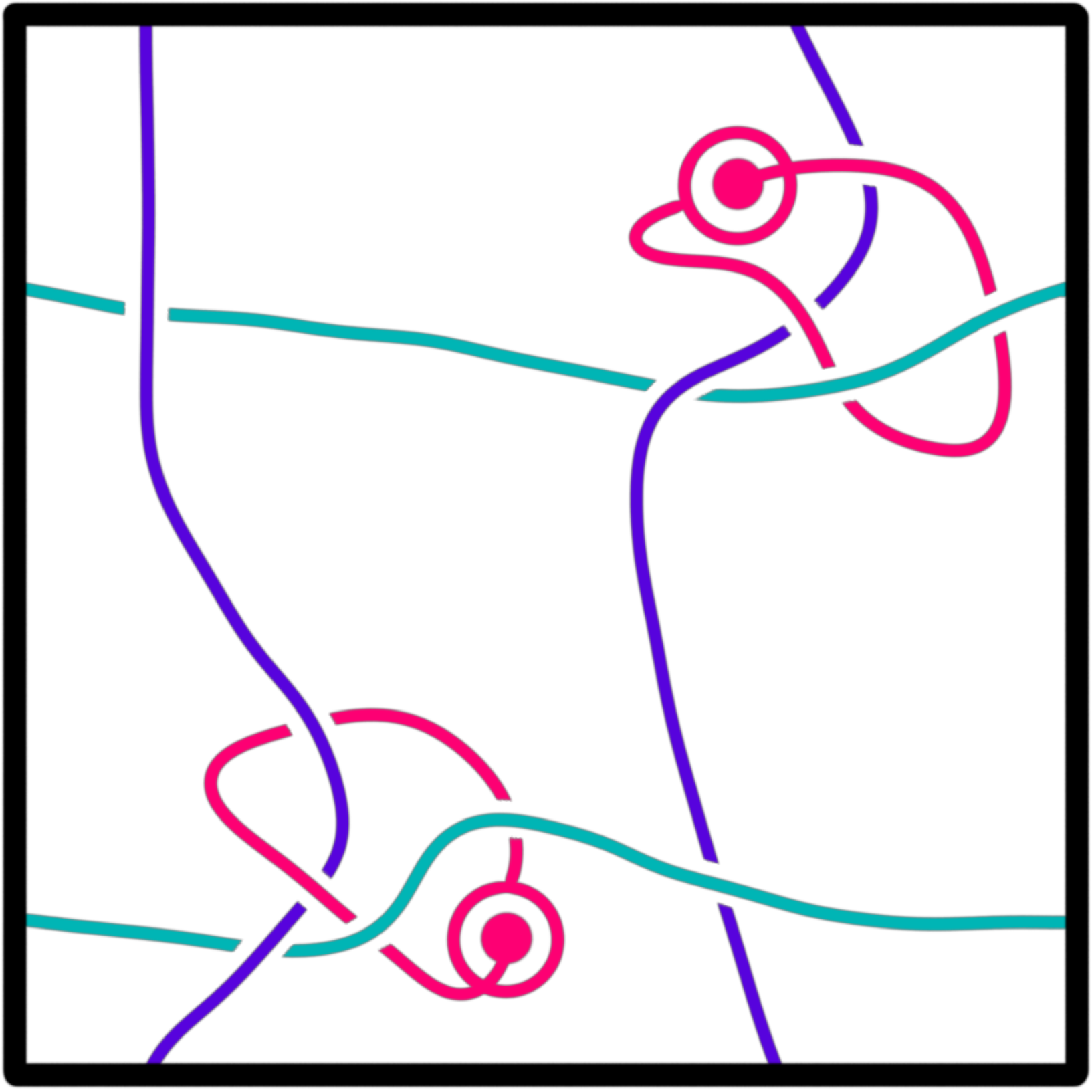}
        \caption{}
        \label{fig:untangling_pi_plus_4,4,5_1_usual_2}
    \end{subfigure}
    \begin{subfigure}[b]{0.19\textwidth}
        \centering
        \includegraphics[width=0.975\textwidth]{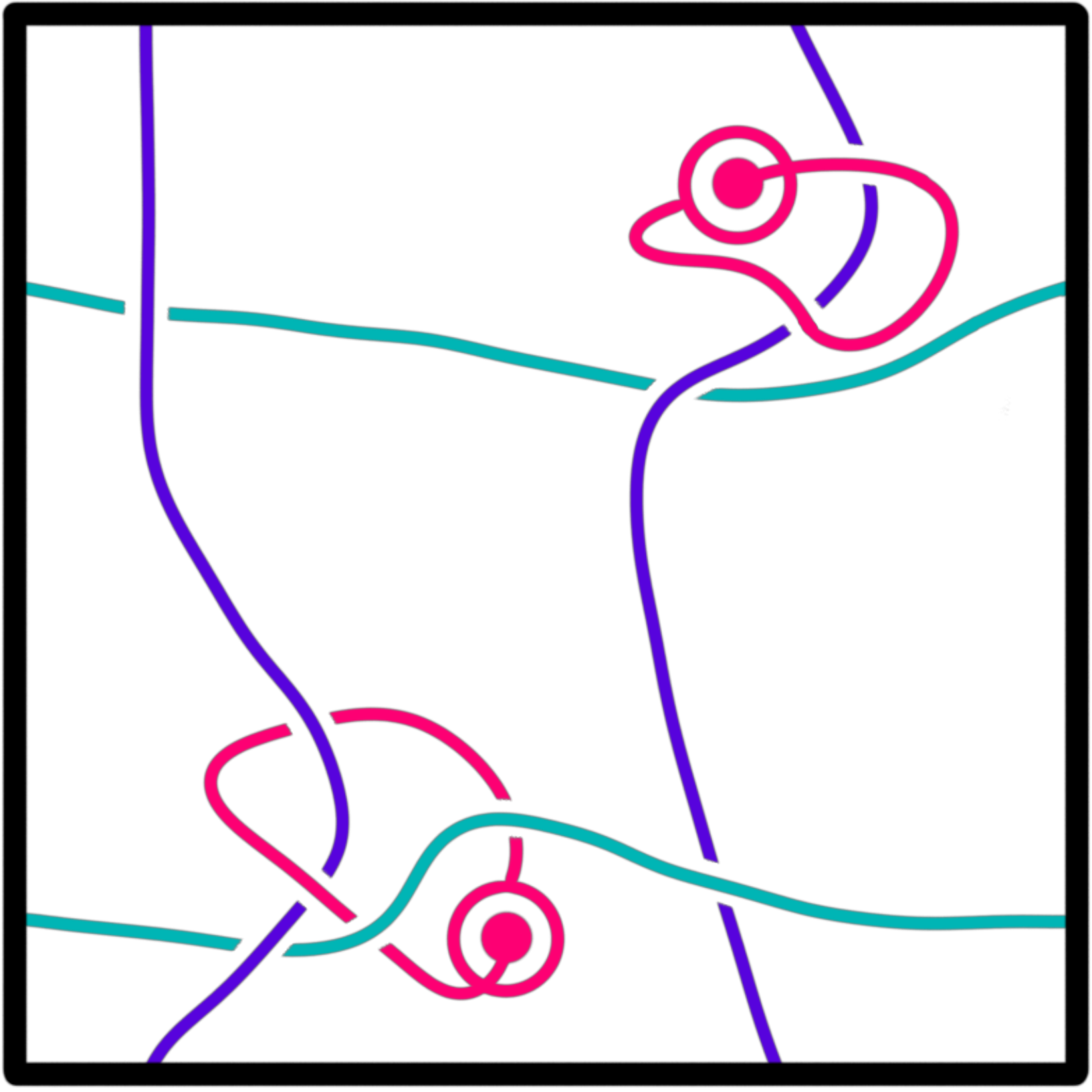}
        \caption{}
        \label{fig:untangling_pi_plus_4,4,5_1_usual_3}
    \end{subfigure}
    \begin{subfigure}[b]{0.19\textwidth}
        \centering
        \includegraphics[width=0.975\textwidth]{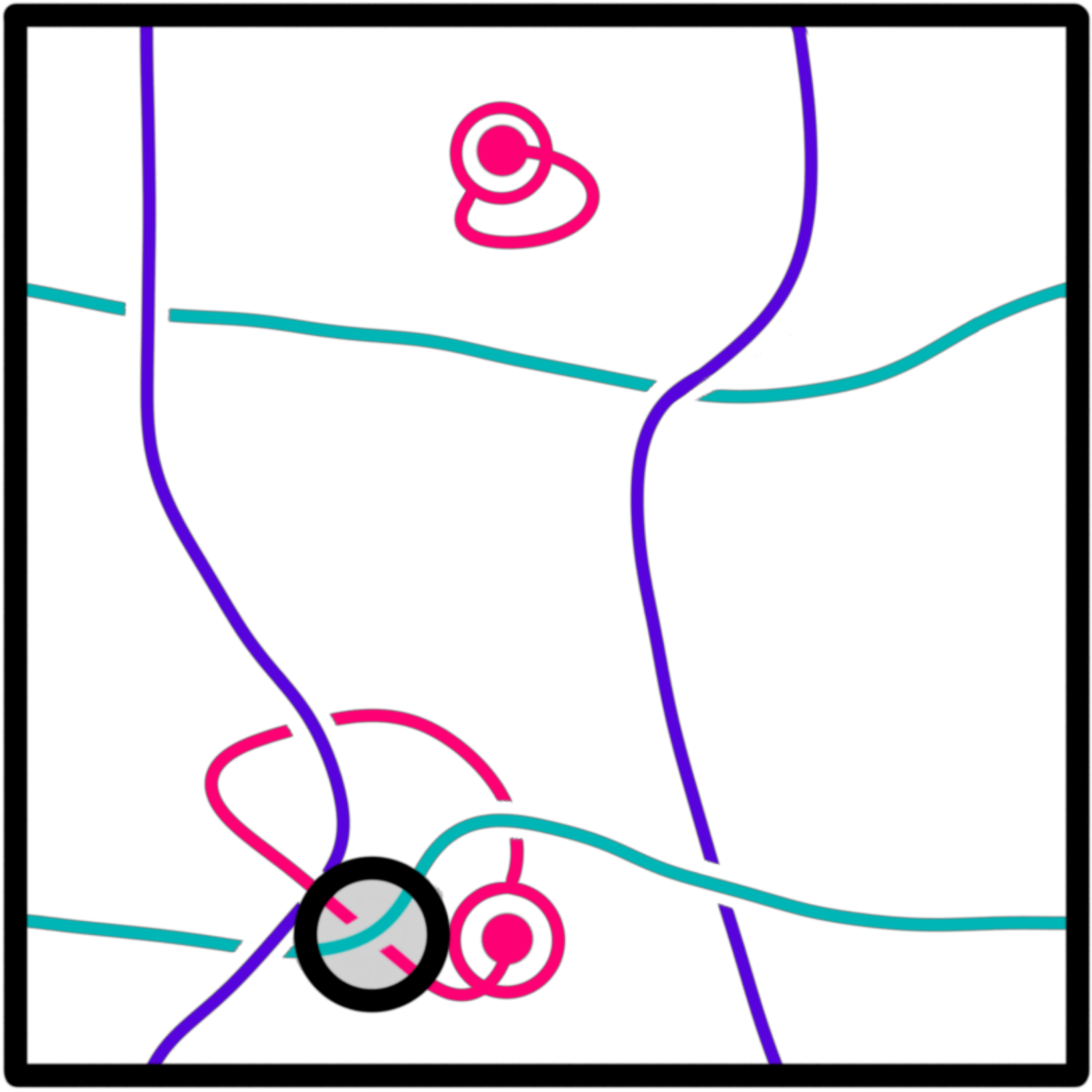}
        \caption{}
        \label{fig:untangling_pi_plus_4,4,5_1_usual_4_1st}
    \end{subfigure}

    \vskip\baselineskip
    
    \begin{subfigure}[b]{0.19\textwidth}
        \centering
        \includegraphics[width=0.975\textwidth]{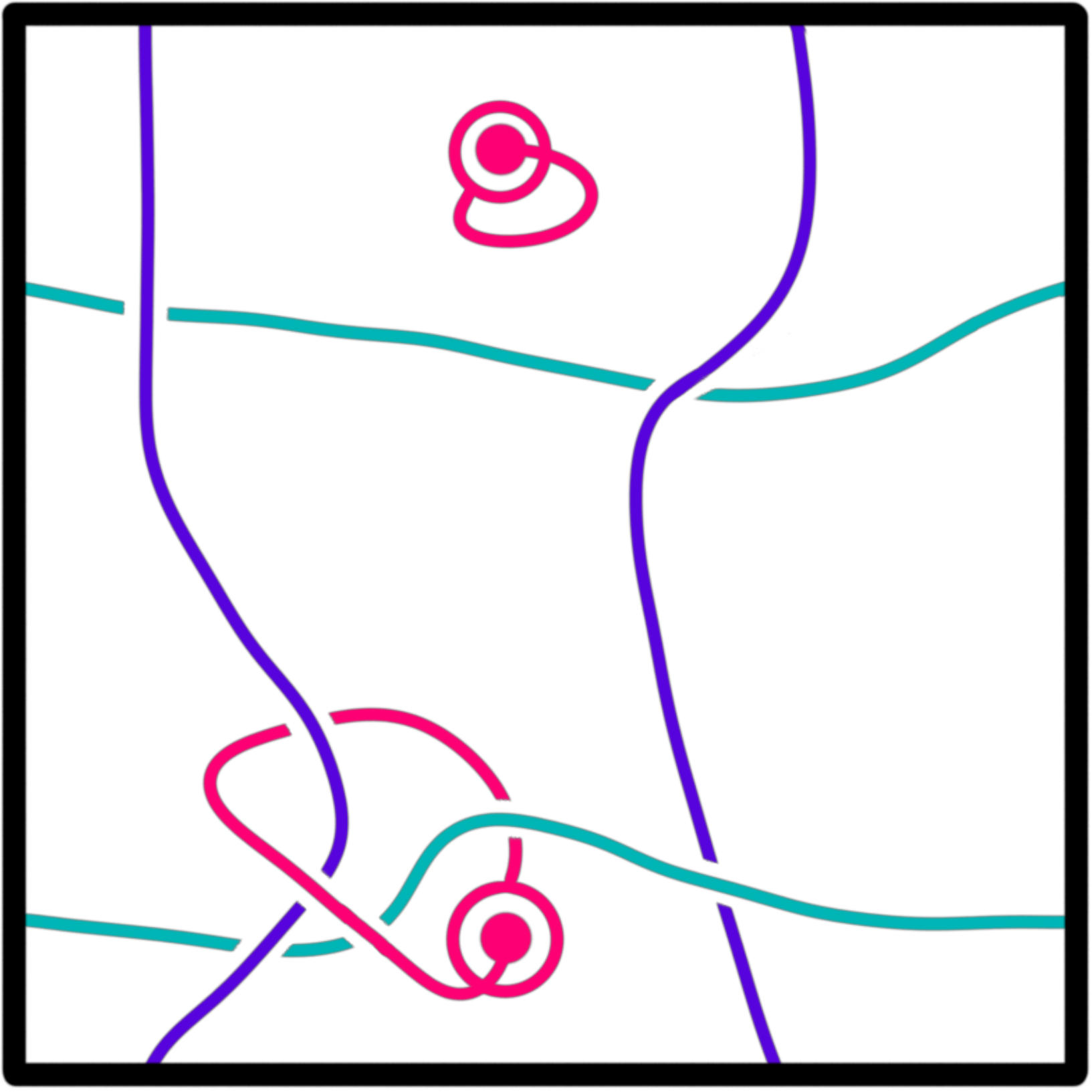}
        \caption{}
        \label{fig:untangling_pi_plus_4,4,5_1_usual_5_1st}
    \end{subfigure}
    \begin{subfigure}[b]{0.19\textwidth}
        \centering
        \includegraphics[width=0.975\textwidth]{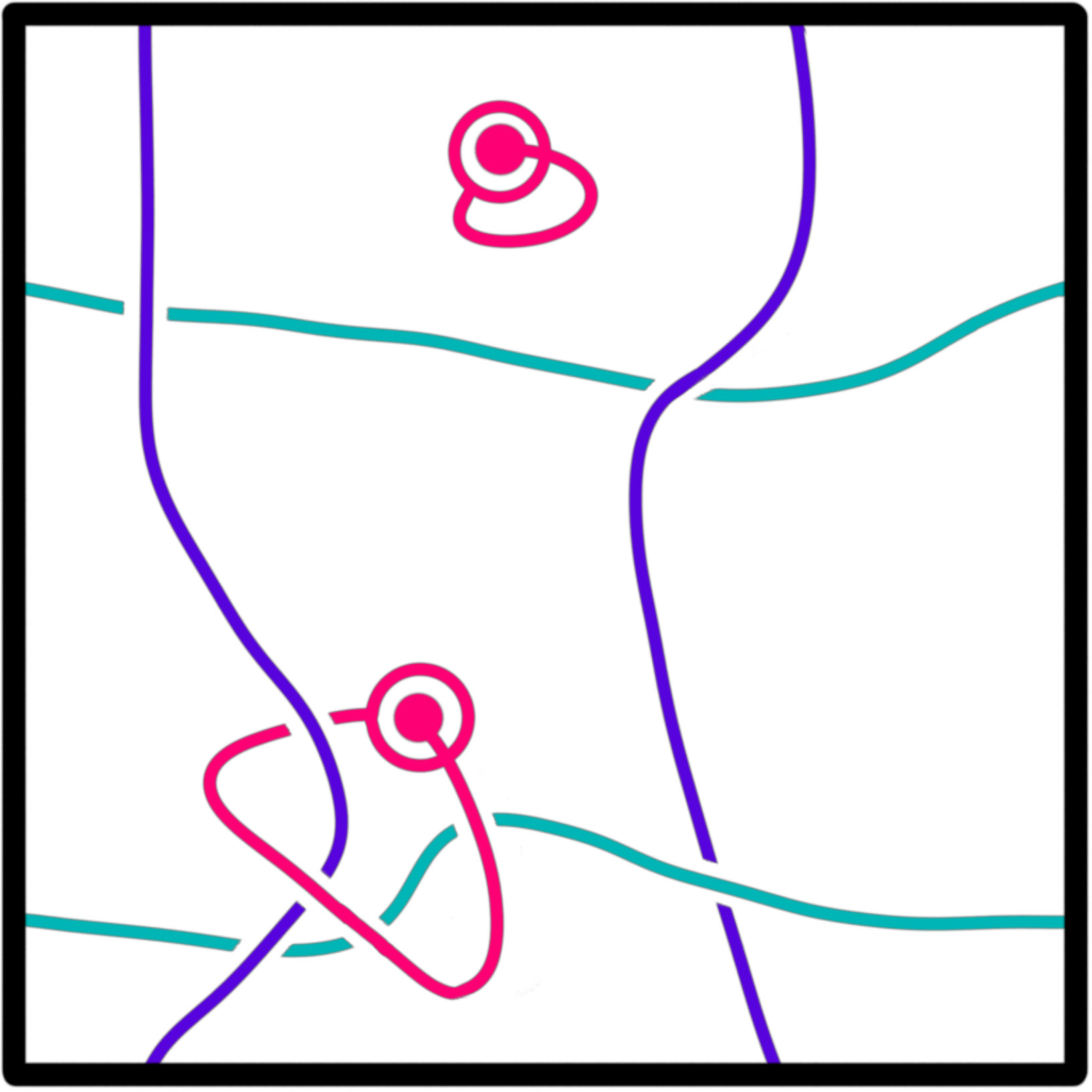}
        \caption{}
        \label{fig:untangling_pi_plus_4,4,5_1_usual_6_1st}
    \end{subfigure}
    \begin{subfigure}[b]{0.19\textwidth}
        \centering
        \includegraphics[width=0.975\textwidth]{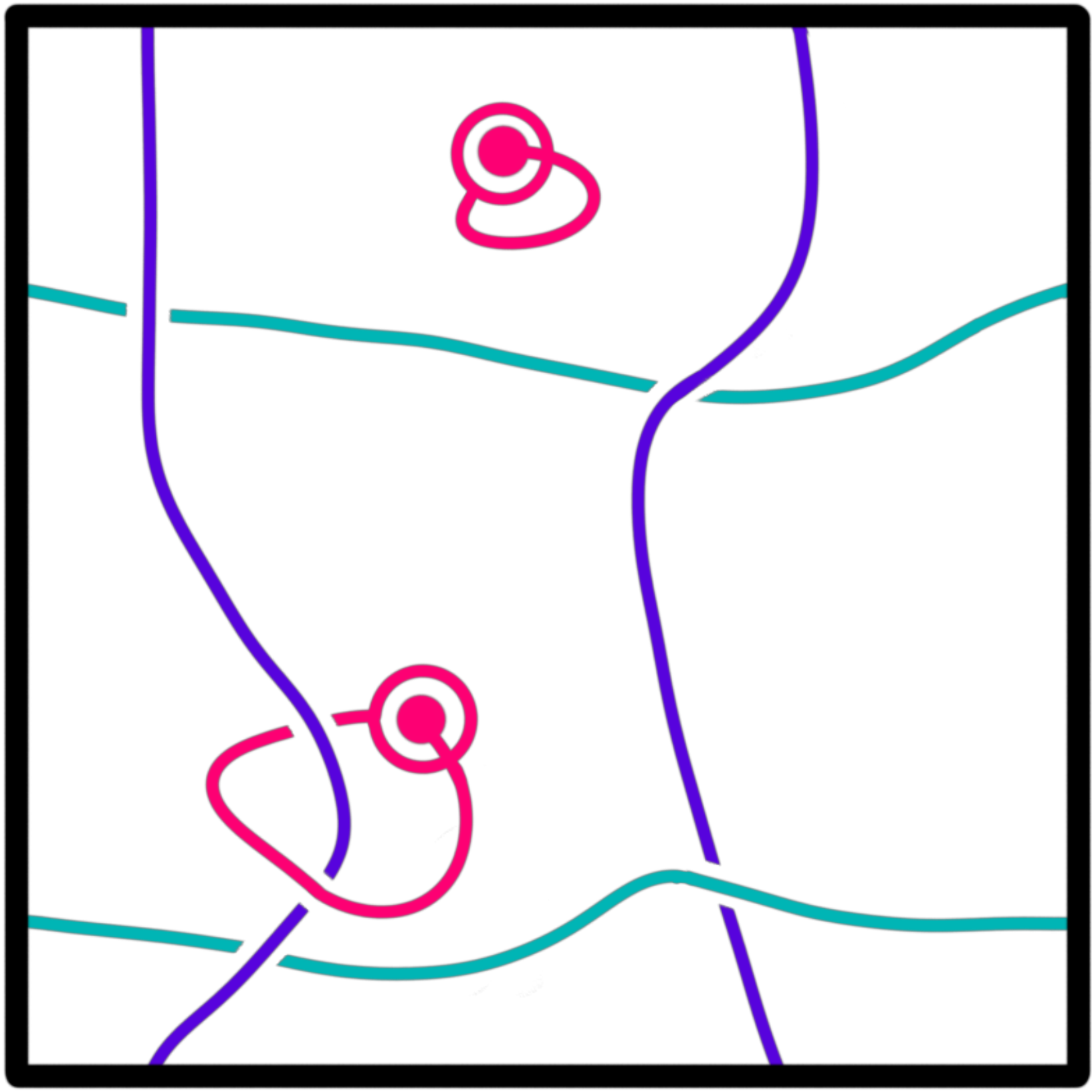}
        \caption{}
        \label{fig:untangling_pi_plus_4,4,5_1_usual_7_1st}
    \end{subfigure}
    \begin{subfigure}[b]{0.19\textwidth}
        \centering
        \includegraphics[width=0.975\textwidth]{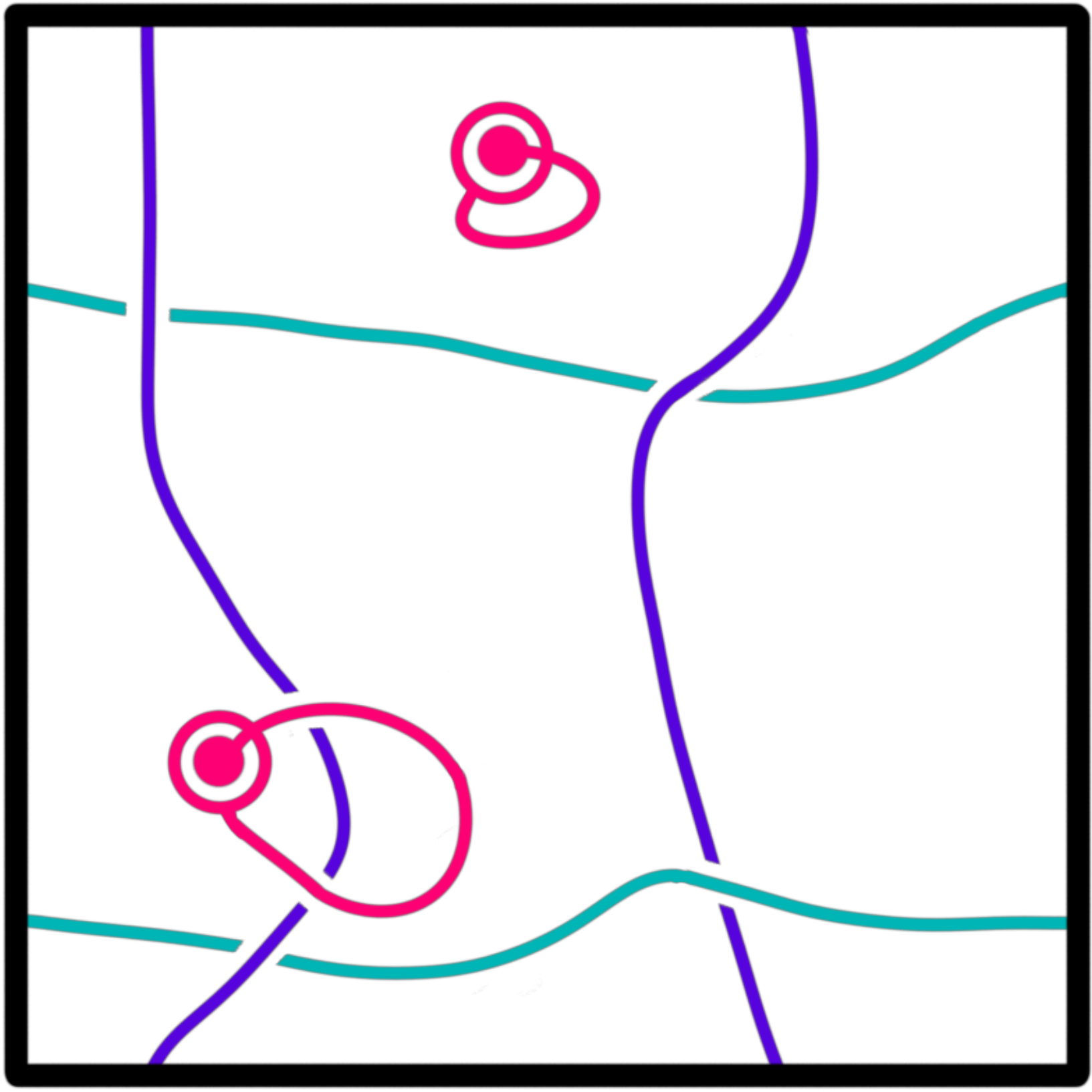}
        \caption{}
        \label{fig:untangling_pi_plus_4,4,5_1_usual_8_1st}
    \end{subfigure}
    \begin{subfigure}[b]{0.19\textwidth}
        \centering
        \includegraphics[width=0.975\textwidth]{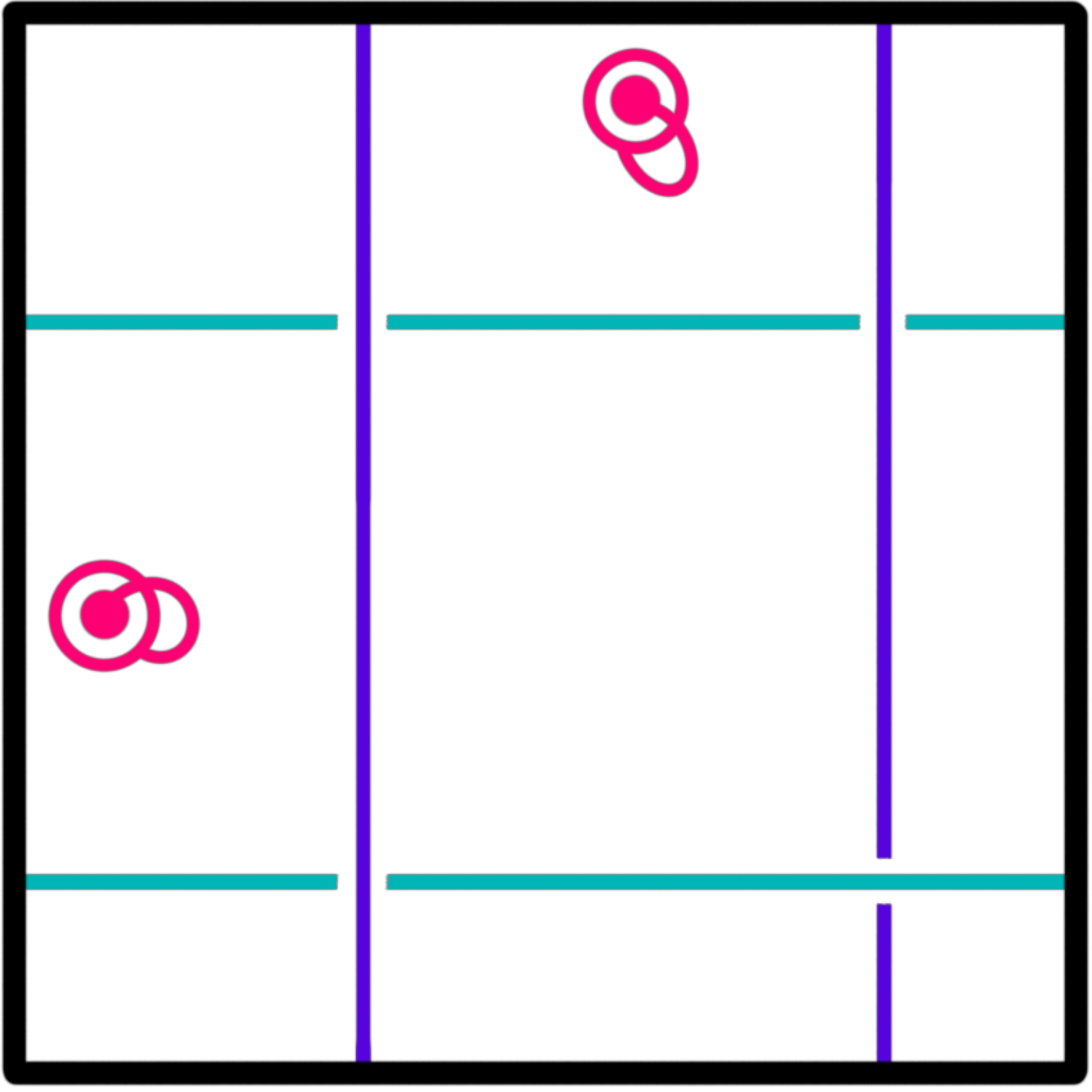}
        \caption{}
        \label{fig:untangling_pi_plus_4,4,5_1_usual_9_1st}
    \end{subfigure}
    \caption{Computing an upper bound for the untangling number of the 3-periodic tangle whose unit cell is shown in Fig.~\ref{fig:pi_plus_u2_uc_xyz_and_tridia}. Each diagram is obtained from the one before by either applying a crossing change on the crossing highlighted by a black circle, or by deforming the curves with a sequence of $R$-moves. The final diagram shown in (j) represents the $\Pi^{+}$ rod packing. Two crossing changes were needed to reach a ground state, which means that the untangling number is at most 2.}
    \label{fig:untangling_pi_plus_4,4,5_1_usual_first_time}
\end{figure}

\bigbreak

\begin{figure}[hbtp]
    \centering
    \begin{subfigure}[b]{0.19\textwidth}
        \centering
        \includegraphics[width=0.975\textwidth]{untangling_pi_plus_4,4,5_1_usual_1.pdf}
        \caption{}
        \label{fig:untangling_pi_plus_4,4,5_1_usual_1_2nd}
    \end{subfigure}
    \begin{subfigure}[b]{0.19\textwidth}
        \centering
        \includegraphics[width=0.975\textwidth]{untangling_pi_plus_4,4,5_1_usual_2.pdf}
        \caption{}
        \label{fig:untangling_pi_plus_4,4,5_1_usual_2_2nd}
    \end{subfigure}
    \begin{subfigure}[b]{0.19\textwidth}
        \centering
        \includegraphics[width=0.975\textwidth]{untangling_pi_plus_4,4,5_1_usual_3.pdf}
        \caption{}
        \label{fig:untangling_pi_plus_4,4,5_1_usual_3_2nd}
    \end{subfigure}
    \begin{subfigure}[b]{0.19\textwidth}
        \centering
        \includegraphics[width=0.975\textwidth]{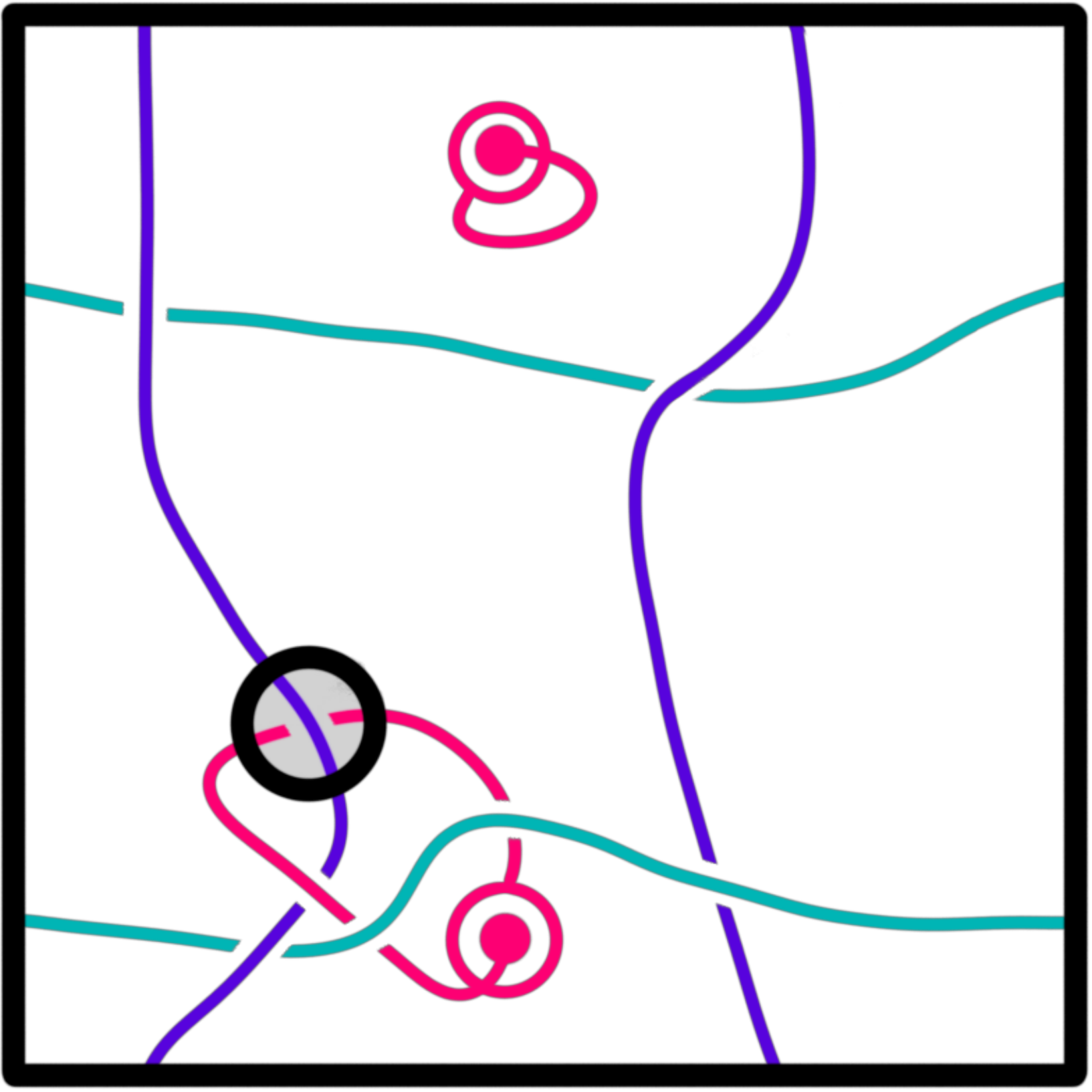}
        \caption{}
        \label{fig:untangling_pi_plus_4,4,5_1_usual_4_2nd}
    \end{subfigure}    
    \begin{subfigure}[b]{0.19\textwidth}
        \centering
        \includegraphics[width=0.975\textwidth]{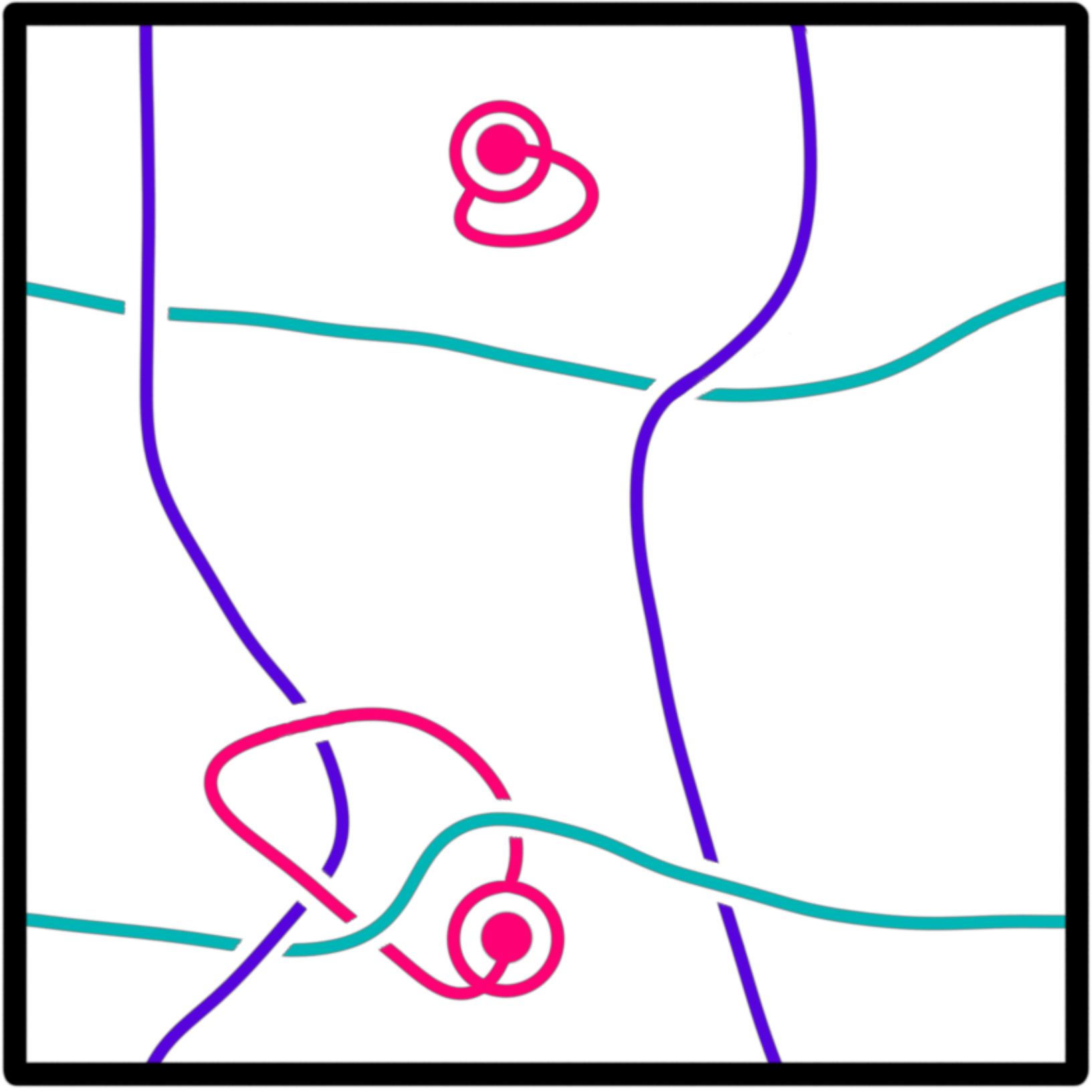}
        \caption{}
        \label{fig:untangling_pi_plus_4,4,5_1_usual_5_2nd}
    \end{subfigure}

    \vskip\baselineskip
    
    \begin{subfigure}[b]{0.19\textwidth}
        \centering
        \includegraphics[width=0.975\textwidth]{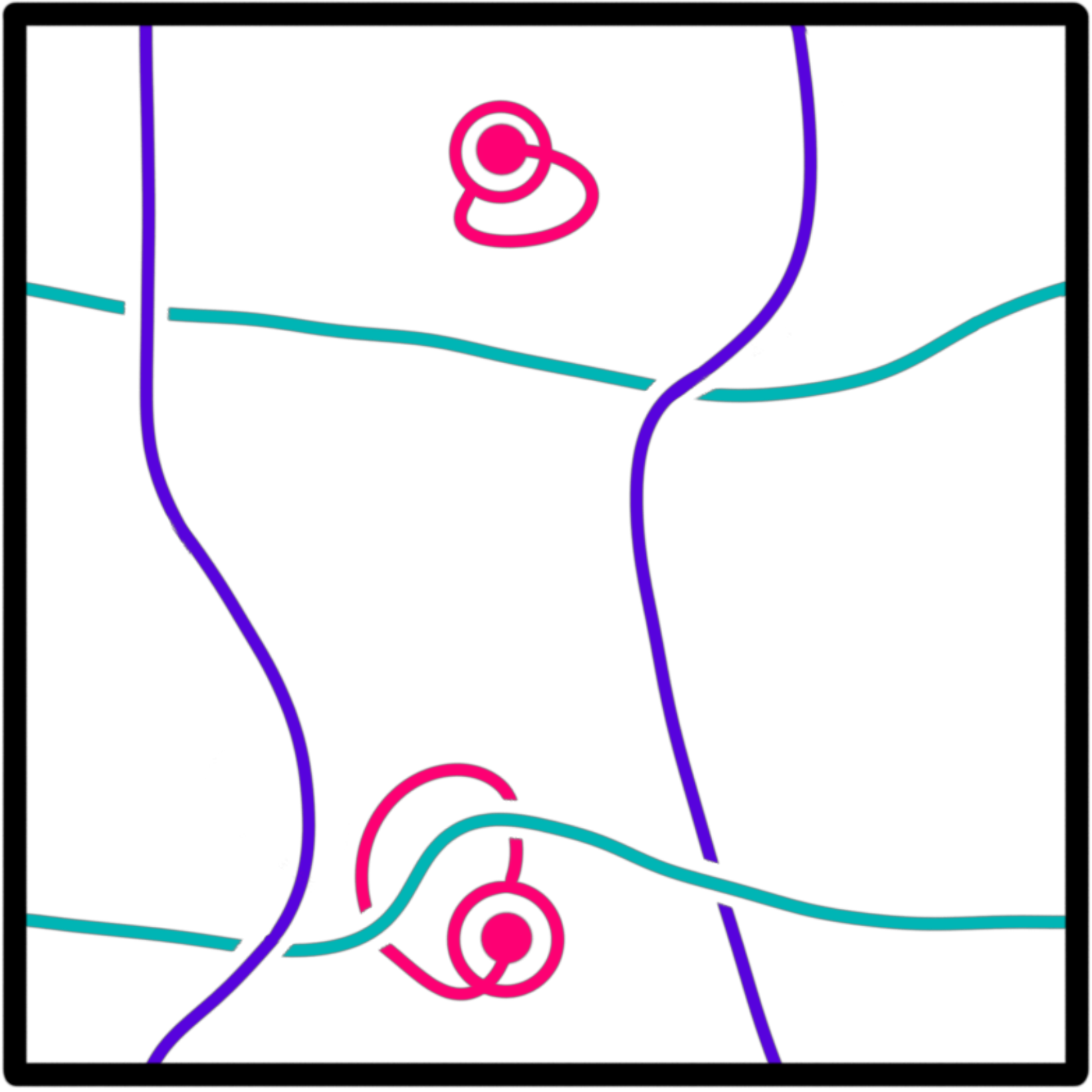}
        \caption{}
        \label{fig:untangling_pi_plus_4,4,5_1_usual_6_2nd}
    \end{subfigure}
    \begin{subfigure}[b]{0.19\textwidth}
        \centering
        \includegraphics[width=0.975\textwidth]{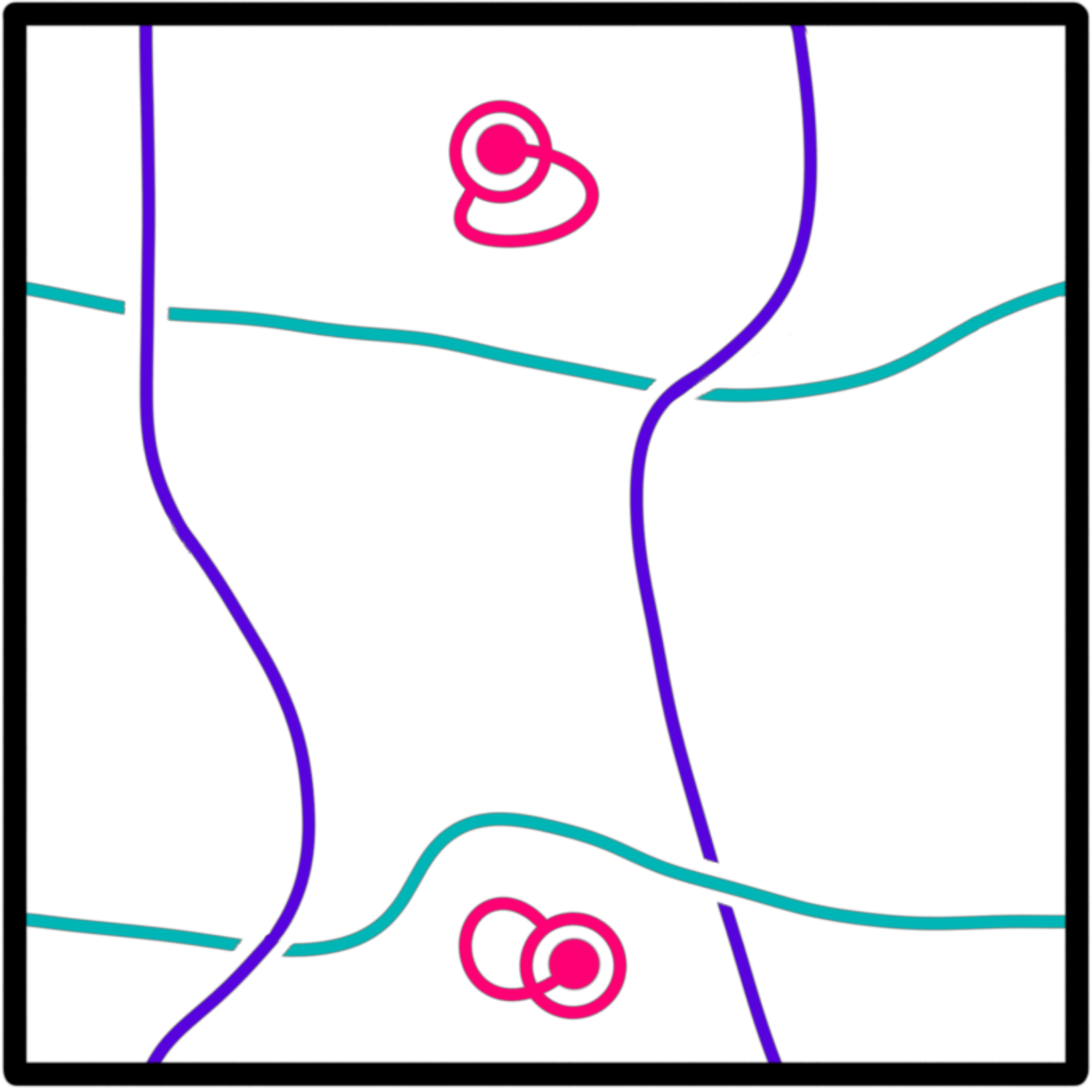}
        \caption{}
        \label{fig:untangling_pi_plus_4,4,5_1_usual_7_2nd}
    \end{subfigure}
    \begin{subfigure}[b]{0.19\textwidth}
        \centering
        \includegraphics[width=0.975\textwidth]{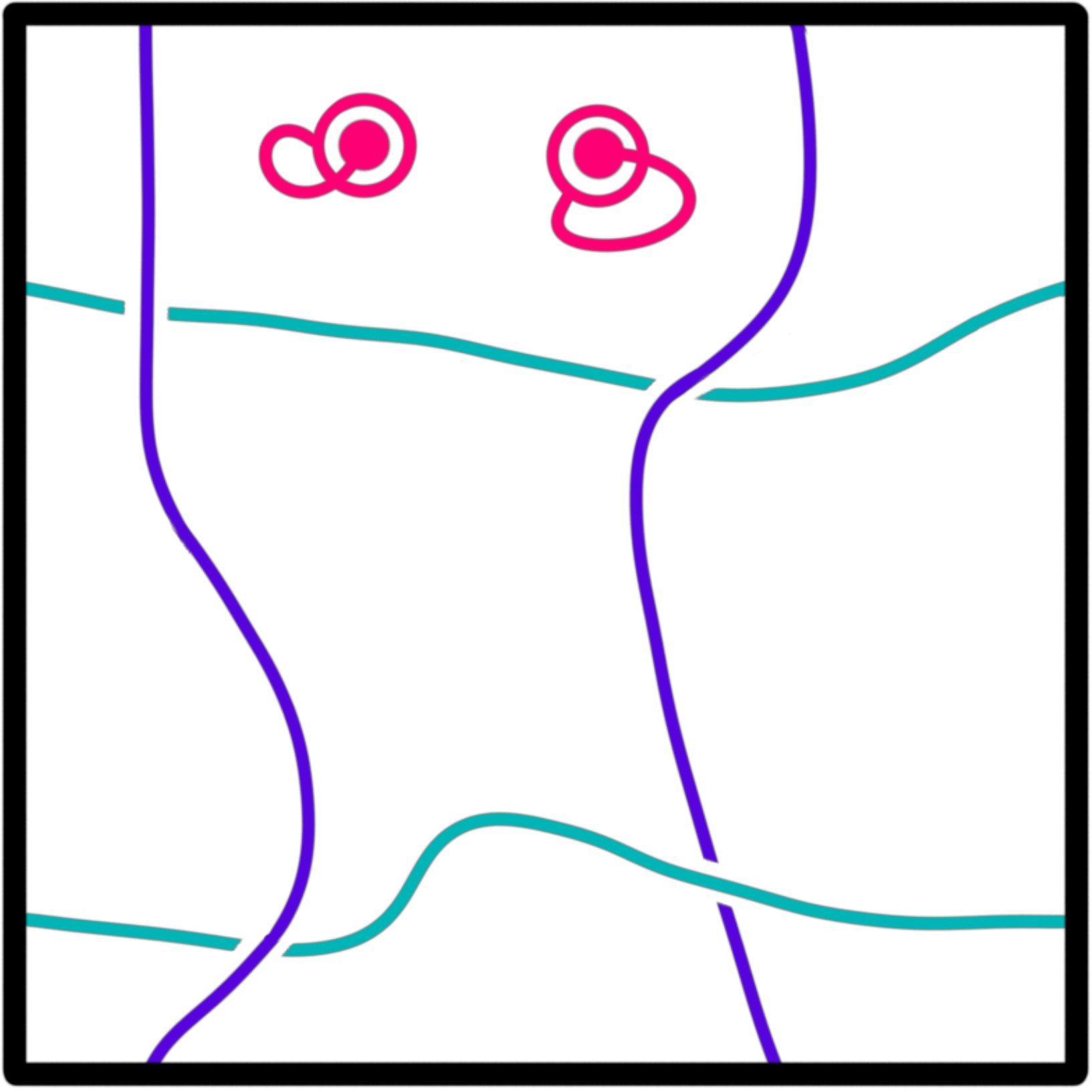}
        \caption{}
        \label{fig:untangling_pi_plus_4,4,5_1_usual_8_2nd}
    \end{subfigure}
    \begin{subfigure}[b]{0.19\textwidth}
        \centering
        \includegraphics[width=0.975\textwidth]{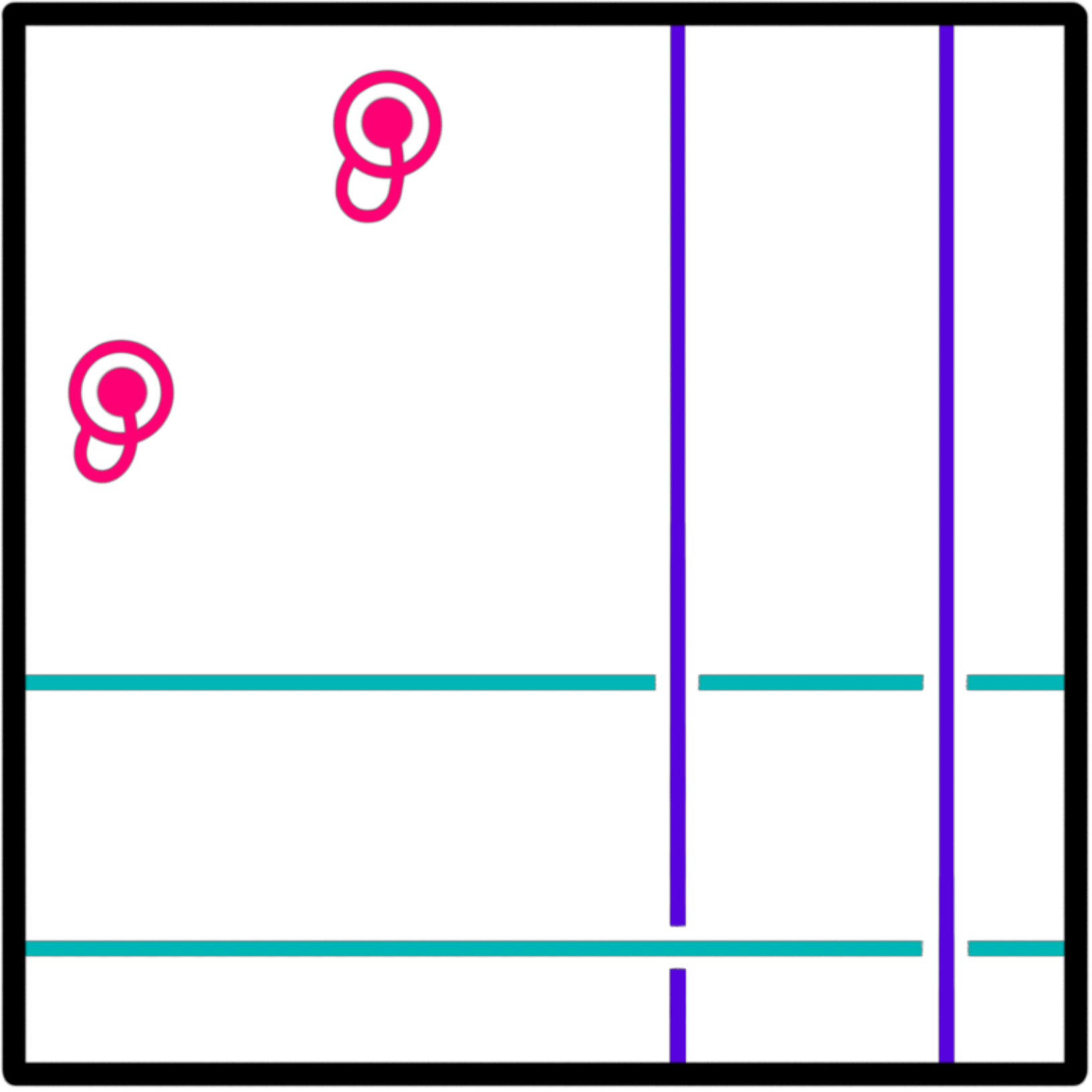}
        \caption{}
        \label{fig:untangling_pi_plus_4,4,5_1_usual_9_2nd}
    \end{subfigure}
    \begin{subfigure}[b]{0.19\textwidth}
        \centering
        \includegraphics[width=0.975\textwidth]{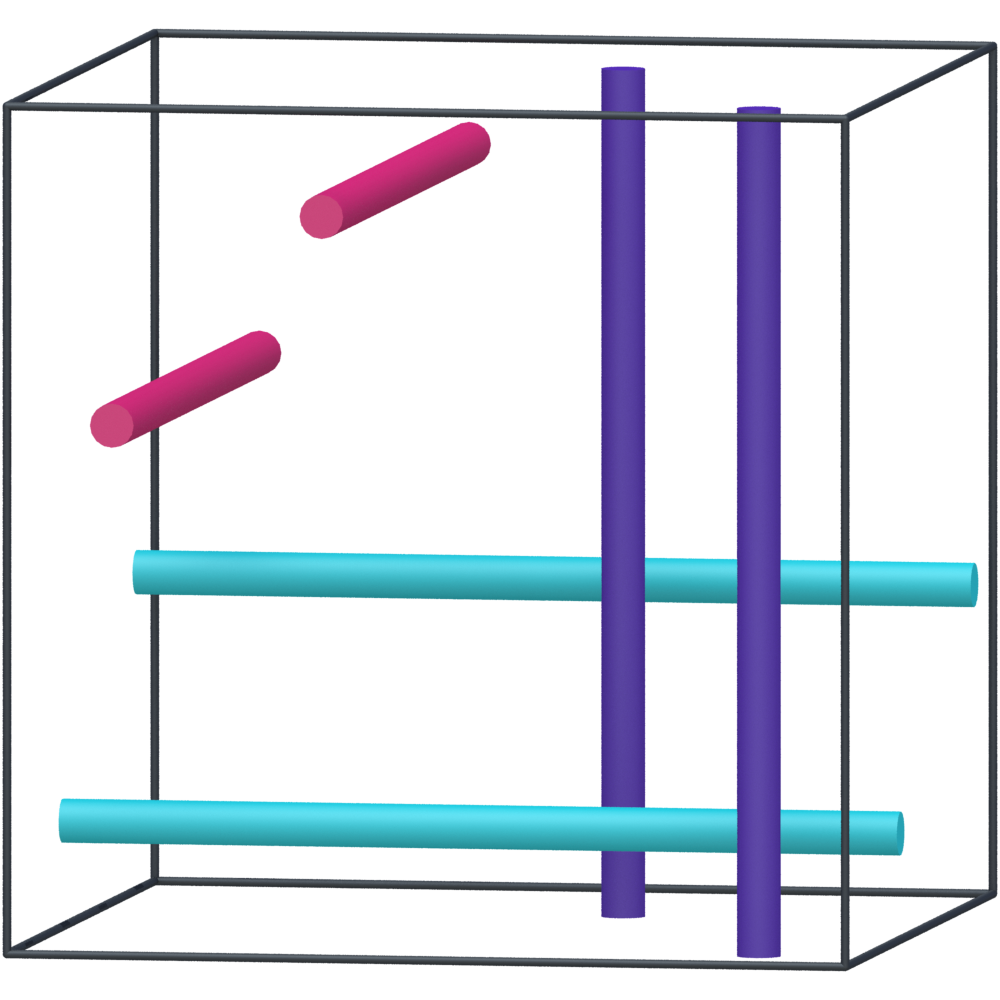}
        \caption{}
        \label{fig:pi_sharp_uc}
    \end{subfigure}
    \caption{Recomputing the upper bound for the untangling number of the 3-periodic tangle whose unit cell is shown in Fig.~\ref{fig:pi_plus_u2_uc_xyz_and_tridia}. By changing the crossings where the crossing changes are applied, one obtains a nearest ground state that is not the $\Pi^{+}$ rod packing obtained in Fig.~\ref{fig:untangling_pi_plus_4,4,5_1_usual_first_time}, but has the same minimum crossing number triplet $(4,4,4)$, and belongs to its $\mathcal{G}$-family.}
    \label{fig:untangling_pi_plus_4,4,5_1_usual_second_time}
\end{figure}

Note that there can indeed be more than one nearest ground state for a given 3-periodic tangle. For example, in Fig.~\ref{fig:untangling_pi_plus_4,4,5_1_usual_second_time} we recompute the upper bound that we previously computed in Fig.~\ref{fig:untangling_pi_plus_4,4,5_1_usual_first_time}, by changing the crossings where the crossing changes are applied. This results in a ground state whose unit cell is shown in Fig.~\ref{fig:pi_sharp_uc}, which is not the $\Pi^{+}$ rod packing previously obtained, but has the same minimum crossing number triplet $(4,4,4)$, and belongs to its $\mathcal{G}$-family.

\bigbreak

The untangling number as well as the nearest ground state also truly depend on the chosen unit cell. See the example displayed in Fig.~\ref{fig:untangling_pi_plus_4,4,5_1_usual_vol2}, where we choose a unit cell that is the double of that in Fig.~\ref{fig:pi_plus_u2_uc}. The upper bound is doubled, that is, four crossing changes are required. Moreover, by mixing the crossings where these crossing changes are applied, we obtain a nearest ground state that is a mix of the two ground states obtained in Fig.~\ref{fig:untangling_pi_plus_4,4,5_1_usual_first_time} and Fig.~\ref{fig:untangling_pi_plus_4,4,5_1_usual_second_time}, which would not be a ground state with respect to the unit cell of smaller volume.

\begin{figure}[hbtp]
    \centering
    \begin{subfigure}[b]{\textwidth}
        \centering
        \includegraphics[width=0.37\textwidth]{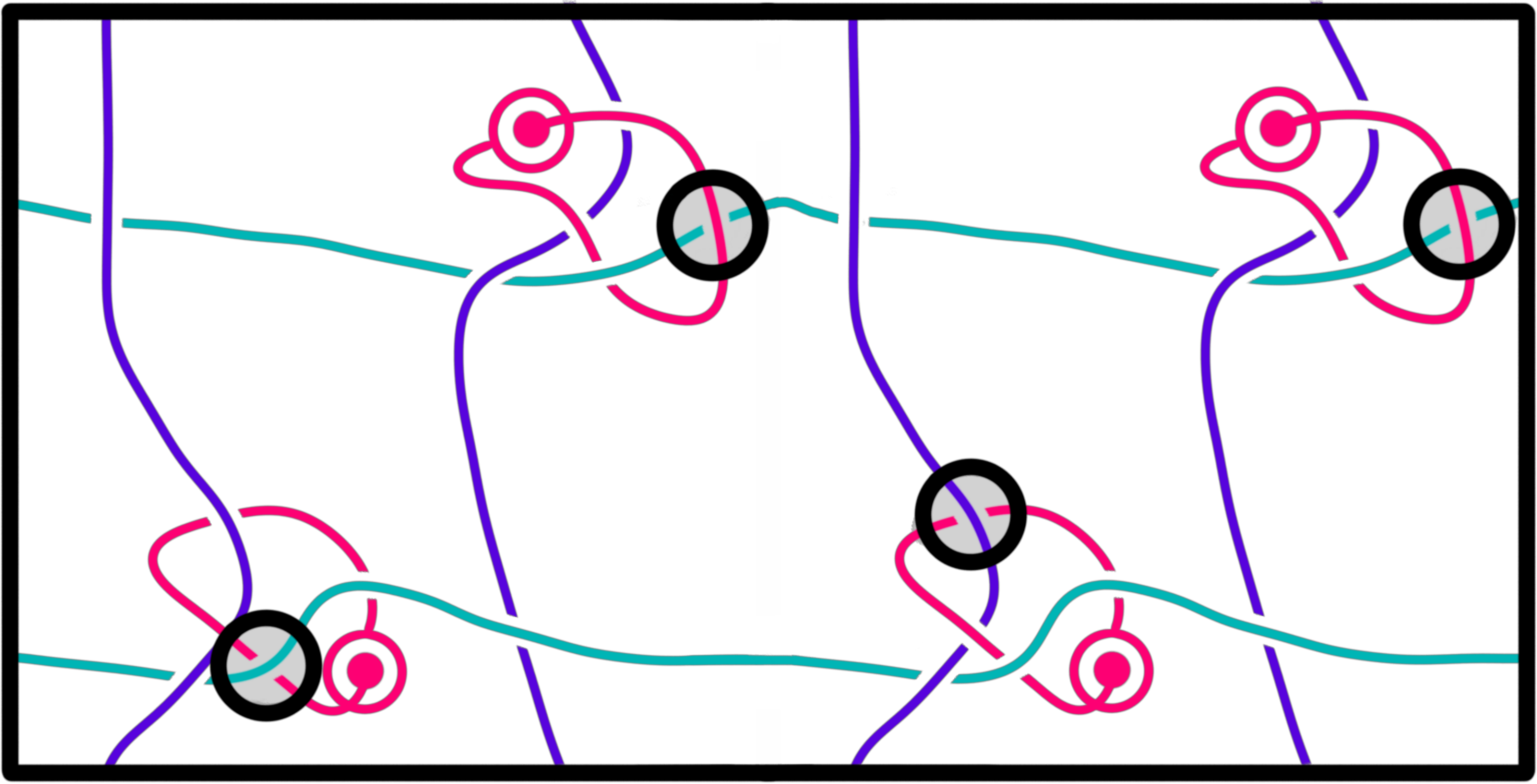}
        \caption{}
        \label{fig:pi_plus_u2_vol2_untangling}
    \end{subfigure}

    \vskip\baselineskip
    
    \begin{subfigure}[b]{\textwidth}
        \centering
        \includegraphics[width=0.37\textwidth]{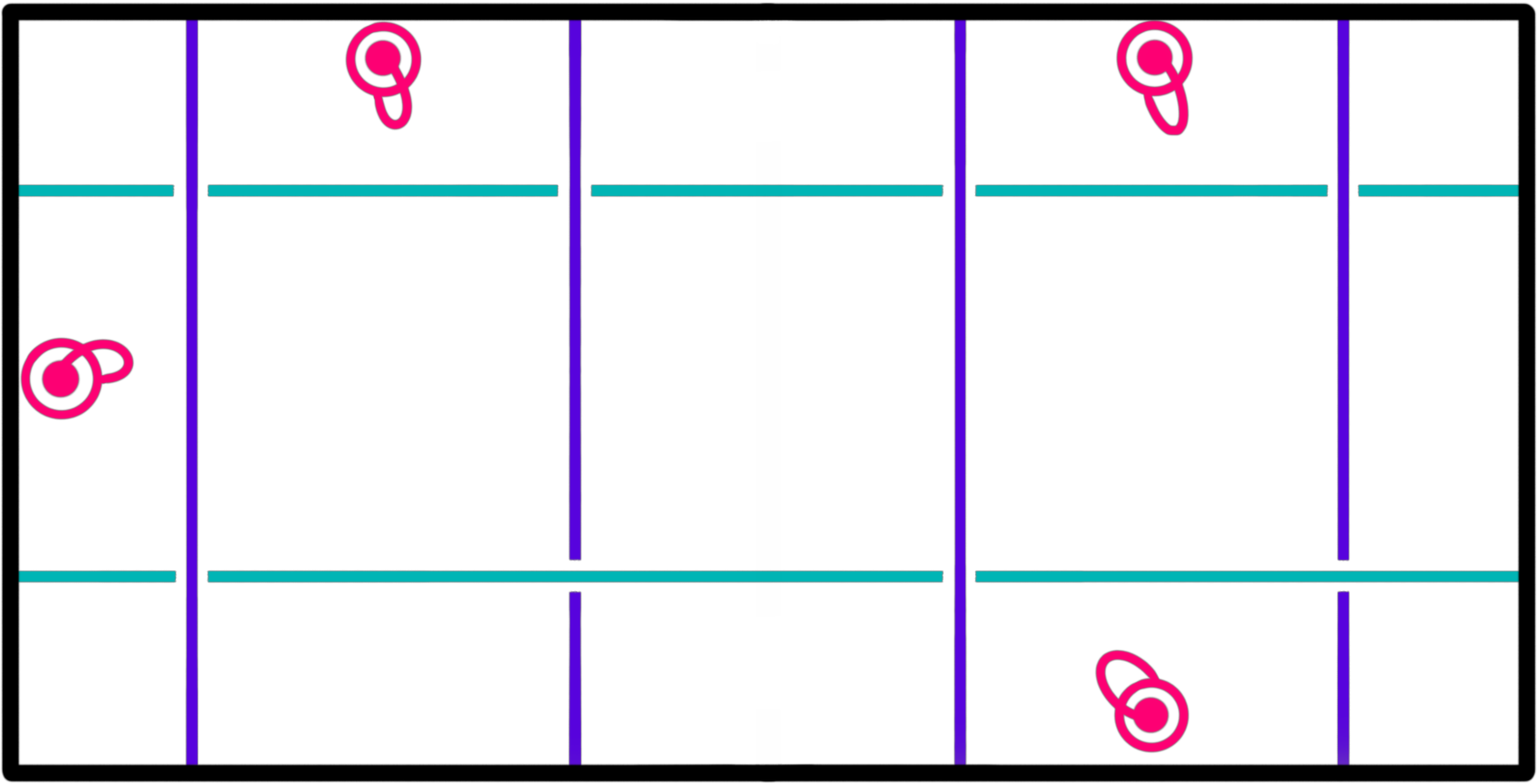}
        \caption{}
        \label{fig:pi_plus_pi_sharp_mix}
    \end{subfigure}
    \caption{Computing an upper bound for the untangling number of 3-periodic tangle associated to the unit cell displayed in Fig.~\ref{fig:pi_plus_u2_uc_xyz_and_tridia}, but with a different unit cell. Doubling the size of the unit cell doubles the upper bound that was obtained in Fig.~\ref{fig:untangling_pi_plus_4,4,5_1_usual_first_time} and Fig.~\ref{fig:untangling_pi_plus_4,4,5_1_usual_second_time}, and diversifying the crossings where the crossing changes are applied yield a nearest ground state that would not be one with respect to the unit cell of smaller volume.}
    \label{fig:untangling_pi_plus_4,4,5_1_usual_vol2}
\end{figure}

\bigbreak

In light of this, one may wish to stabilise the untangling number, which can be done by considering the minimum of all untangling numbers among all unit cells of a given 3-periodic tangle, leading to the following definition of an invariant.

\begin{definition}\label{def:min_untangling_numb}
The \emph{minimum untangling number} of a 3-periodic tangle $K$ is defined to be the least untangling number $u(K,U)$ over all unit cells $U$ of $K$.
\end{definition}

\bigbreak

\begin{figure}[hbtp]
    \centering
    \begin{subfigure}[b]{0.19\textwidth}
        \centering
        \includegraphics[width=0.975\textwidth]{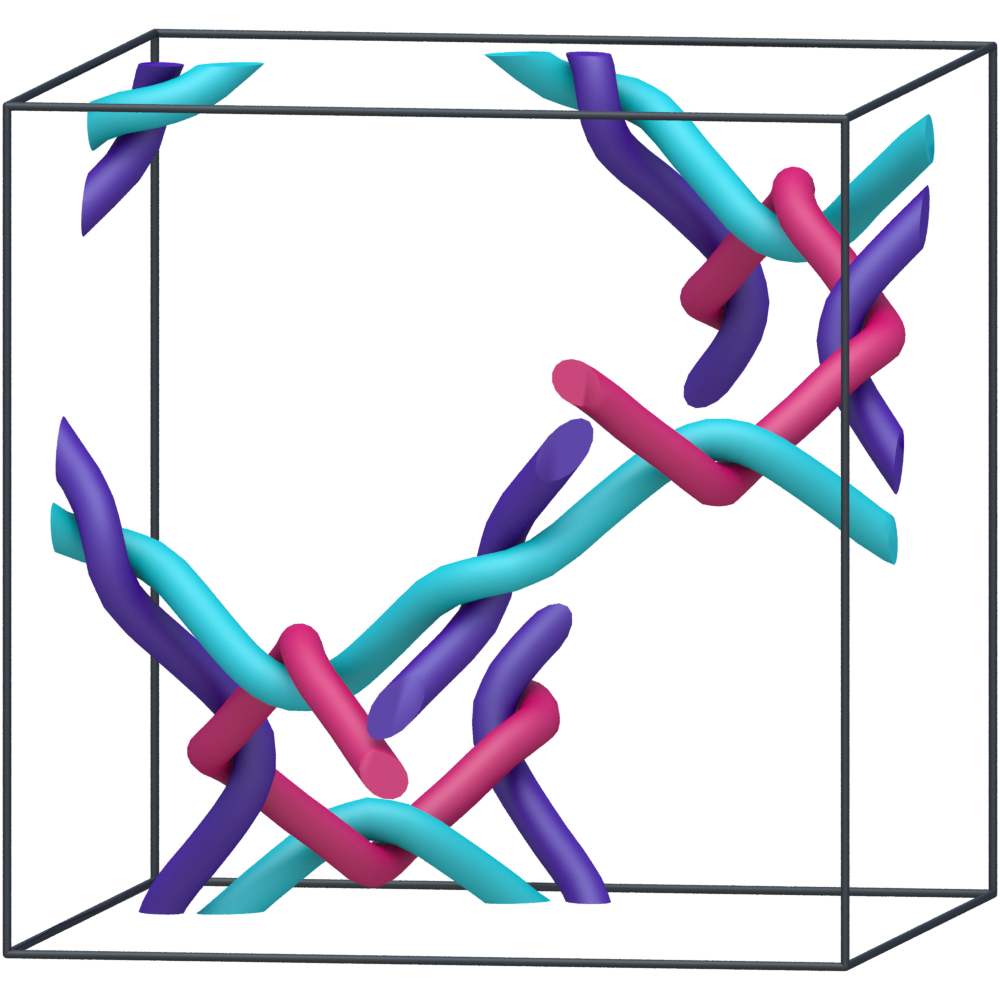}
        \caption{}
        \label{fig:1p6on2_to_the_6_uc}
    \end{subfigure}
    \begin{subfigure}[b]{0.19\textwidth}
        \centering
        \includegraphics[width=0.975\textwidth]{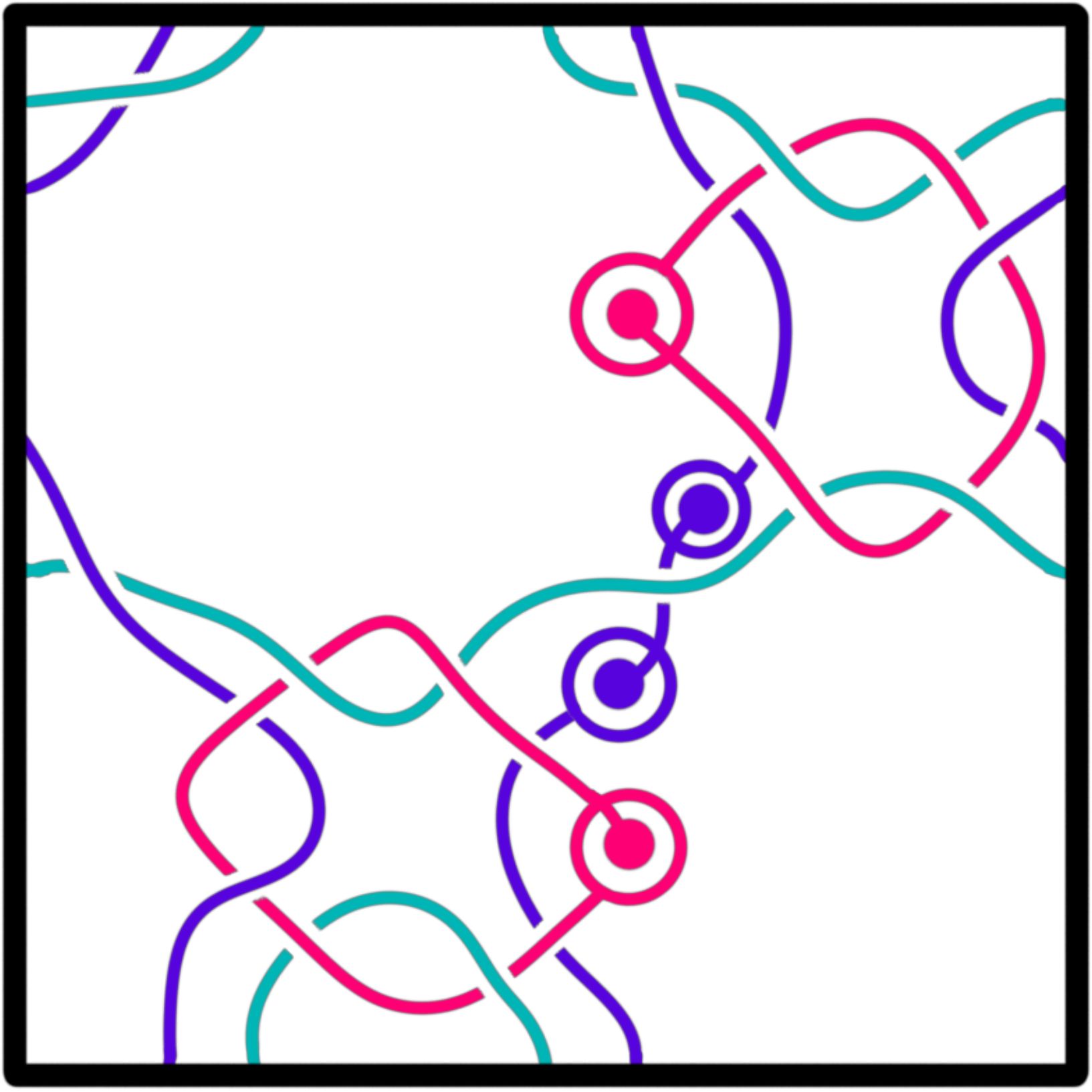}
        \caption{}
        \label{fig:1p6on2_to_the_6_dia_1}
    \end{subfigure}
    \begin{subfigure}[b]{0.19\textwidth}
        \centering
        \includegraphics[width=0.975\textwidth]{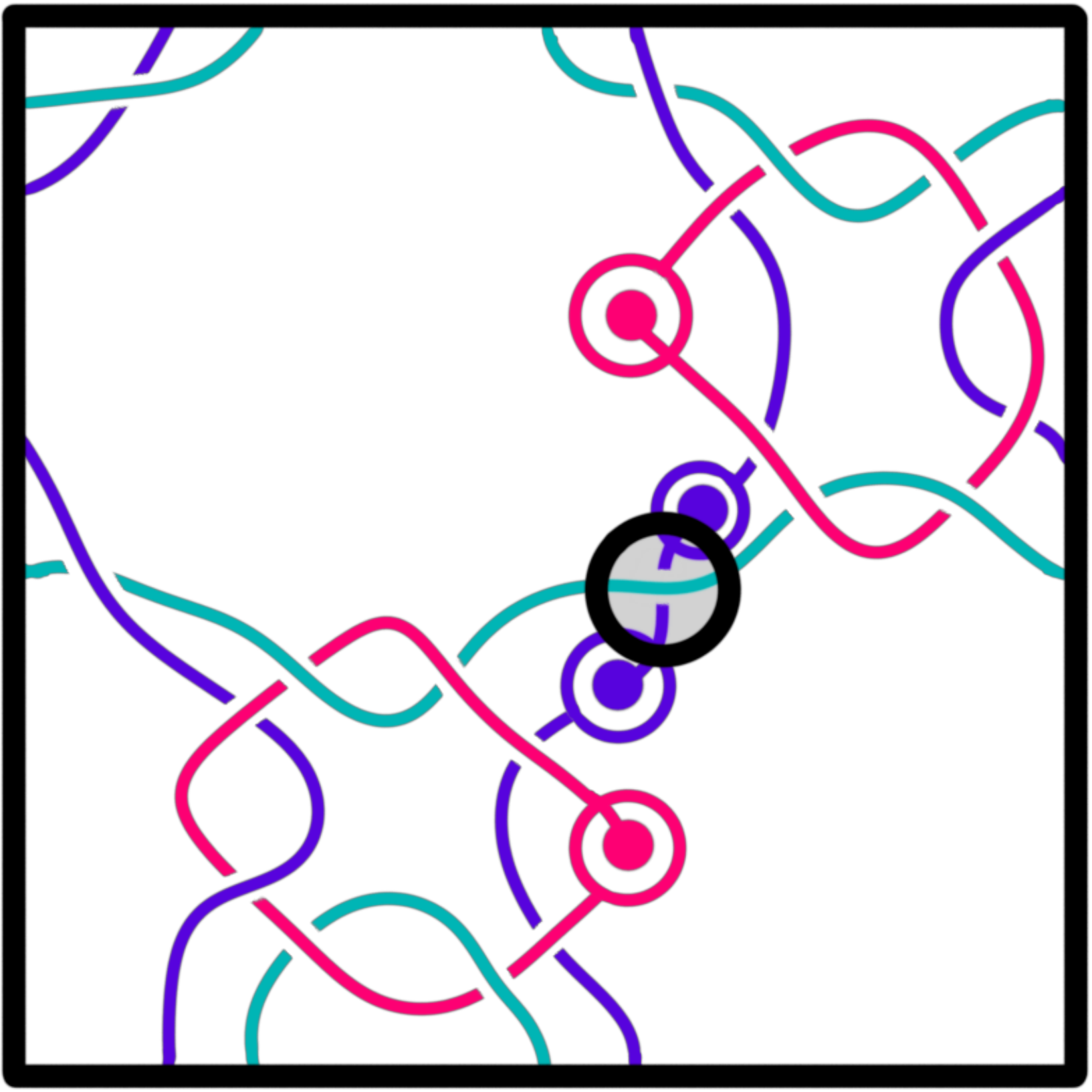}
        \caption{}
        \label{fig:1p6on2_to_the_6_dia_2}
    \end{subfigure}
    \begin{subfigure}[b]{0.19\textwidth}
        \centering
        \includegraphics[width=0.975\textwidth]{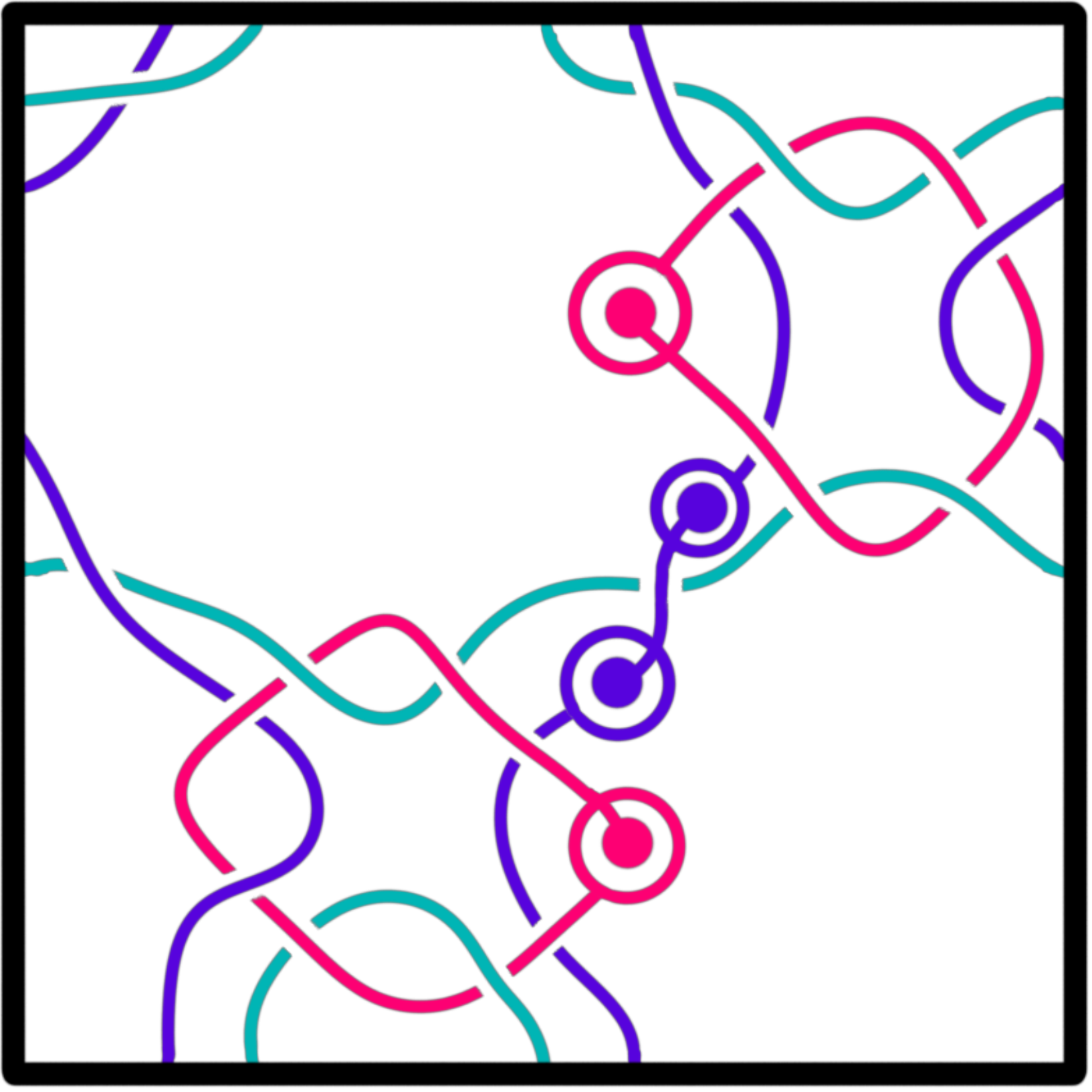}
        \caption{}
        \label{fig:1p6on2_to_the_6_dia_3}
    \end{subfigure}
    \begin{subfigure}[b]{0.19\textwidth}
        \centering
        \includegraphics[width=0.975\textwidth]{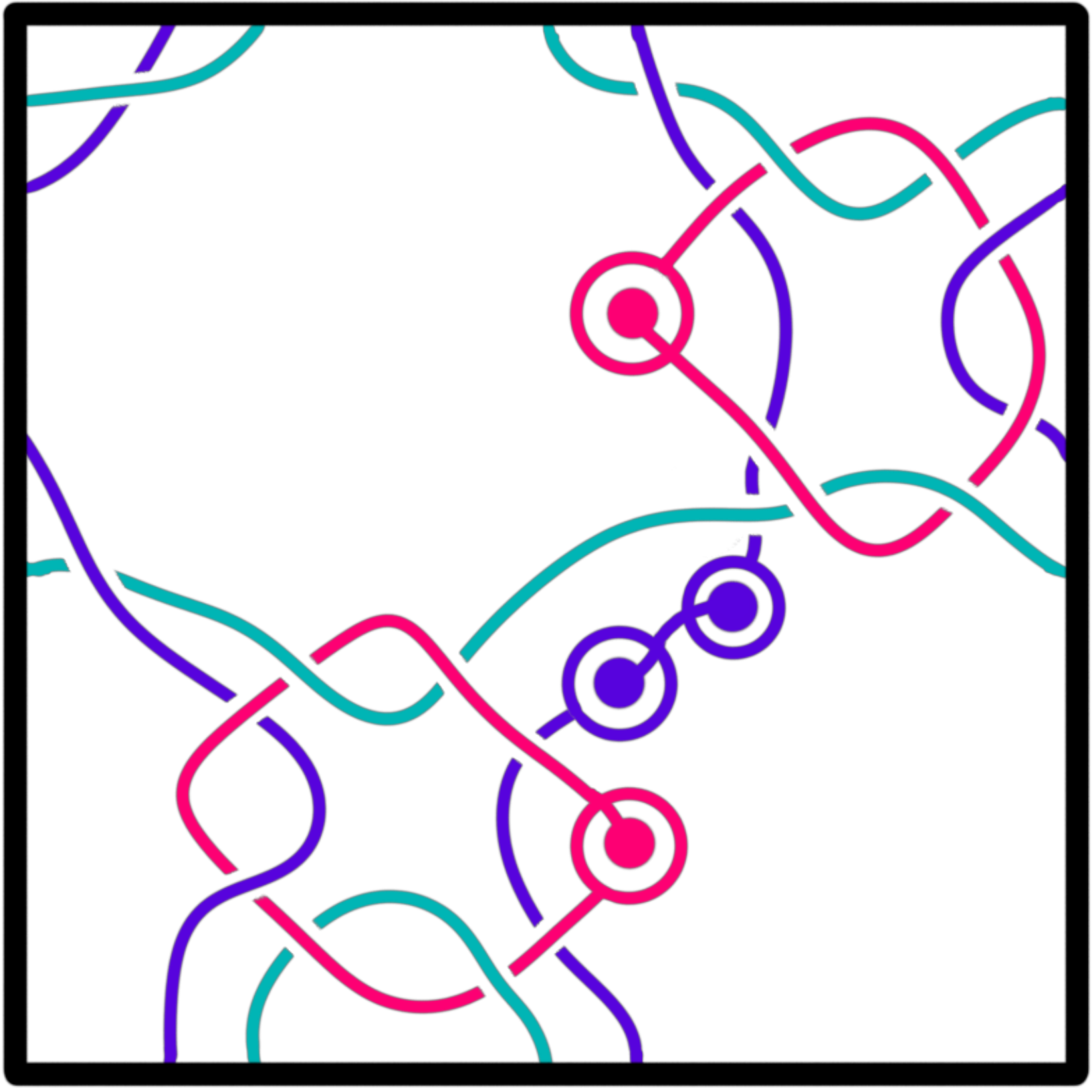}
        \caption{}
        \label{fig:1p6on2_to_the_6_dia_4}
    \end{subfigure}

    \vskip\baselineskip
    
    \begin{subfigure}[b]{0.19\textwidth}
        \centering
        \includegraphics[width=0.975\textwidth]{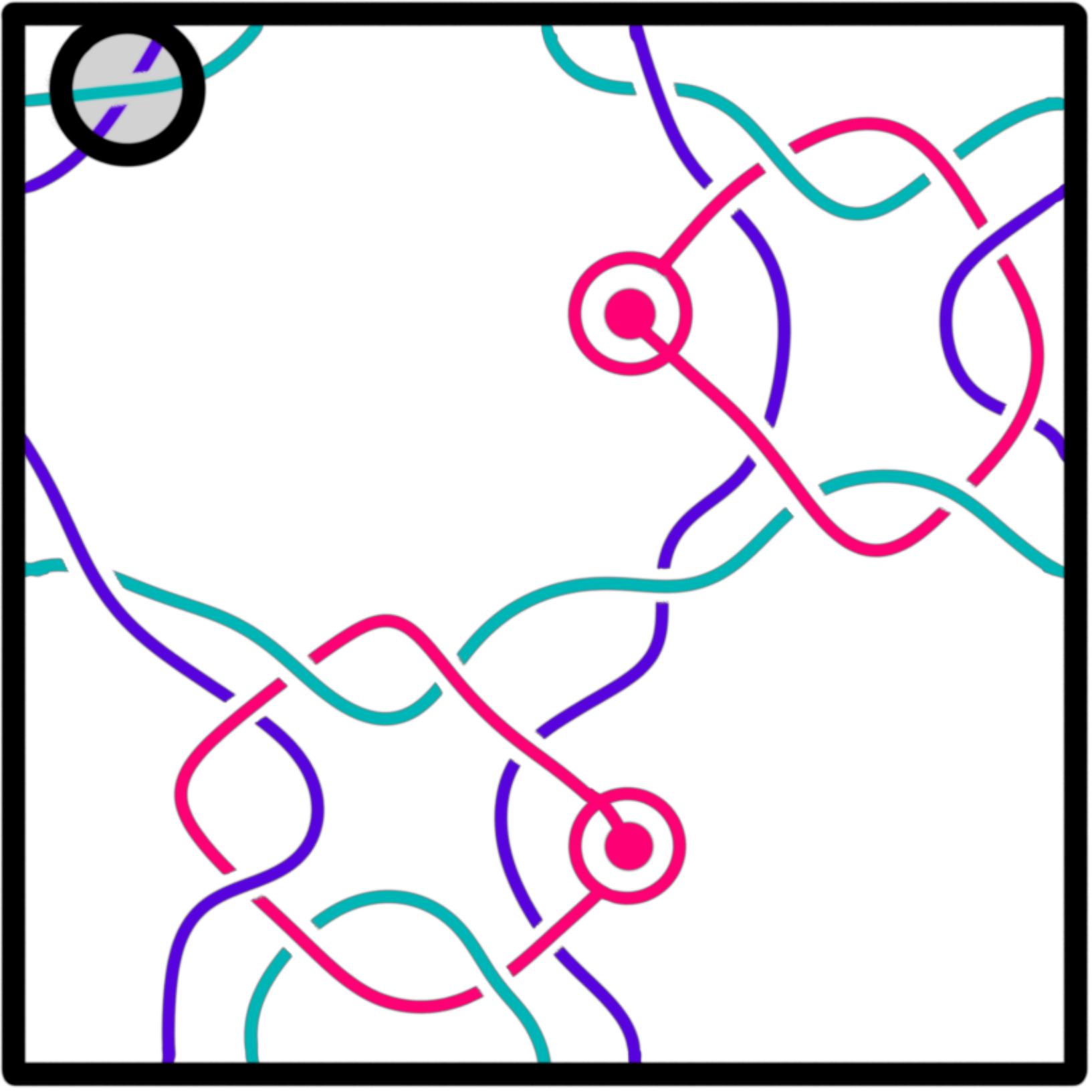}
        \caption{}
        \label{fig:1p6on2_to_the_6_dia_5}
    \end{subfigure}
    \begin{subfigure}[b]{0.19\textwidth}
        \centering
        \includegraphics[width=0.975\textwidth]{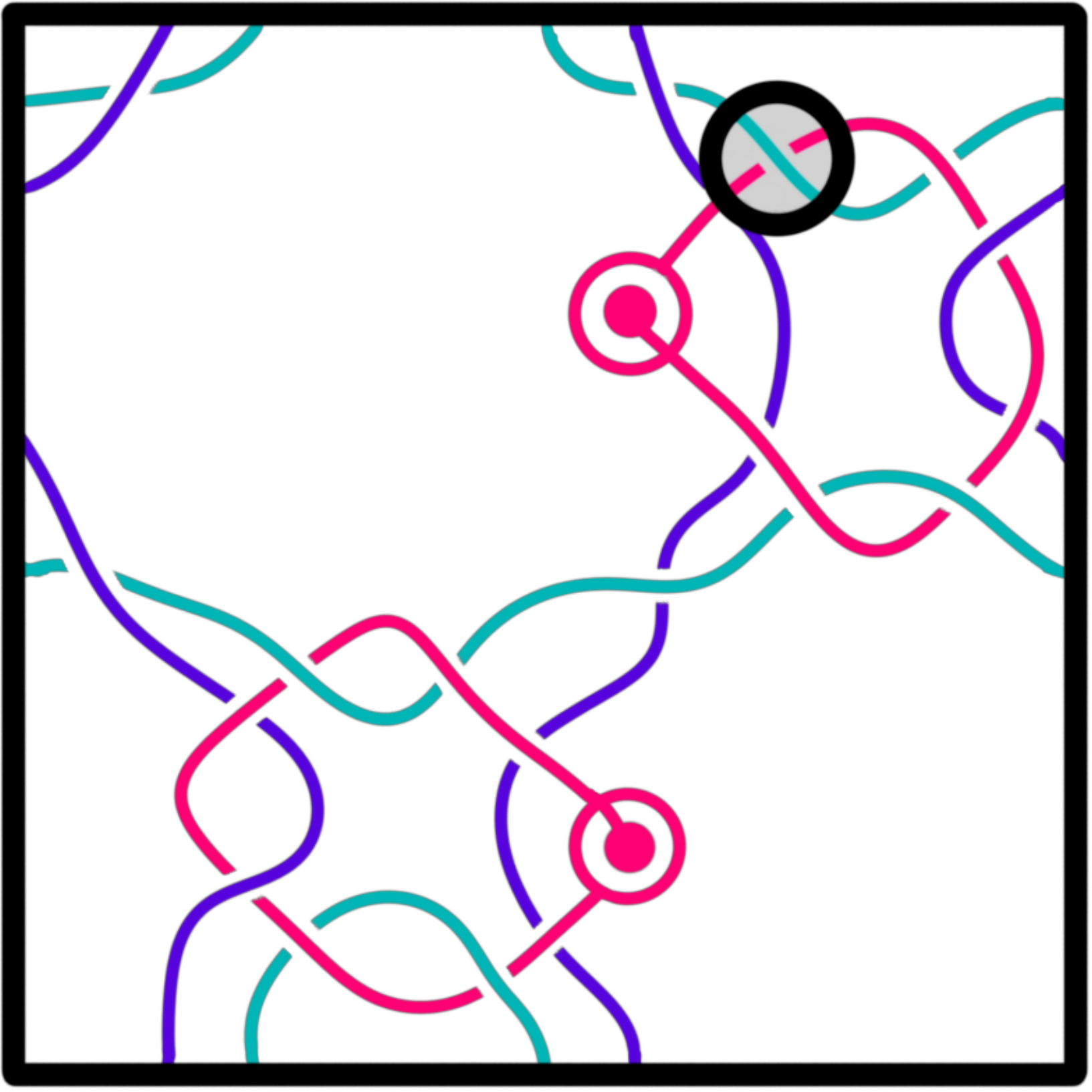}
        \caption{}
        \label{fig:1p6on2_to_the_6_dia_6}
    \end{subfigure}
    \begin{subfigure}[b]{0.19\textwidth}
        \centering
        \includegraphics[width=0.975\textwidth]{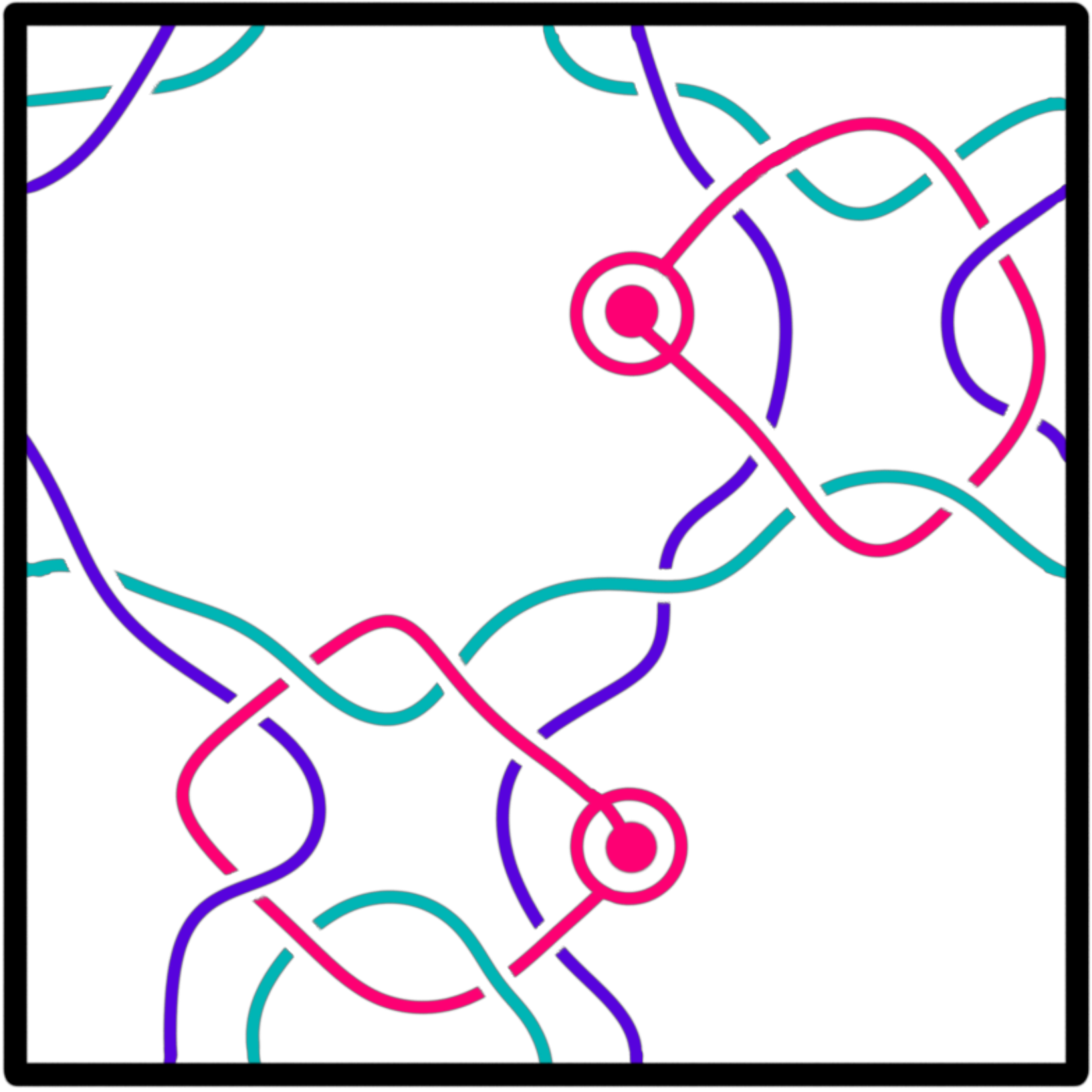}
        \caption{}
        \label{fig:1p6on2_to_the_6_dia_7}
    \end{subfigure}
    \begin{subfigure}[b]{0.19\textwidth}
        \centering
        \includegraphics[width=0.975\textwidth]{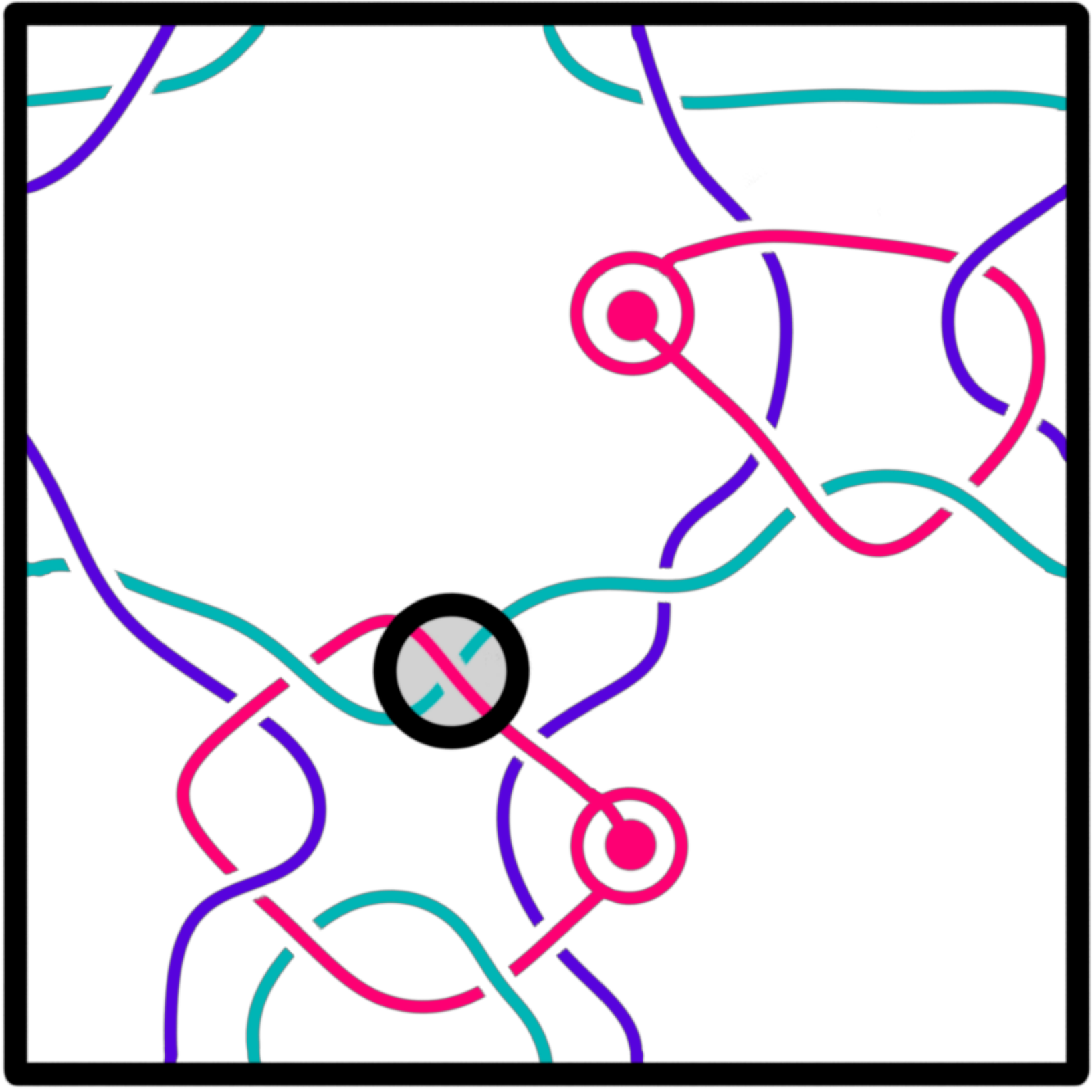}
        \caption{}
        \label{fig:1p6on2_to_the_6_dia_8}
    \end{subfigure}
    \begin{subfigure}[b]{0.19\textwidth}
        \centering
        \includegraphics[width=0.975\textwidth]{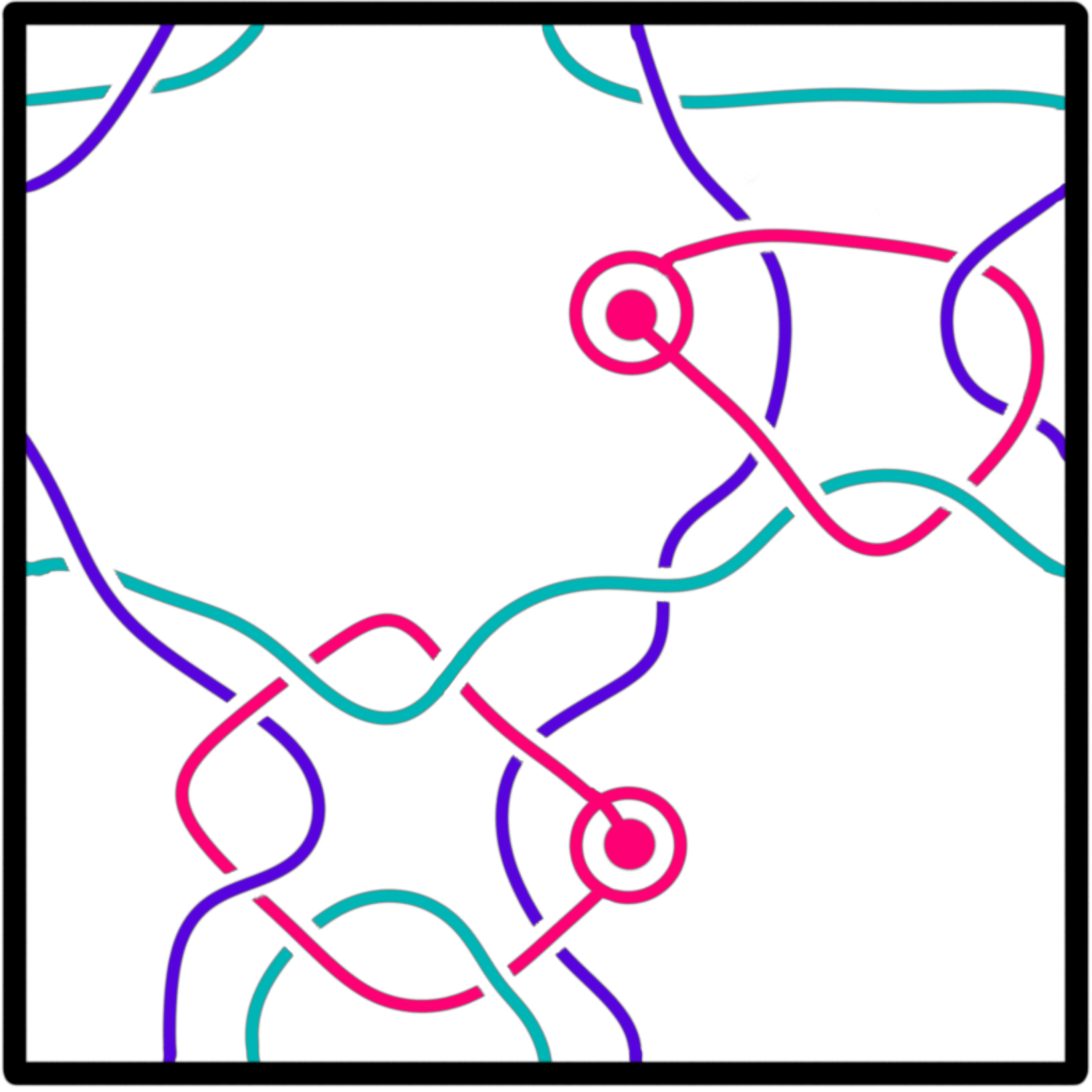}
        \caption{}
        \label{fig:1p6on2_to_the_6_dia_9}
    \end{subfigure}

    \vskip\baselineskip

    \begin{subfigure}[b]{0.19\textwidth}
        \centering
        \includegraphics[width=0.975\textwidth]{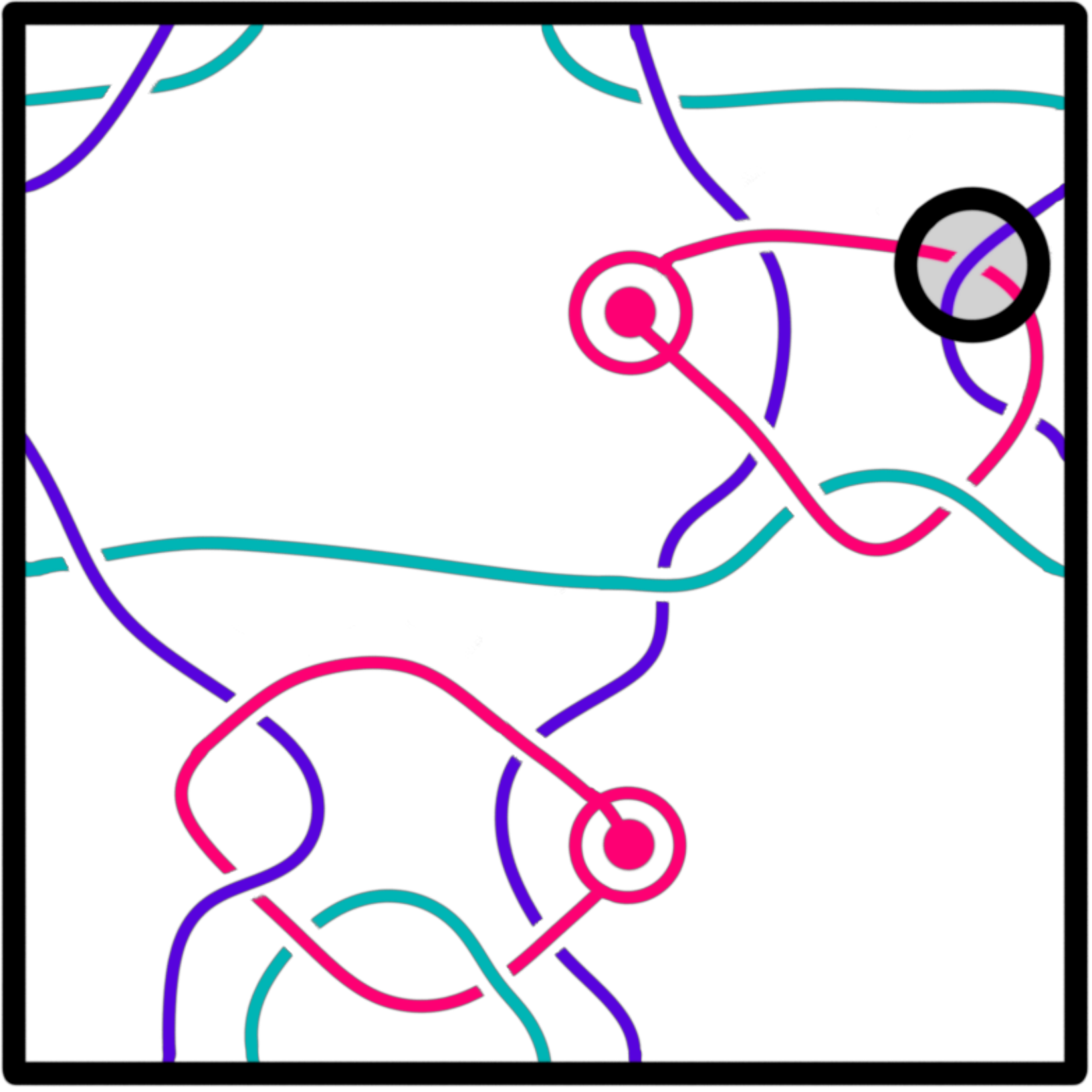}
        \caption{}
        \label{fig:1p6on2_to_the_6_dia_10}
    \end{subfigure}
    \begin{subfigure}[b]{0.19\textwidth}
        \centering
        \includegraphics[width=0.975\textwidth]{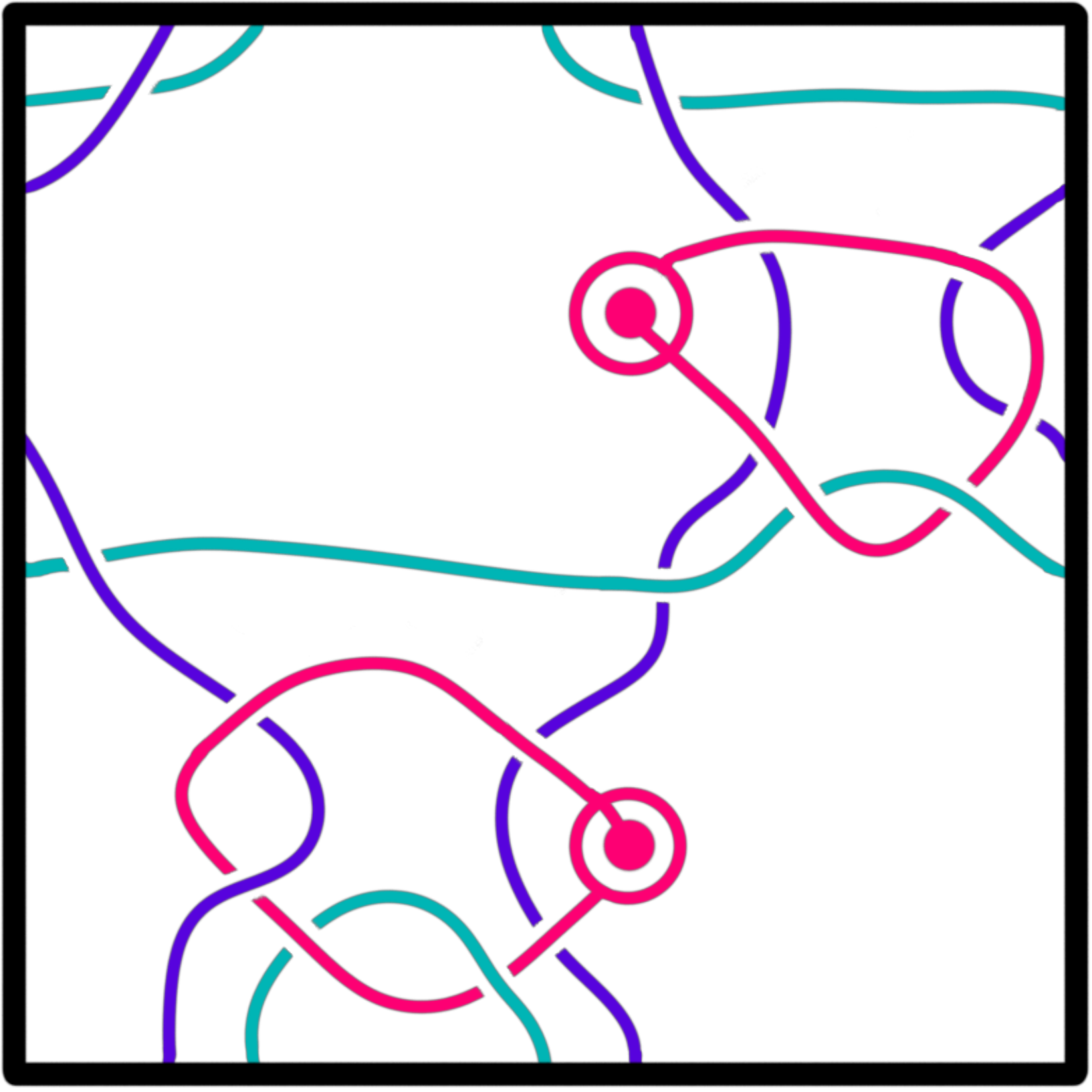}
        \caption{}
        \label{fig:1p6on2_to_the_6_dia_11}
    \end{subfigure}
    \begin{subfigure}[b]{0.19\textwidth}
        \centering
        \includegraphics[width=0.975\textwidth]{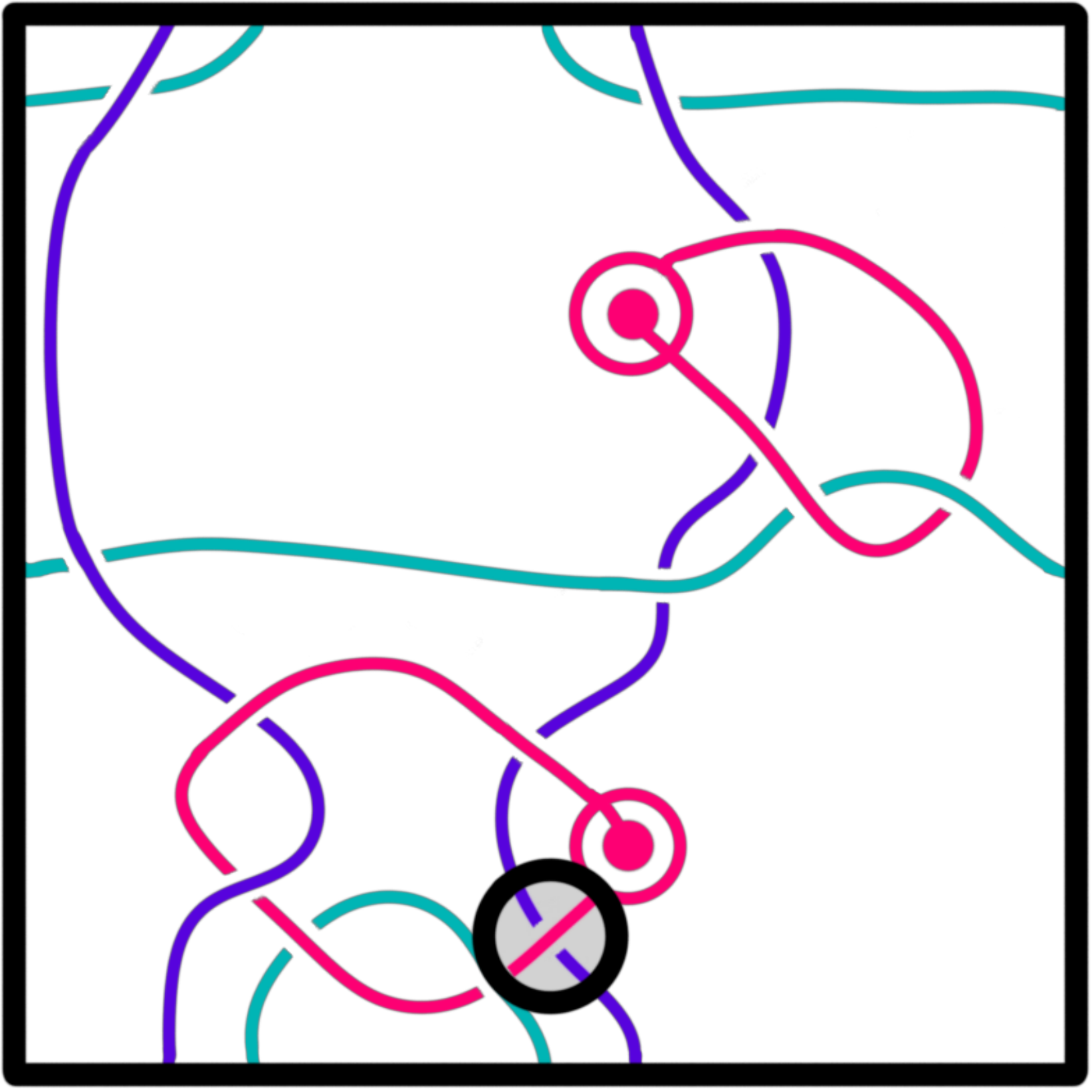}
        \caption{}
        \label{fig:1p6on2_to_the_6_dia_12}
    \end{subfigure}
    \begin{subfigure}[b]{0.19\textwidth}
        \centering
        \includegraphics[width=0.975\textwidth]{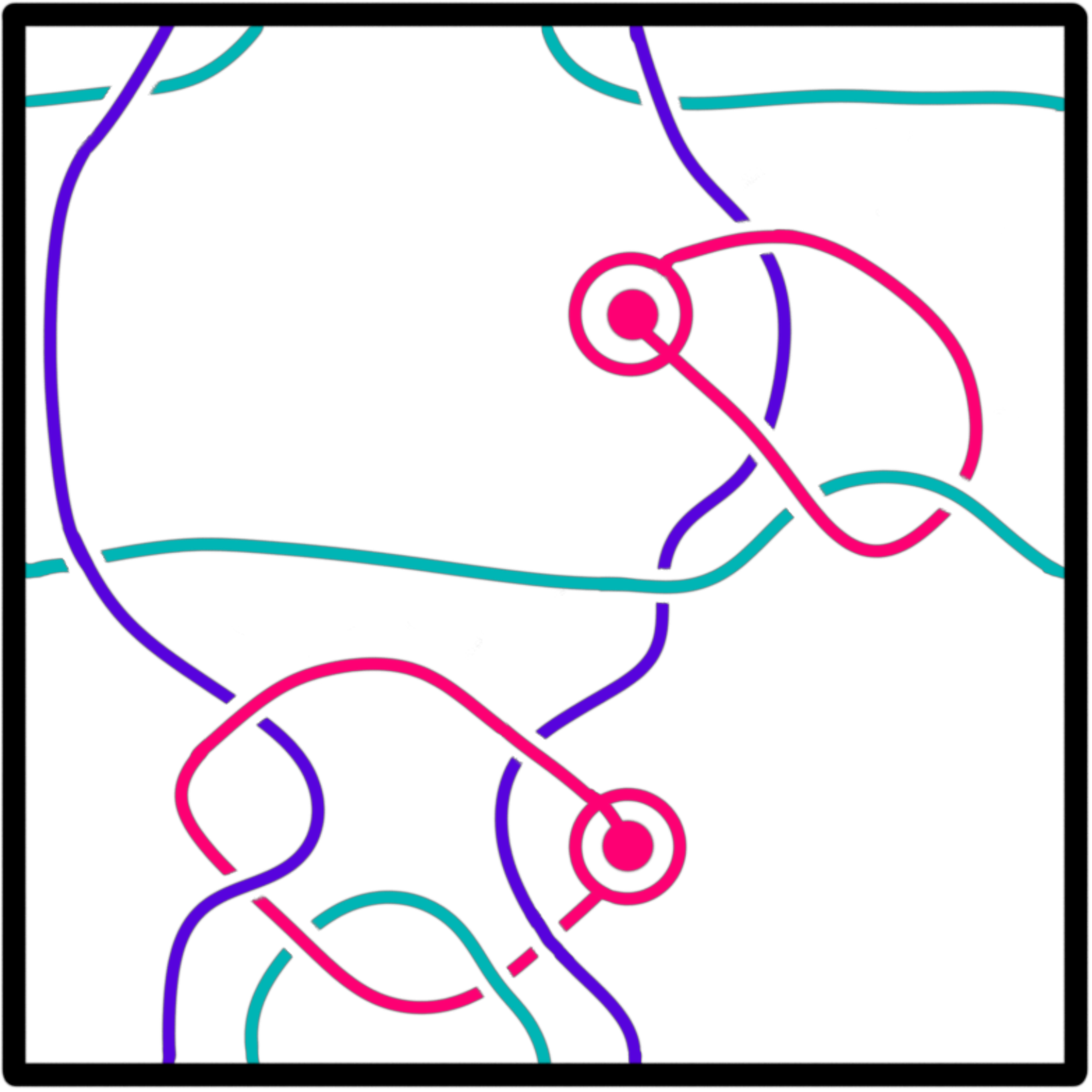}
        \caption{}
        \label{fig:1p6on2_to_the_6_dia_13}
    \end{subfigure}
    \begin{subfigure}[b]{0.19\textwidth}
        \centering
        \includegraphics[width=0.975\textwidth]{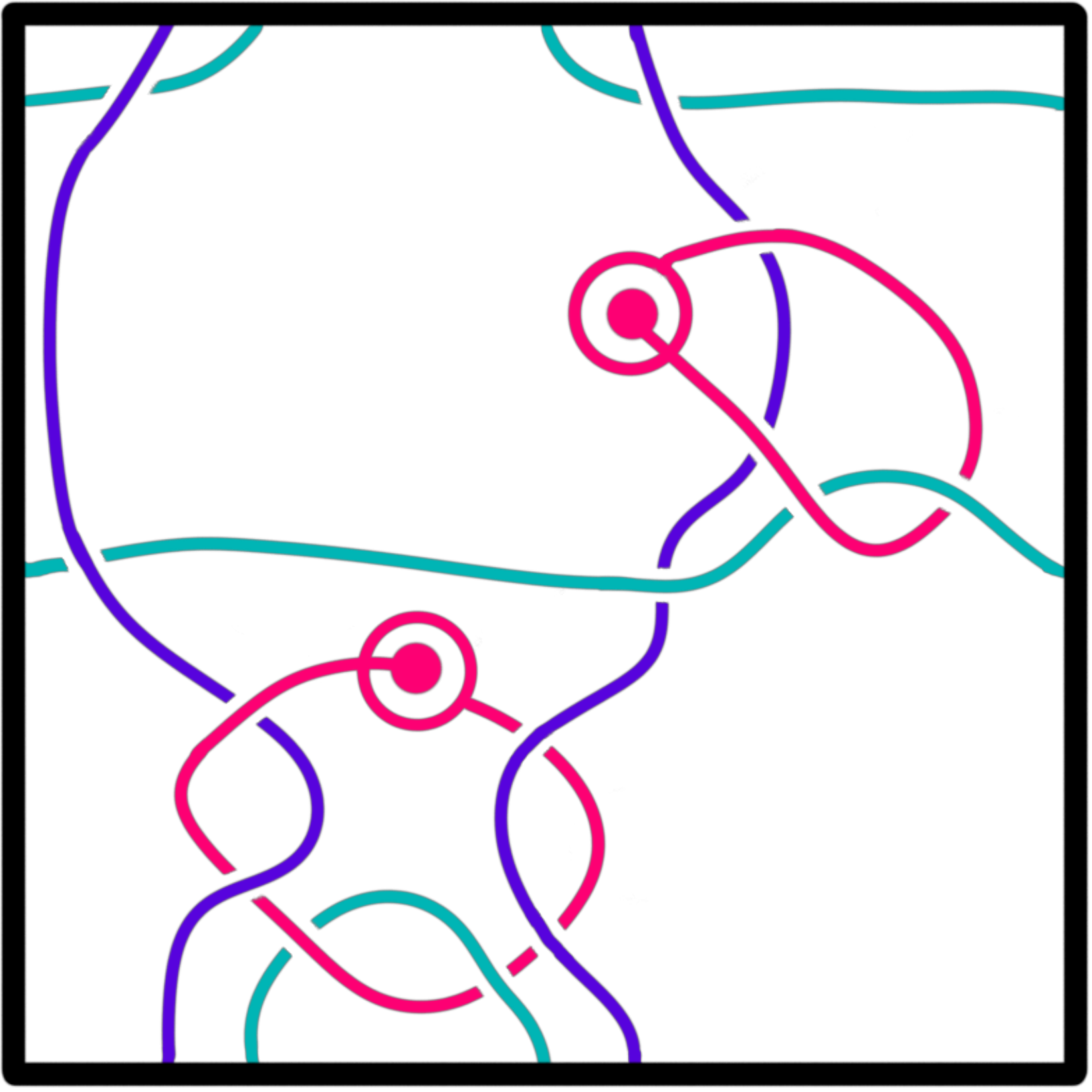}
        \caption{}
        \label{fig:1p6on2_to_the_6_dia_14}
    \end{subfigure}

    \vskip\baselineskip

    \begin{subfigure}[b]{0.19\textwidth}
        \centering
        \includegraphics[width=0.975\textwidth]{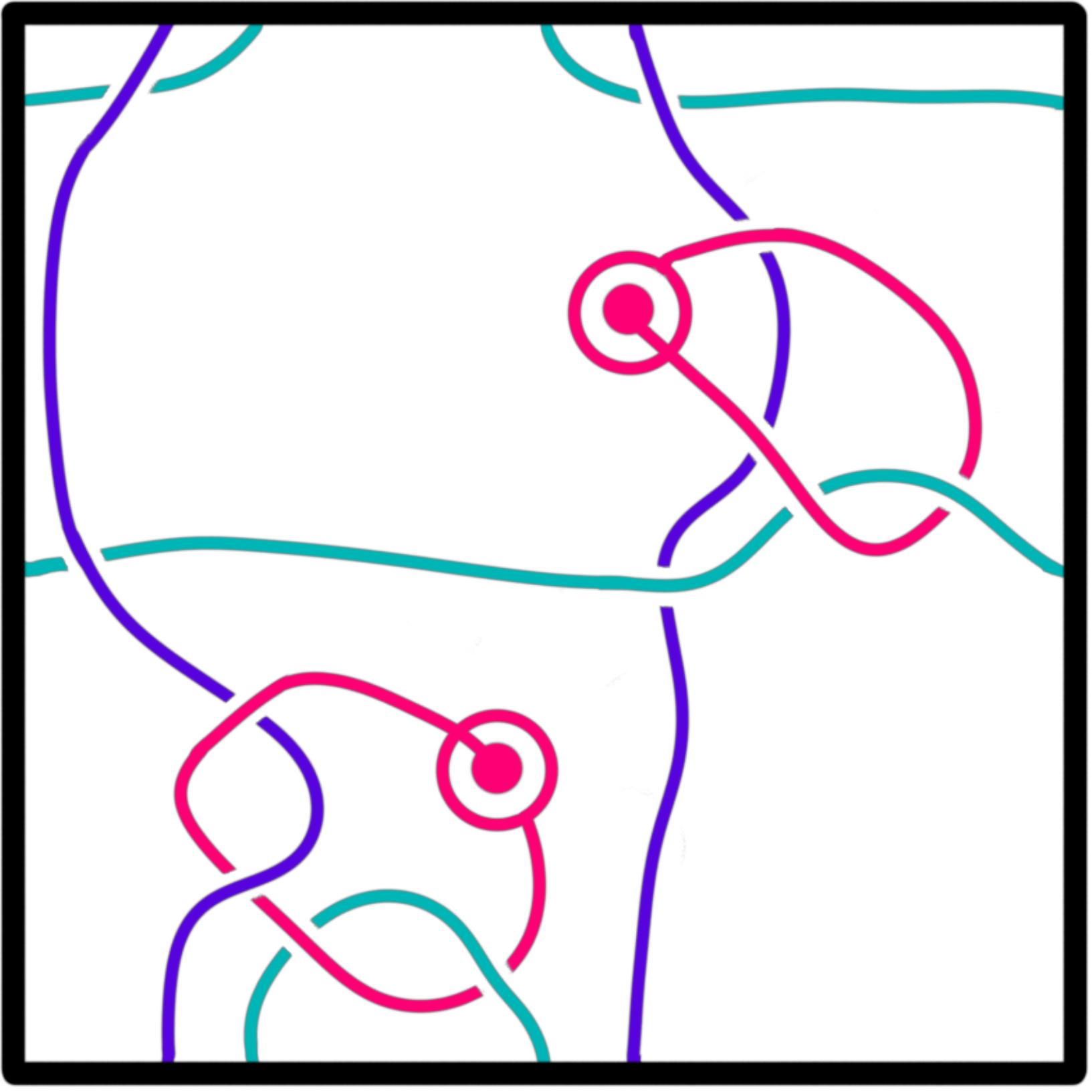}
        \caption{}
        \label{fig:1p6on2_to_the_6_dia_15}
    \end{subfigure}
    \begin{subfigure}[b]{0.19\textwidth}
        \centering
        \includegraphics[width=0.975\textwidth]{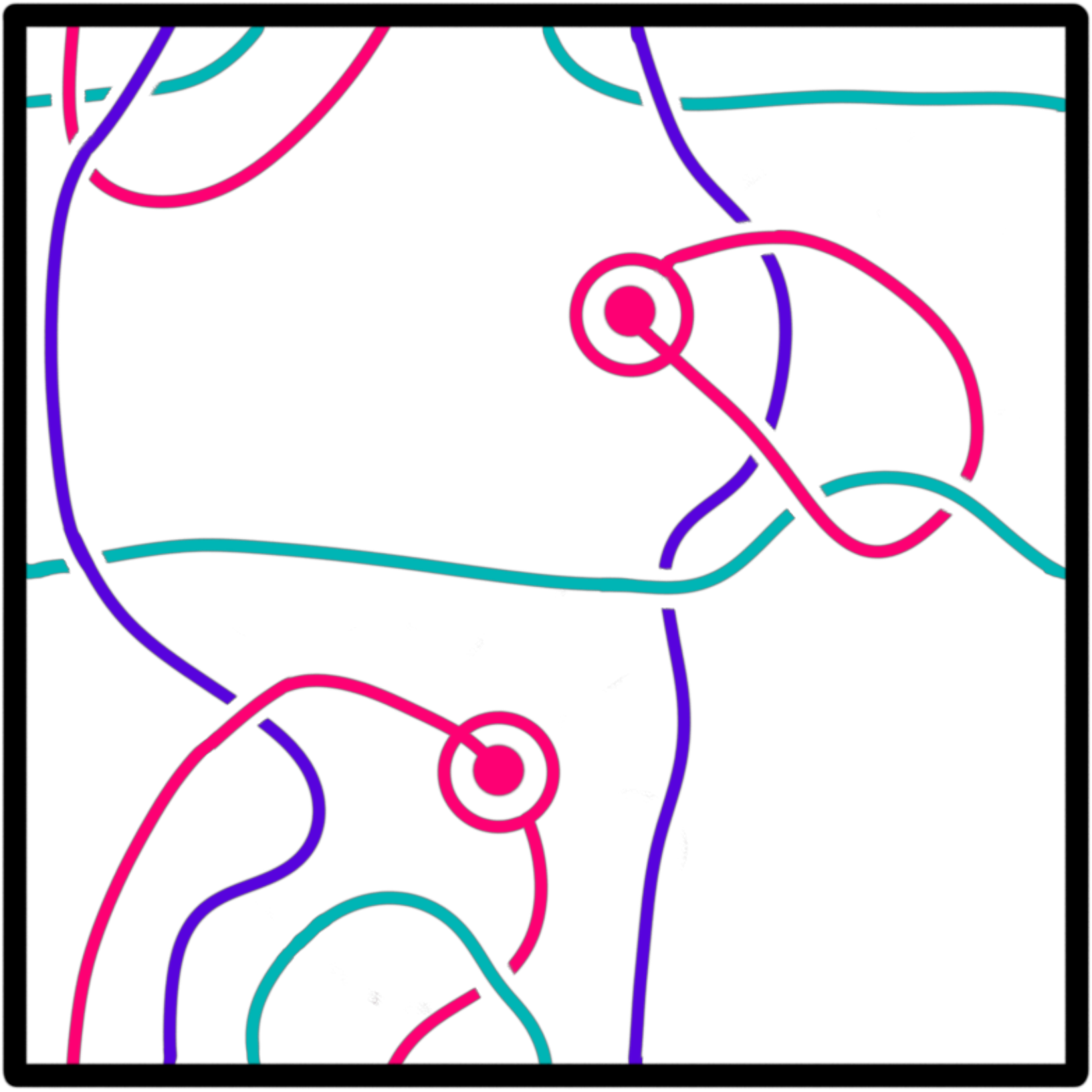}
        \caption{}
        \label{fig:1p6on2_to_the_6_dia_16}
    \end{subfigure}
    \begin{subfigure}[b]{0.19\textwidth}
        \centering
        \includegraphics[width=0.975\textwidth]{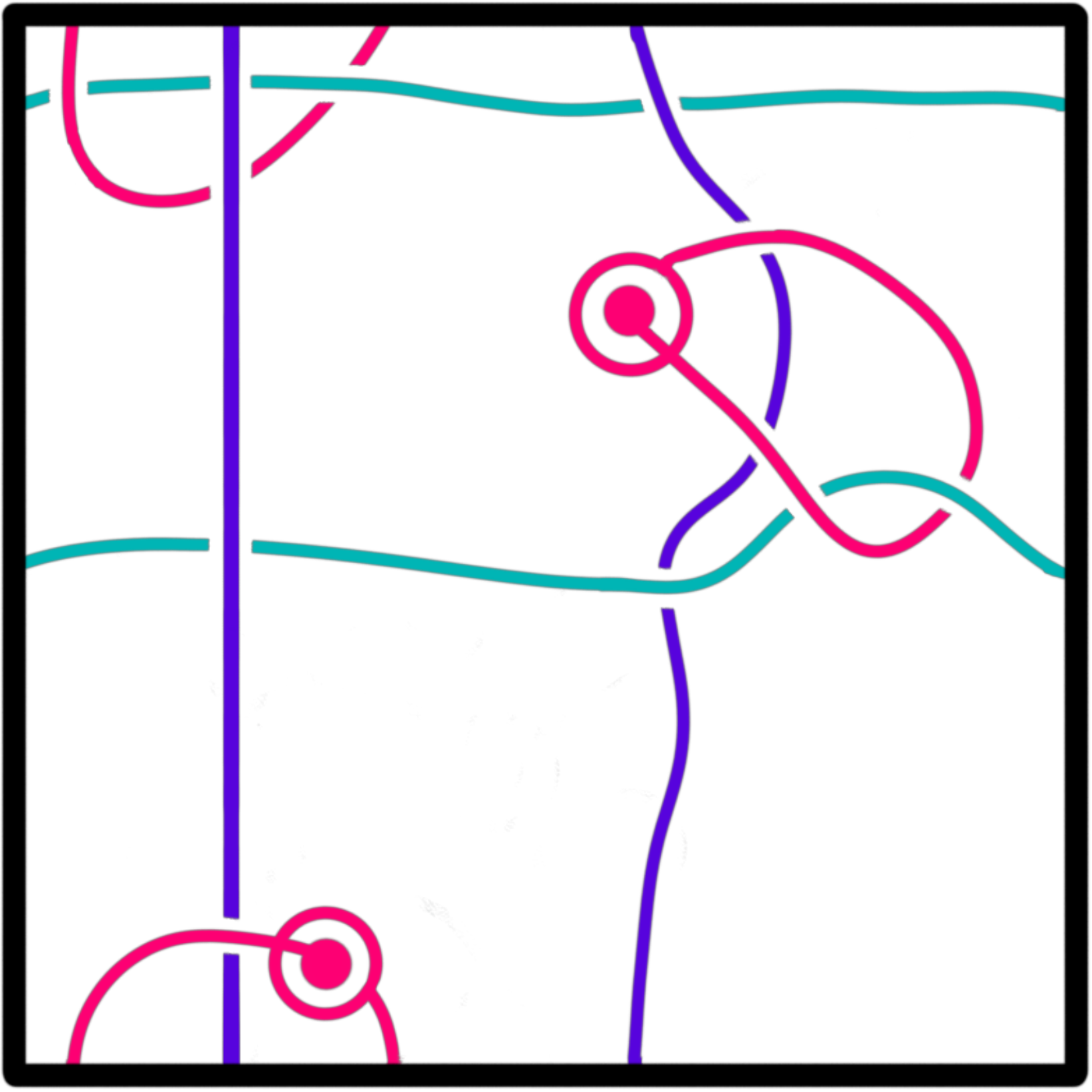}
        \caption{}
        \label{fig:1p6on2_to_the_6_dia_17}
    \end{subfigure}
    \begin{subfigure}[b]{0.19\textwidth}
        \centering
        \includegraphics[width=0.975\textwidth]{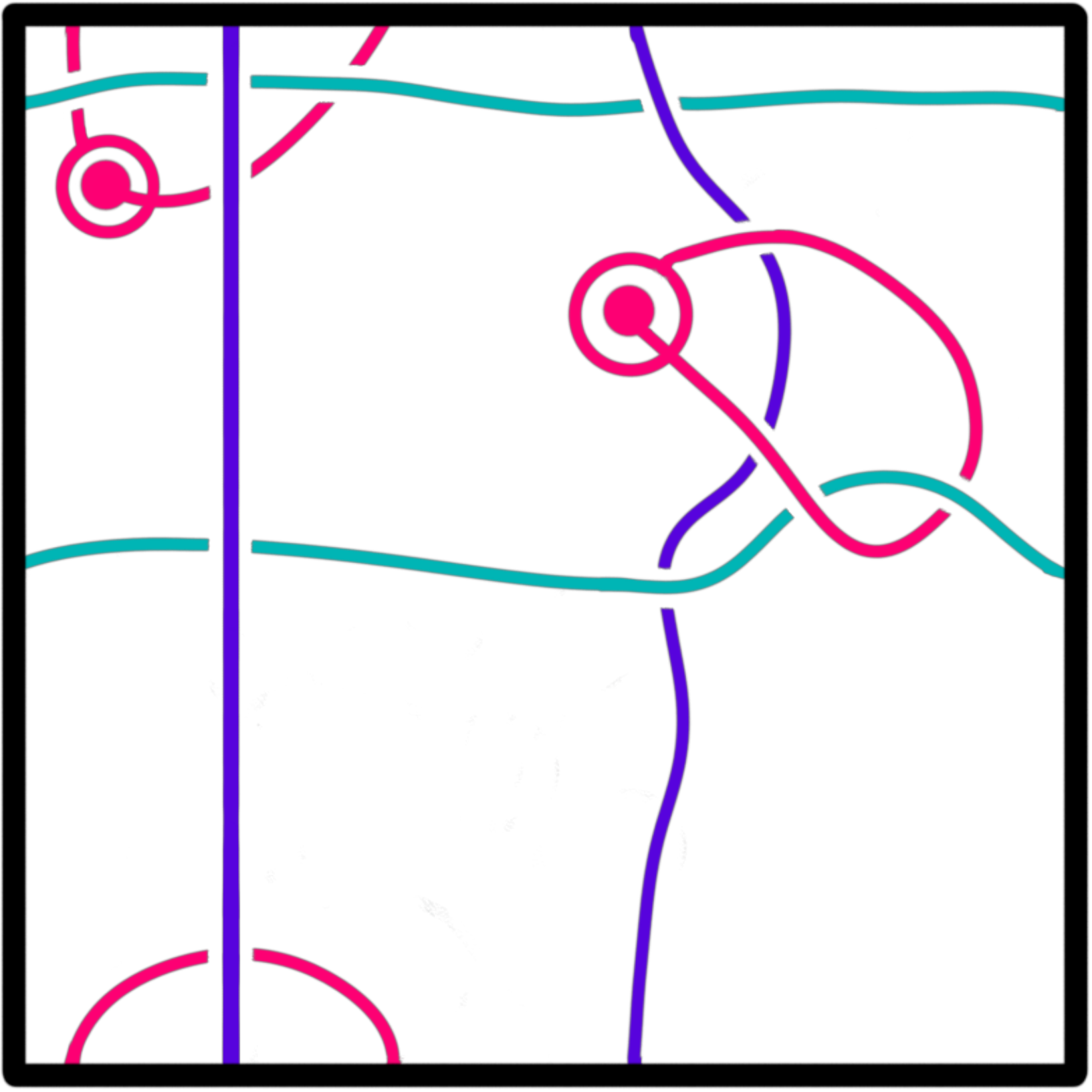}
        \caption{}
        \label{fig:1p6on2_to_the_6_dia_18}
    \end{subfigure}
    \begin{subfigure}[b]{0.19\textwidth}
        \centering
        \includegraphics[width=0.975\textwidth]{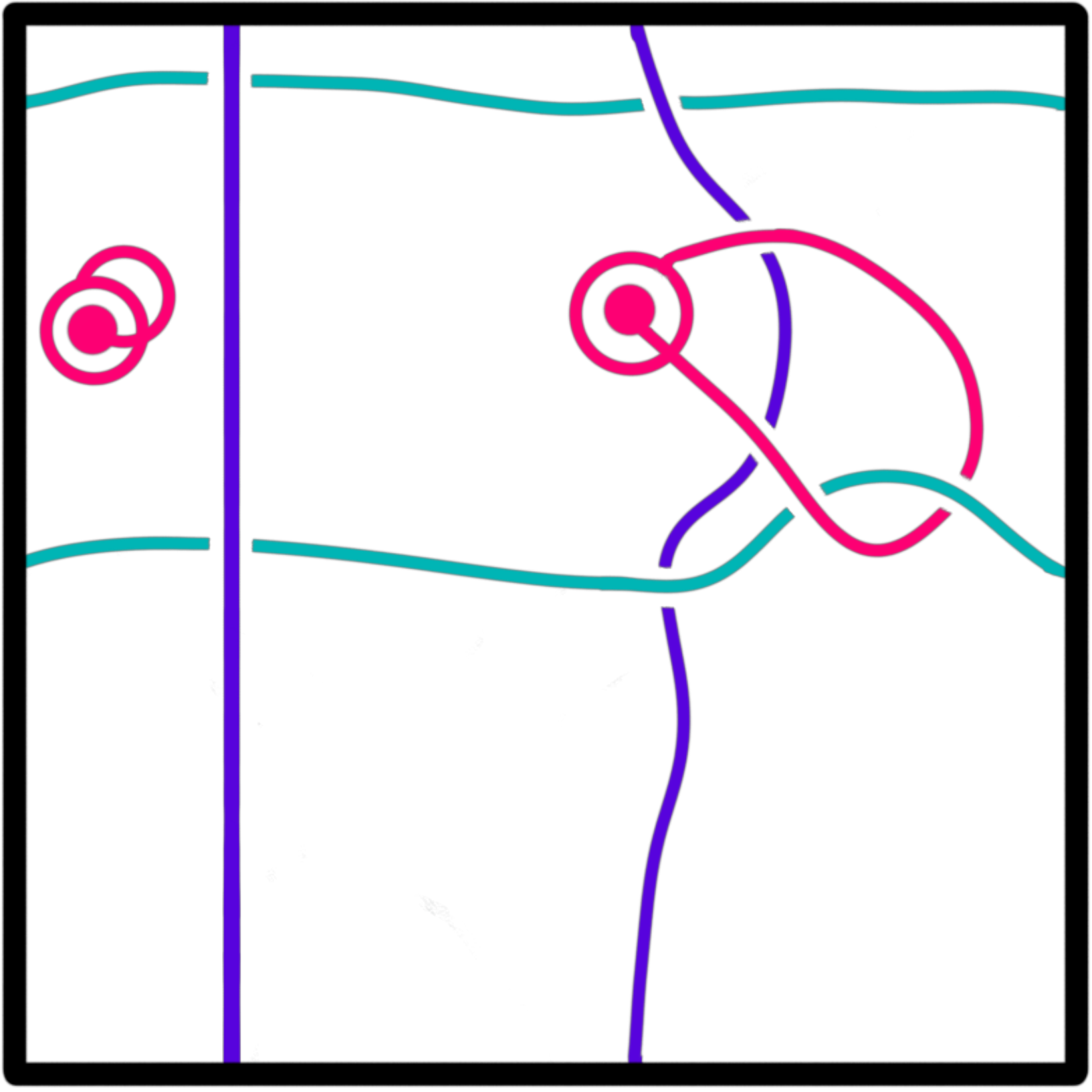}
        \caption{}
        \label{fig:1p6on2_to_the_6_dia_19}
    \end{subfigure}

    \vskip\baselineskip

    \begin{subfigure}[b]{0.19\textwidth}
        \centering
        \includegraphics[width=0.975\textwidth]{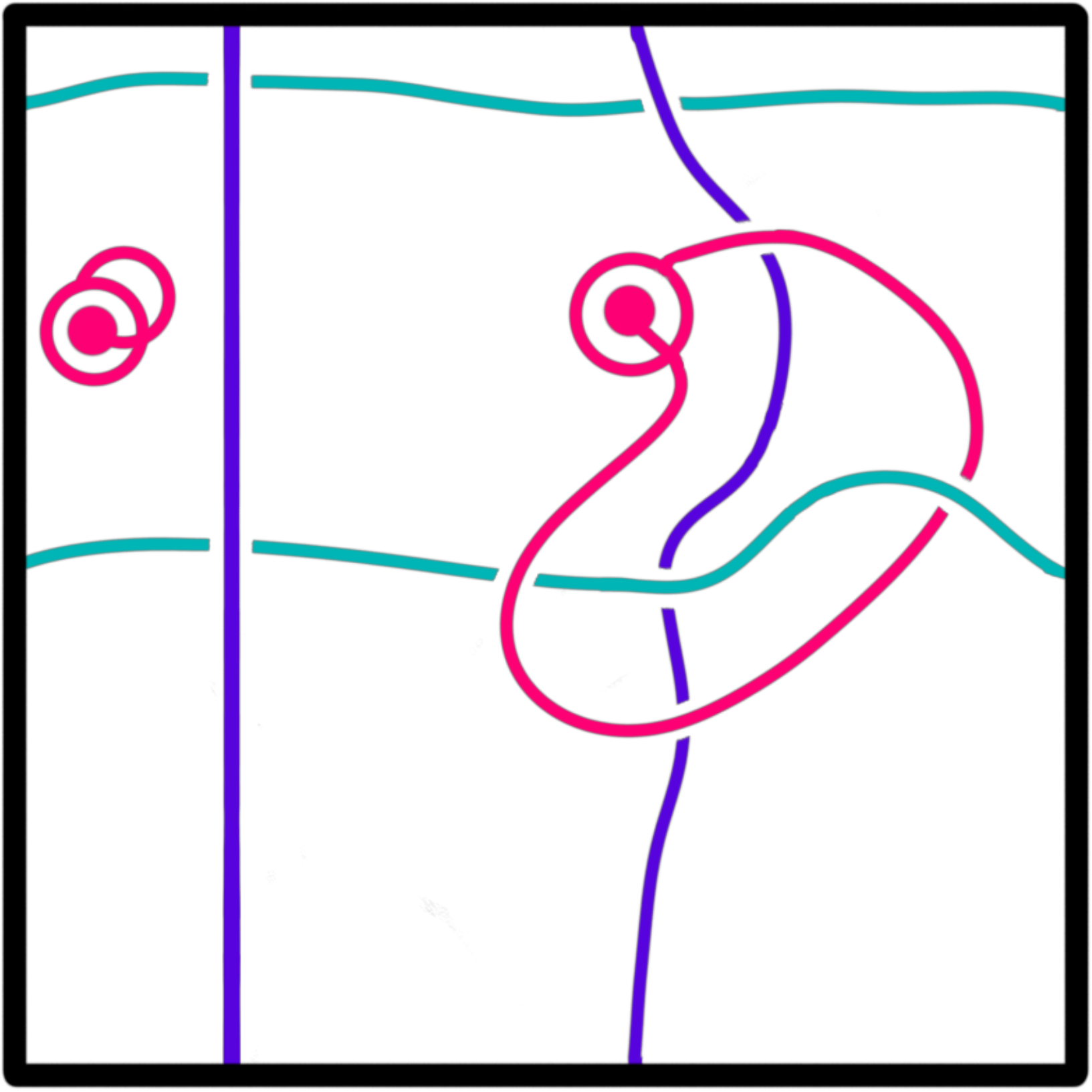}
        \caption{}
        \label{fig:1p6on2_to_the_6_dia_20}
    \end{subfigure}
    \begin{subfigure}[b]{0.19\textwidth}
        \centering
        \includegraphics[width=0.975\textwidth]{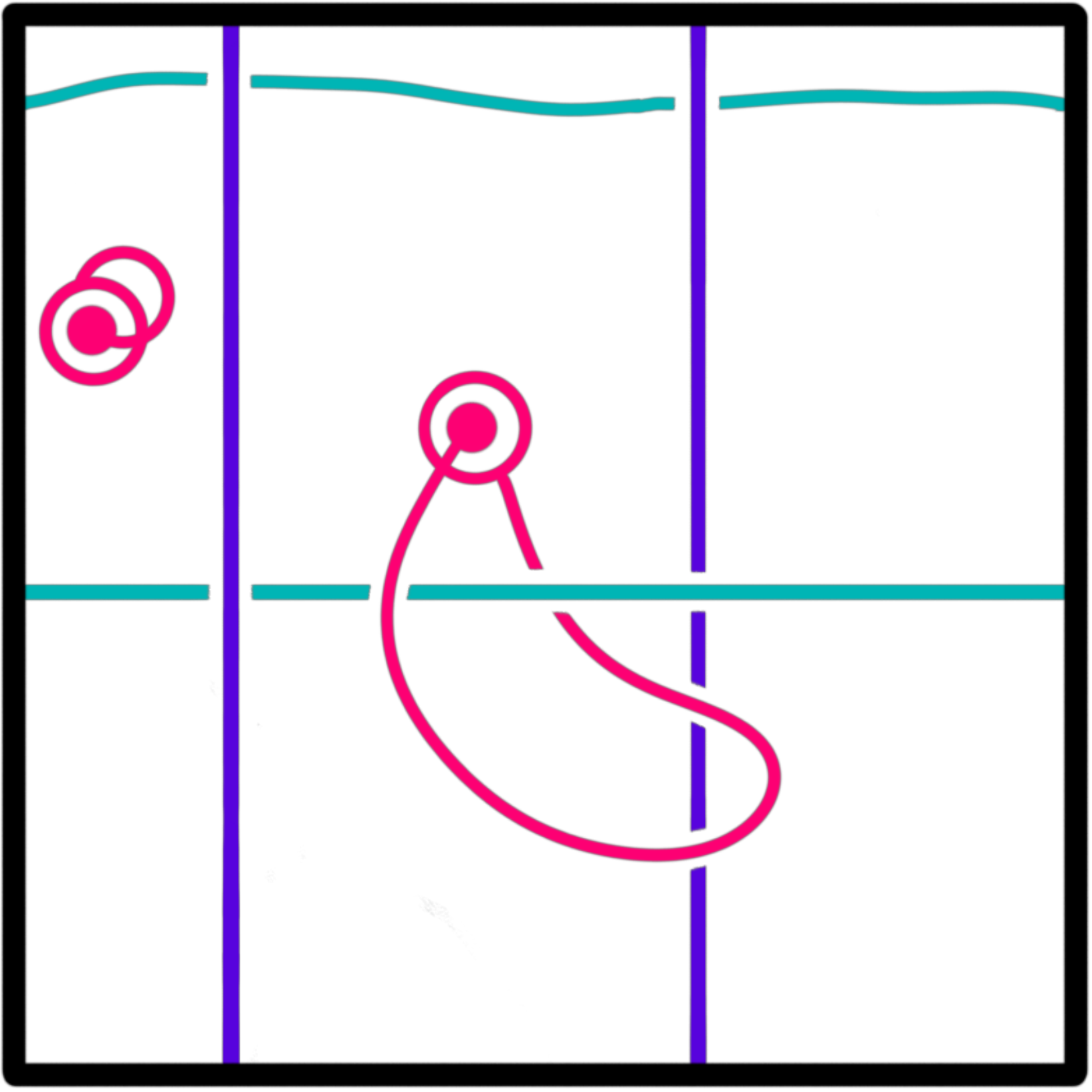}
        \caption{}
        \label{fig:1p6on2_to_the_6_dia_21}
    \end{subfigure}
    \begin{subfigure}[b]{0.19\textwidth}
        \centering
        \includegraphics[width=0.975\textwidth]{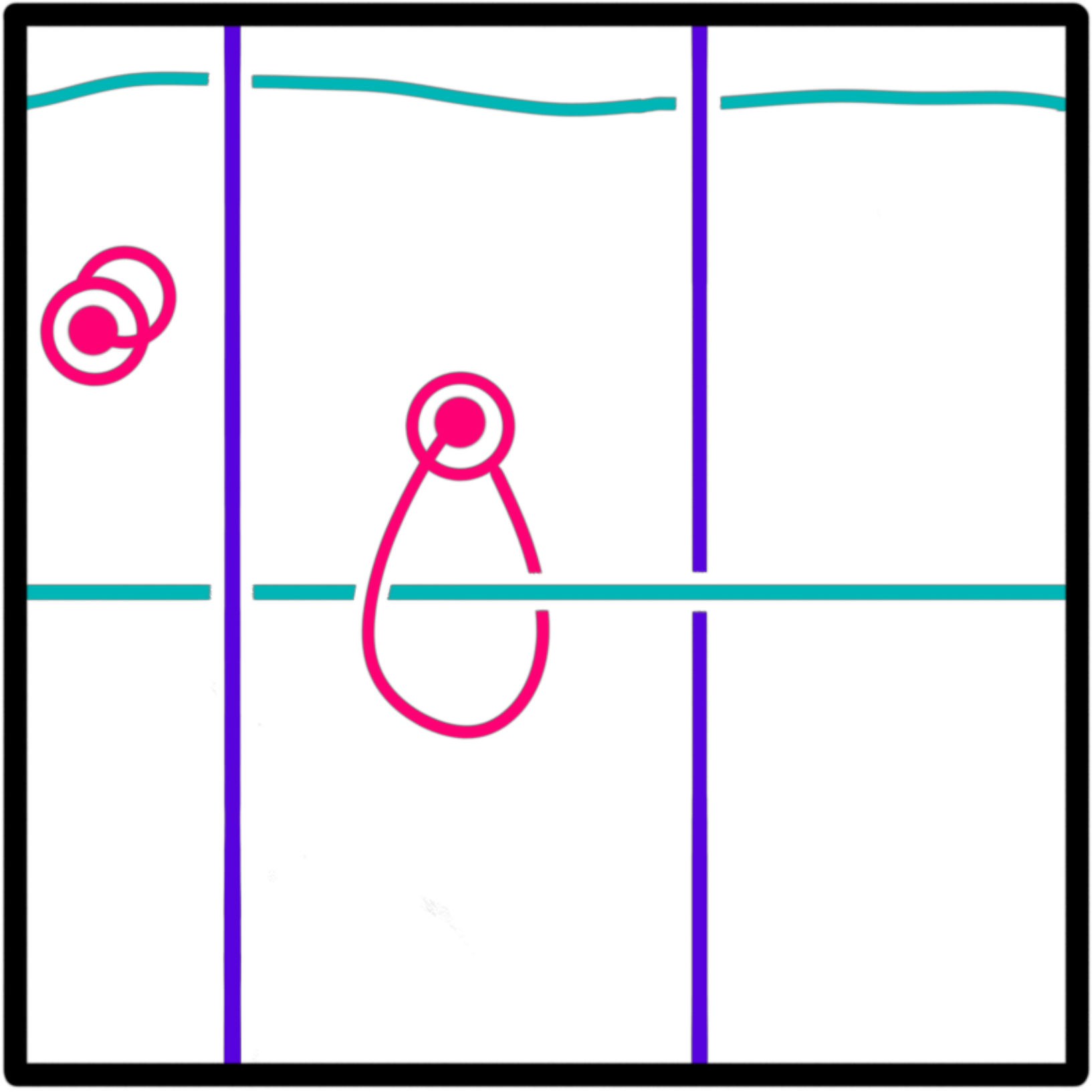}
        \caption{}
        \label{fig:1p6on2_to_the_6_dia_22}
    \end{subfigure}
    \begin{subfigure}[b]{0.19\textwidth}
        \centering
        \includegraphics[width=0.975\textwidth]{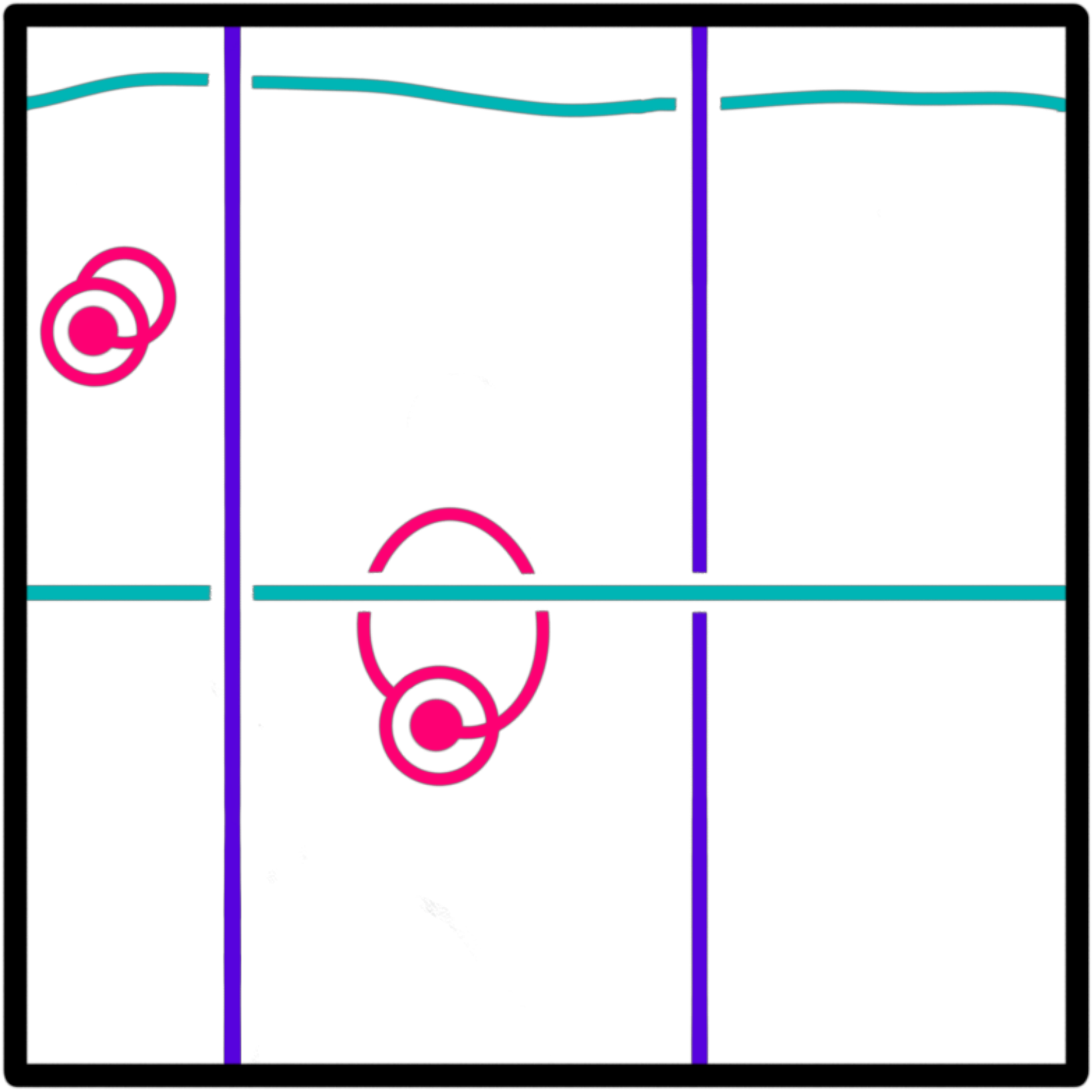}
        \caption{}
        \label{fig:1p6on2_to_the_6_dia_23}
    \end{subfigure}
    \begin{subfigure}[b]{0.19\textwidth}
        \centering
        \includegraphics[width=0.975\textwidth]{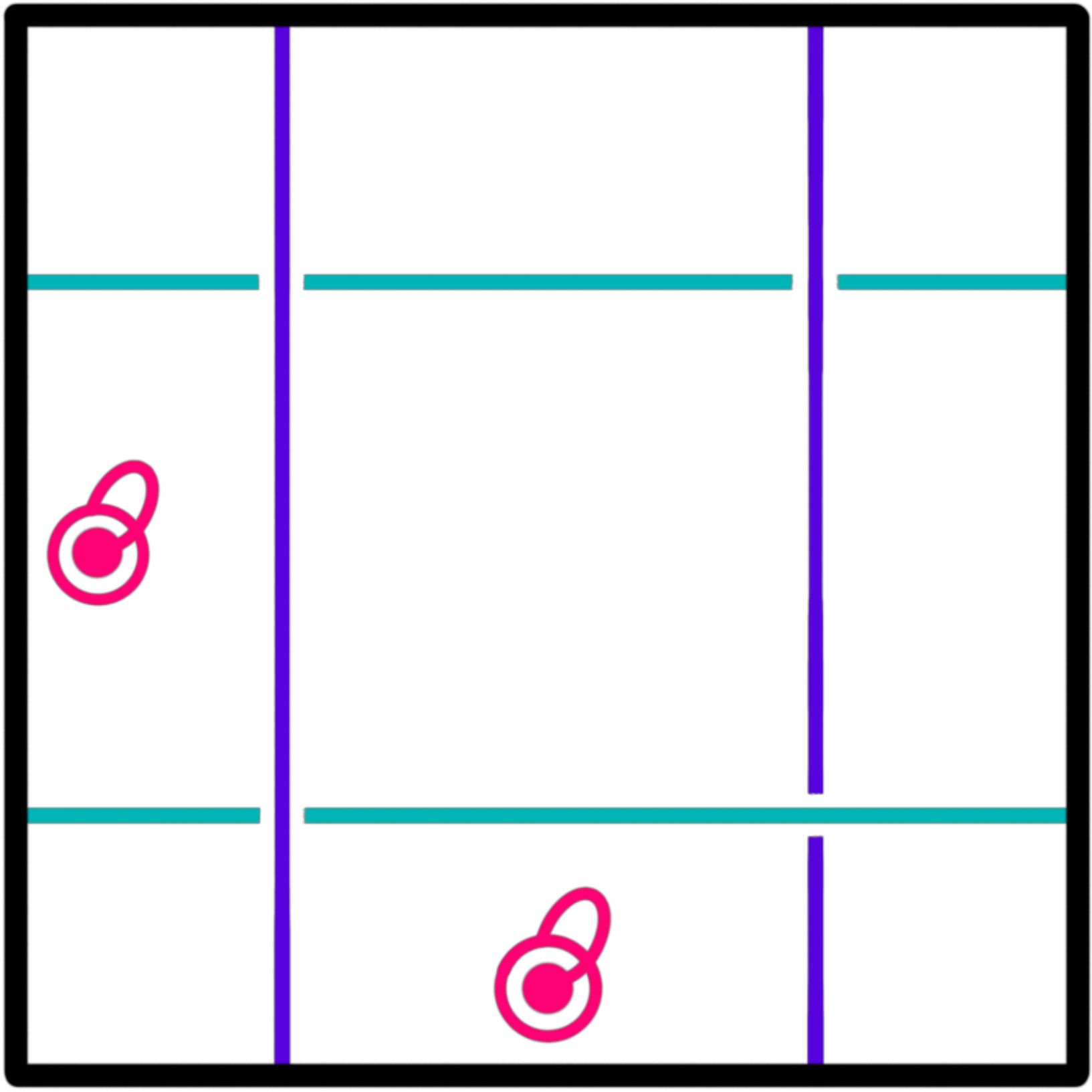}
        \caption{}
        \label{fig:1p6on2_to_the_6_dia_24}
    \end{subfigure}
    
    \caption{Computing an upper bound for the untangling number of the structure shown in Fig.~\ref{fig:1p6on2_to_the_6_extended}. The initial diagram of (b) is obtained from the unit cell displayed in (a). From (b) to (y), each diagram is obtained from the one before either by applying one crossing change on the crossing highlighted by a black circle or by applying some $R$-moves. The final diagram shown in (y) represents the $\Pi^{+}$ rod packing. A total of six crossings changes were applied, which means that the untangling number cannot be more than 6.}
    \label{fig:1p6on2_to_the_6_untangling}
\end{figure}

We end this section with the computation of an upper bound for the untangling number of the 3-periodic tangle shown in Fig.~\ref{fig:1p6on2_to_the_6_extended}. One unit cell of the structure is displayed in Fig.~\ref{fig:1p6on2_to_the_6_uc}. From Fig.~\ref{fig:1p6on2_to_the_6_dia_1} to Fig.~\ref{fig:1p6on2_to_the_6_dia_24}, a total of six crossing changes, applied on the crossings highlighted by black circles, and some additional $R$-moves are performed to reach a ground state. This means that untangling number is at most 6. The final diagram displayed in Fig.~\ref{fig:1p6on2_to_the_6_dia_24} represents the $\Pi^{+}$ rod packing.

\section{Computability of the untangling number}\label{sec:computability}
As already noted earlier, computing the untangling number is generally difficult, and in practice one can typically only obtain upper bounds. For instance, for the 3-periodic tangle considered in Fig.~\ref{fig:1p6on2_to_the_6_untangling}, it is difficult to ascertain that there is no faster way to reach a ground state than the one given in the computation. This situation is analogous to the computation of the unknotting number, where additional techniques have been developed over the past decades to establish lower bounds \cite{Murasugi1965,nakanishi_I,ma-qiu}. The present work is primarily definitional and does not provide such criteria, however, it lays the groundwork necessary for developing them in future research. An example of an interesting direction to explore would be the definition of an analogue to the \emph{Ma-Qiu index} given in \cite{ma-qiu}. It is obtained from a modified Wirtinger presentation of the fundamental group of the complement of a knot. Indeed, there is interest in the study of complements of unit cells of 3-periodic tangles in the 3-torus, and some work on the subject has already been carried out in \cite{purcell2024,purcell2025}. To obtain such a lower bound, a presentation of the fundamental group of complements of unit cells must be developed. Such a presentation is traditionally derived from a diagram, which is now available for 3-periodic tangles \cite{ANDRIAMANALINA2025109346}, but obtaining the presentation is not straightforward in this setting. Indeed, due to the periodic nature of the unit cells, there is no privileged position (traditionally above the curve components for classical knots) at which to place the base point of the loops used to compute the fundamental group of the unit cell. More generally, many lower bounds for the unknotting number are obtained from diagrammatic knot invariants. Such invariants are usually computed by assigning orientations to the curves and signs to crossings. However, extending such invariants to 3-periodic tangles remains challenging because of the $R_4$ move, shown in Fig.~\ref{fig:Reid_move_IV}, which reverses the type of a crossing.

While the development of rigorous lower bounds remains an important direction for future research from a theoretical perspective, we argue that empirical upper bounds are often sufficient for most simple cases, from a practical standpoint. For example, the computation in Fig.~\ref{fig:untangling_pi_plus_4,4,5_1_usual_first_time} likely coincides with the actual untangling number of the structure. Indeed, one can observe that the components of the unit cell shown in Fig.~\ref{fig:pi_plus_u2_uc} can be grouped into two independent sets of three components each. Each group is tangled, and the crossing changes shown in Fig.~\ref{fig:untangling_pi_plus_4,4,5_1_usual_1} and Fig.~\ref{fig:untangling_pi_plus_4,4,5_1_usual_4_1st} are necessary to respectively untangle them. This suggests that the untangling number of the structure cannot be less than 2. The untangling numbers of many simple examples can often be empirically determined in the same manner. For more complex cases, such as that in Fig.~\ref{fig:1p6on2_to_the_6_untangling}, computational approaches like those developed in \cite{Applebaum18082025} may assist in determining the untangling number. However, further concepts, such as the computational representation of a crossing diagram, must first be developed to enable this.

\section{Conclusion}\label{sec:conclusion}
This article describes the untangling number of 3-periodic tangles, a measure of entanglement complexity, as well as the ground states, which are the least tangled 3-periodic tangles. These ideas use the knot-theoretic diagrammatic representations of 3-periodic entanglements introduced and rigorously investigated in \cite{ANDRIAMANALINA2025109346}, a setting in which the crossing number, necessary for the development of the methods outlined here, can be defined.

We established in Theorem~\ref{thm:homotopy_ground_states} a sufficient criterion for a 3-periodic tangle to be a ground state, based on the minimum number of intersections between homotopy classes of curves. As an application of this criterion, we showed that rod packings are typical ground states.

The constructions developed here provide the necessary foundations for future research on lower bounds for the untangling number, as well as computational approaches for determining its value.

\begin{credits}
\subsubsection{\ackname}
This work was funded by the Deutsche Forschungsgemeinschaft (DFG, German Research Foundation) - Project number 468308535, and partially supported by RIKEN iTHEMS and by the JSPS KAKENHI Grant-in-Aid for Early-Career Scientists 25K17246. We also acknowledge the education licence for Houdini/SideFX, which was used for visualisation.

\subsubsection{\discintname}
The authors have no competing interests to declare that are
relevant to the content of this article.

\end{credits}
%
%

\bibliographystyle{splncs04}
\bibliography{bibliography}

\end{document}